\documentclass[12pt]{article}

\usepackage{authblk}
\usepackage{amsmath,amsfonts,amssymb,lmodern,geometry,enumerate}
\usepackage[font=small,labelfont=bf]{caption}
\usepackage{tkz-euclide}
\usepackage{subcaption}
\usepackage[T1]{fontenc}
\usepackage[latin1]{inputenc}
\usepackage[english]{babel}
\usepackage{lmodern}
\usepackage{scalefnt}
\usepackage{dsfont}
\usepackage{stmaryrd}
\usepackage{color}
\usepackage{bm,bbm}
\usepackage{mathrsfs,url,color}
\usepackage{breakcites} 
\usepackage{textcase}
\usepackage{graphicx}
\usepackage{wrapfig}
\usepackage{flushend,cuted}
\usepackage{bm}
\usepackage{tabularx}
\usepackage{color}
\usepackage{indentfirst}
\usepackage{amssymb}
\usepackage{xparse}
\usepackage{tikz}
\usepackage{pgfplots}
\pgfplotsset{compat=newest}
\usepackage{mdwlist}
\usepackage{amsmath}
\usepackage{amsthm}
\usepackage{nameref}
\usetikzlibrary[intersections,
positioning,
petri,
backgrounds,
fit,
decorations.pathmorphing,
arrows,
arrows.meta,
bending,
calc,
intersections,
through,
backgrounds,
shapes.geometric,
quotes,
matrix,
trees,
shapes.symbols,
graphs,
math,
patterns,
external,
scopes,
matrix,
lindenmayersystems,
shapes.callouts,
shapes.misc,
angles,
shapes.arrows,
shadings]

\newtheorem{thm}{Theorem}[section]

\newtheorem{defi}[thm]{Definition}
\newtheorem{lma}[thm]{Lemma}
\newtheorem{cor}[thm]{Corollary}
\newtheorem{ap}{Application}

\theoremstyle{remark}
\newtheorem{re}{Remark}[section]
\newcommand{\law}{\mathscr{L}}

\newcommand{\Pro}{\mathbb{P}} 
\newcommand{\prob}{\Pro}
\newcommand{\E}{\mathbb{E}}
\newcommand{\mean}{\E}
\newcommand{\real}{\mathbb{R}}

\newcommand{\Z}{\mathbb{Z}}
\newcommand{\C}{\mathbb{C}}

\newcommand{\indep}{\perp \!\!\! \perp}

\newcommand{\bone}{{\bf 1}}

\newcommand{\cn}{{\cal N}}

\def\d{{\delta}}
\def\e{{\epsilon}}
\def\iid{{i.i.d.}}

\def\Ceka{\v Cekan\-avi\v cius}

\def\equald{\stackrel{\mbox{\scriptsize{d}}}{=}}

\def\im{{\mbox{\sl i}}}
\def\ime{{\mbox{\sl\scriptsize i}}}

\def\scrB{{\mathscr{B}}}

\def\scrF{{\mathscr{F}}}

\def\[{\left[}
\def\]{\right]}
\def\({\left(}
\def\){\right)}

\newcommand{\cov}{{\rm Cov}}
\def\var{{\rm Var}}

\def\ignore#1{}
\def\Ref#1{(\ref{#1})}
\def\qed{\hfill\hbox{${\vcenter{\vbox{
					\hrule height 0.4pt\hbox{\vrule width 0.4pt height 6pt
						\kern5pt\vrule width 0.4pt}\hrule height 0.4pt}}}$}}

\newcounter{con}
\stepcounter{con}
\newcommand{\qcon}[1]{\addtocounter{con}{1}}

\newcounter{cproofa}

\newcounter{cproofb}

\newcounter{cproofc}

\numberwithin{equation}{section}

\ignore{\makeatletter
	\renewcommand\section{\@startsection {section}{1}{\z@}%
		{-3.5ex \@plus -1ex \@minus -.2ex}%
		{1.3ex \@plus.2ex}%
		{\center\small\sc\MakeTextUppercase}}
	\def\subsection#1{\@startsection {subsection}{2}{0pt}%
		{-3.5ex \@plus -1ex \@minus -.2ex}%
		{1ex \@plus.2ex}%
		{\bf\mathversion{bold}}{#1}}
	\def\subsubsection#1{\@startsection{subsubsection}{3}{0pt}%
		{\medskipamount}%
		{-10pt}%
		{\normalsize\itshape}{\kern-2.2ex. #1.}}
	\makeatother}
\allowdisplaybreaks[4]

\begin{document}

\title{\sc\bf\large\MakeUppercase{
			Normal approximation in total variation for statistics in geometric probability
}}

\author[1]{Tianshu Cong\thanks{{\sf{email: tcong1@student.unimelb.edu.au.}} Work supported by a Research Training Program Scholarship and a Xing Lei Cross-Disciplinary PhD Scholarship in Mathematics and Statistics at the University of Melbourne..}}
\author[1]{Aihua Xia\thanks{{\sf{email: aihuaxia@unimelb.edu.au}}. Work supported by the Australian Research Council Grants Nos DP150101459 and DP190100613.}}

\affil[1]{%
School of Mathematics and Statistics, the University of Melbourne, Parkville VIC 3010, Australia} 

\date{\today}
	
\maketitle
\vskip-1cm
\begin{abstract}
		We use Stein's method to establish the rates of normal approximation in terms of the total variation distance for a large class of sums of score functions of marked Poisson point processes on $\mathbb{R}^d$. As in the study under the weaker Kolmogorov distance, the score functions are assumed to satisfy stabilizing and
		moment conditions. At the cost of an additional non-singularity condition for score functions, we show that the rates are in line with those under the Kolmogorov distance. We demonstrate the use of the theorems in four applications: Voronoi tessellation, $k$-nearest neighbours, timber volume and maximal layers.
\end{abstract}

	
	
\section{Introduction}
	
Limit theorems of functionals of Poisson point processes initiated in \cite{AB93} have been of considerable interest in the literature, see, e.g., \cite{S12,S16,LSY19} and references therein. The key element leading to the success is the stabilization introduced in~\cite{PY01,PY05}. The main character of the stabilization is that insertion of a point into a Poisson point process only induces a local effect in some sense hence there is little change in the functionals. However, adding an additional point to the Poisson point process results in the Palm process of the Poisson point process at the point \cite[Chapter~10]{Kallenberg83} and it is shown in \cite{CX04,CRX20} that the magnitude of the difference between a point process and its Palm processes is directly linked to the accuracy of Poisson and normal approximations of the point process. This is also the fundamental reason why the limit theorems in the above mentioned papers can be established.
	
The normal approximation theory is generally quantified in terms of the Kolmogorov distance $d_K$: for two random variables $X_1$ and $X_2$ with distributions $F_1$ and $F_2$,
$$d_K(X_1,X_2):=d_K(F_1,F_2):=\sup_{x\in\real}|F_1(x)-F_2(x)|.$$ 
The well-known Berry-Esseen Theorem \cite{Berry41,Esseen42} states that if $X_i$, $1\le i\le n,$ are independent and identically distributed ($\iid$) random variables with mean 0 and variance 1, 
define $Y_n=\frac{\sum_{i=1}^nX_i}{\sqrt{n}}$, $Z\sim N(0,1)$, where $\sim$ denotes ``is distributed as'', then 
$$d_K(Y_n,Z)\le \frac{C\mean|{X_1}|^3}{\sqrt{n}}.$$
The Kolmogorov distance $d_K(F_1,F_2)$ measures the maximum difference between the distribution functions $F_1$ and $F_2$, but it does not
tell much about the difference between the probabilities $\prob(X_1\in A)$ and $\prob(X_2\in A)$ for a non-interval Borel set $A\subset\real$, e.g., $A=\cup_{i\in \Z}(2i,2i+0.5]$, where $\Z$ denotes the set of all integers.  Such difference is reflected in the total variation distance $d_{TV}(F_1,F_2)$ defined by
$$d_{TV}(X_1,X_2):=d_{TV}(F_1,F_2):=\sup_{A\in\scrB(\real)}|F_1(A)-F_2(A)|,$$
where $\scrB(\real)$ stands for the Borel $\sigma$-algebra on $\real$. If $F_i$'s are absolutely continuous, that is, for arbitrary $A$ in  $\scrB(\real)$,   $F_i(A):=\int_AF_i'(x)dx$, then the definition is equivalent to 
$$d_{TV}(F_1,F_2)=\frac12\sup_f\left|\int f(x)F_1'(x)dx-\int f(x) F_2'(x)dx\right|,$$
where the supremum is taken over all measurable functions $f$ on $(\real,\scrB(\real))$ such that $\|f\|:=\sup_{x\in \real} |f(x)|\le1$.
	
	Although central limit theorems in the total variation have been studied in some special circumstances (see, e.g., \cite{DF87,MM07,BC16}), it is generally believed that the total variation distance is too strong for normal approximation, see, e.g., \cite{Ceka00,CL10,Fang14}. For example, the total variation distance between any discrete distribution and any normal distribution is always 1. To recover central limit theorems in the total variation, a common approach is to discretize the distribution of interest and approximate it with a simple discrete distribution, e.g.,
	translated Poisson \cite{Rollin05,Rollin07}, centered binomial \cite{Rollin08}, discretized normal \cite{CL10,Fang14} and a family of polynomial type distributions \cite{GX06}. The multivariate versions of these approximations are investigated by~\cite{BLX18}.
	
	By discretizing a distribution $F$ of interest, we essentially group the probability of an area and put it at one point in the area, hence the information of $F(A)$ for a general set $A\in\scrB(\real)$ is completely lost. 
	In this paper, we consider the normal approximation in the total variation to the sum of random variables under various circumstances.
	
	{\bf An inspiring example:} \cite[p.~146]{Feller71}.
	Let $\{X_i: \ i\ge 1\}$'s be $\iid$ random variables taking values $0$ and $1$ with equal probability, then
	$X=\sum_{k=1}^\infty 2^{-k}X_k$ has uniform distribution on $(0,1)$.
	If we separate the even and odd terms into $U=\sum_{k=1}^\infty 2^{-2k}X_{2k}$ and $V=\sum_{k=1}^\infty 2^{-(2k-1)}X_{2k-1}$, then $U$ and $V$ are independent, {$2U\equald V$}, but both $U$ and $V$ have singular distributions.
	Now, we can construct mutually independent random variables $\{U_i,\ V_i: \, i\ge 1\}$ such that  $U_i\equald U-\mean U$ and $V_i\equald V-\mean V$.
	Consider $\xi_1=U_1+V_1$, $\xi_2=-V_1-U_2$, $\xi_3=U_2+V_2$, $\dots$, then $\{\xi_i\}$ is a sequence of 1-dependent and identically distributed random variables  having the uniform distribution on $(-0.5,0.5)$. One can easily verify that $\sum_{i=1}^n\xi_i$ does not converge to normal as $n\to\infty$, hence stronger conditions are needed to ensure normal approximation for the sum of dependent random variables. 
	
	Under the Kolmogorov distance, user-friendly conditions are usually formulated to ensure that the variance of the sum becomes large as $n\to\infty$. In the context of functionals of Poisson point processes, a typical condition to guarantee the variance of the sum converging to infinity is to assume nondegeneracy \cite{PY01,XY15}, that is, the conditional variance of the sum given the information outside a local region is away from $0$. Under the total variation distance, we use a non-singular condition instead of the nondegeneracy to ensure that the distribution of the functional is diffuse enough for a proper normal approximation for any Borel sets. This condition is almost necessary because it is an essential ingredient in the special case of the sum of $\iid$ random variables, see \cite{BC16}\footnote{We thank Vlad Bally for bringing their work to our attention.} for a brief review of the development for the CLT in total variation distance. 
	
	The Lebesgue decomposition theorem \cite[p.~134]{Halmos74} ensures that any distribution function $F$ on $\real$ can be represented as
	\begin{equation}F=(1-\alpha_F) F_s+\alpha_FF_a,\label{decom1}\end{equation}
	where $\alpha_F\in[0,1]$, $F_s$ and $F_a$ are two distribution functions such that, with respect to the Lebesgue measure on $\real$, $F_a$ is absolutely continuous and $F_s$ is singular \cite[p.~126]{Halmos74}.
	
	\begin{defi} A distribution function $F$ on $\real$ is said to be non-singular if $\alpha_F>0$. A random variable is said to be non-singular if its distribution function is non-singular.
	\end{defi}
	
	Recalling that two measures on the same measurable space $\mu_1\le\mu_2$ if $\mu_1(A)\le\mu_2(A)$ for all measurable sets $A$. We can see that a random variable $X$ is non-singular if and only if there exists a sub-probability measure $\mu\neq 0$ such that $\mu \le \law(X)$ and there exists a function $f$ on $\mathbb{R}$ satisfying that
		\begin{equation*}\mu(A)=\int_A f(x)dx,\mbox{ for all }A\in\scrB(\real). \end{equation*}
	
	In this paper, we demonstrate that many of the limit theorems of functionals of Poisson point processes with respect to the Kolmogorov distance in the literature, e.g., \cite{PY01,PY05,S12,S16}, still hold under the total variation distance. In Section~\ref{Generalresults}, we give definitions of the concepts, state the conditions and present the main theorems. In Section~\ref{Applications}, these theorems are applied to establish error bounds of normal approximation for statistics in Voronoi tessellation, $k$-nearest neighbours, timber volume and maximal layers. The proofs of the main results in Section~\ref{Generalresults} rely on a number of preliminaries and lemmas which are given in Section~\ref{Preliminaries}. For the ease of reading, all proofs are postponed to Section~\ref{Theproofs}.

	\section{General results}\label{Generalresults}
	
	We consider the functionals of a marked point process with a Poisson point process in $\mathbb{R}^d$ as its ground process and each point carries a mark in a measurable space $(T,\mathscr{T})$ independently of other marks, where $\mathscr{T}$ is a $\sigma$-algebra on $T$.
	More precisely,  let $\bm{S}:=\mathbb{R}^d\times T$ be equipped with the product $\sigma$-field $\mathscr{S}:=\mathscr{B}(\mathbb{R}^d)\times \mathscr{T}$, where $\mathscr{B}(\mathbb{R}^d)$ is the Borel $\sigma-$algebra of $\mathbb{R}^d$.  
	We use~$\bm{C}_{\bm{S}}$ to denote the space of all locally finite non-negative integer valued measures $\xi$, often called a {\it configuration}, on~$\bm{S}$ such 
	that $\xi(\{{x}\}\times T)\le 1$ for all ${x}\in\mathbb{R}^d$. The space~$\bm{C}_{\bm{S}}$ is endowed with the $\sigma$-field 
	$\mathscr{C}_{\bm{S}}$ generated by the vague topology \cite[p.~169]{Kallenberg83}. 
	A {\it marked point process\/} $\Xi$ is a measurable mapping from
	$(\Omega,\mathscr{F},\prob)$ to $(\bm{C}_{\bm{S}},\mathscr{C}_{\bm{S}})$ \cite[p.~49]{Kallenberg17}. The induced simple point process 
	$\bar\Xi(\cdot):=\Xi(\cdot\times T)$ is called the {\it ground process\/} \cite[p.~3]{Daley08} or 
	projection  \cite[p.~17]{Kallenberg17} of 
	the marked point process~$\Xi$ {on $\bm{S}$}. The functionals we study in the paper are defined on $\Gamma_\alpha:=\[-\frac{1}{2}\alpha^{\frac{1}{d}}, \frac{1}{2}\alpha^{\frac{1}{d}}\]^d$ having the forms
	$$W_\alpha:={\sum_{(x,m)\in\Xi_{\Gamma_\alpha}}\eta(\left(x,m\right), \Xi)}$$
	and
	$$\bar{W}_\alpha:={\sum_{{(x,m)\in\Xi_{\Gamma_\alpha}}}\eta(\left(x,m\right), \Xi_{\Gamma_\alpha},\Gamma_\alpha)=\sum_{{(x,m)\in\Xi_{\Gamma_\alpha}}}\eta(\left(x,m\right), \Xi,\Gamma_\alpha)},$$
	where $\Xi\sim\mathscr{P}_{\lambda,\mathscr{L}_T}$ is a marked Poisson point process having a homogeneous Poisson point process on $\mathbb{R}^d$ with intensity measure $\lambda dx$ as its ground process and $\iid$ marks on $(T,\mathscr{T})$ with the law $\mathscr{L}_T$,  $\Xi_{\Gamma_\alpha}$ is its restricted process {on $\Gamma_\alpha$} defined as $\Xi_A(B\times D):=\Xi((A\cap B)\times D)$ for all $D\in \mathscr{T}$ and $A,B\in \mathscr{B}(\mathbb{R}^d)$. The function $\eta$ is called a {\it score function} {(resp. {\it restricted score function})}, i.e., a measurable function on $\(\bm{S}\times \bm{C}_{\bm{S}},\mathscr{S} \times  \mathscr{C}_{\bm{S}}\)$ to $\(\mathbb{R},\mathscr{B}\(\mathbb{R}\)\)$ {(resp. a function mapping $\(\bm{S}\times \bm{C}_{\bm{S}}\times \mathbb{R}^d\)$ to $\mathbb{R}$ which is $(\mathscr{S}\cap (\Gamma_\alpha\times T) )\times  \mathscr{C}_{\bm{S}\cap (\Gamma_\alpha\times T)}\rightarrow \mathscr{B}(\mathbb{R})$ measurable for fixed the third coordinate)} and it represents the interaction between a point and the configuration. Because the interest is in the values of the score function of the points in a configuration, for convenience, $\eta\((x,m),\mathscr{X}\){~(\mbox{resp.}~\eta\((x,m),\mathscr{X},\Gamma_\alpha\))}$ is understood as $0$ for all $x\in \mathbb{R}^d$ and $\mathscr{X}\in \bm{C}_{\bm{S}}$ such that $(x,m)\notin \mathscr{X}$.
	We consider the score functions satisfying the following four conditions.
	
	
	

	
	\noindent{\it A2.1 Stabilization}
	
	For a locally finite configuration $\mathscr{X}$ and $z\in (\mathbb{R}^d\times T)\cup\{\emptyset\}$, write $\mathscr{X}^{\lbag z \rbag}=\mathscr{X}$ if $z=\emptyset$ and $\mathscr{X}^{\lbag z \rbag}=\mathscr{X}\cup\{z\}$ otherwise. We use $\delta_v$ denote the Dirac measure at $v$. The notion of stabilization is introduced in \cite{PY01} and we adapt it to our setup as follows. 
	\begin{defi}\label{defi4} (unrestricted case)
		A score function $\eta$ on $\mathbb{R}^d\times T$ is range-bound (resp. exponentially stabilizing, polynomially stabilizing of order $\beta>0$) with respect to intensity $\lambda$ and a probability measure $\mathscr{L}_T$ on $T$ if for all $x\in \mathbb{R}^d$, $z\in (\mathbb{R}^d\times T)\cup\{\emptyset\}$, and almost all realizations $\mathscr{X}$ of the homogeneous marked Poisson point process $\Xi\sim\mathscr{P}_{\lambda,\mathscr{L}_T}$, there exists an  $R:=R(x):=R(x,m_x,\mathscr{X}^{\lbag z \rbag})\in(0,\infty)$ (a radius of stabilization), such that for all locally finite $\mathscr{Y}\subset (\mathbb{R}^d\backslash B(x,R))\times T$, where $B(x,R)$ is the ball with centre $x$ and radius $R$, we have 
		$$
		\eta\left(\left(x,m_x\right), \left[\mathscr{X}^{\lbag z \rbag}\cap\left(B(x,R)\times T\right)\right]\cup \mathscr{Y} \right)=\eta\left(\left(x,m_x\right),\mathscr{X}^{\lbag z \rbag}\cap\left(B(x,R)\times T\right) \right)
		$$
		and the tail probability $$\tau(t):=\sup_{x\in \mathbb{R}^d,m_x\in \rm{supp}(\mathscr{L}_T)}\sup_{z\in (\mathbb{R}^d\times T)\cup\{\emptyset\}} \mathbb{P}\left(R(x,m_x,\Xi^{\lbag z \rbag}+\delta_{(x,m_x)})\ge t \right)$$ satisfies that
		$$
		\tau(t)=0\mbox{ for some }t\in \mathbb{R}_+	\ \ \ (\mbox{resp. }\tau(t)\le C_1e^{-C_2t},\ \tau(t)\le C_1t^{-\beta} {~\mbox{for all }t\in \mathbb{R}_+})
		$$
		for some positive constants $C_1$ and $C_2$. 
	\end{defi}
	
	For the functionals with input of restricted marked Poisson point process, we have the following counterpart of stabilization.
	Note that the score function for the restricted input is not affected by points outside $\Gamma_\alpha$.
		
	\begin{defi}\label{defi4r} (restricted case) We say the score function $\eta$ is {range-bound (resp. exponentially stabilizing, polynomially stabilizing of order $\beta>0$)} with respect to intensity $\lambda$ and a probability measure $\mathscr{L}_T$ on $T$ if for $\alpha\in \mathbb{R}_+$, $x\in \Gamma_\alpha$, and $z\in (\Gamma_\alpha \times T)\cup\{\emptyset\}$,  almost all realizations $\mathscr{X}$ of the homogeneous marked Poisson point process $\Xi\sim\mathscr{P}_{\lambda,\mathscr{L}_T}$, there exists a  $\bar{R}:=\bar{R}(x,\alpha):=\bar{R}(x,m_x,\alpha,\mathscr{X}^{\lbag z \rbag})\in(0,\infty)$ (a radius of stabilization), such that for all locally finite $\mathscr{Y}\subset (\Gamma_\alpha\backslash B(x,R))\times T$, we have 
		\begin{align}
			&\eta\left(\left(x,m_x\right), \left[\mathscr{X}_{\Gamma_\alpha}^{{\lbag z \rbag}}\cap\left(B(x,\bar{R})\times T\right)\right]\cup \mathscr{Y},\Gamma_\alpha \right)\nonumber\\
			=&\eta\left(\left(x,m_x\right),\mathscr{X}_{\Gamma_\alpha}^{{\lbag z \rbag}}\cap\left(B(x,\bar{R})\times T\right) ,\Gamma_\alpha\right)\label{defi4.1}
		\end{align}
		and the tail probability $$\bar{\tau}(t):=\sup_{x\in \mathbb{R}^d,m_x\in \rm{supp}(\mathscr{L}_T),\alpha\in \mathbb{R}_+}\sup_{z\in (\Gamma_\alpha\times T)\cup\{\emptyset\}} \mathbb{P}\left(\bar{R}(x,m_x,\alpha, \Xi^{{\lbag z \rbag}}+\delta_{(x,m_x)})\ge t \right)$$ satisfies that
		\begin{equation*}
			\bar{\tau}(t)=0\mbox{ for some }t\in \mathbb{R}_+	\ \ \ (\mbox{resp. }\bar{\tau}(t)\le C_1e^{-C_2t},\ \bar{\tau}(t)\le C_1t^{-\beta} {~\mbox{for all }t\in \mathbb{R}_+})
		\end{equation*}
		for some positive constants $C_1$ and $C_2$. 
	\end{defi}
		
	\noindent{\it A2.2 Translation Invariance}
	
	We write $d(x,A):=\inf\{d(x,y);~y\in A\}$, $A\pm B:=\{x\pm y; \ x\in A,\ y\in B\}$ for $x\in \mathbb{R}^d$ and  $A,B\in\mathscr{B}\(\mathbb{R}^d\)$ and define the shift operator as $\Xi^x(\cdot\times D):=\Xi((\cdot+x)\times D)$ for all $x\in \mathbb{R}^d$, $D\in \mathscr{T}$. 
	
	{\it A2.2.1 Unrestricted Case:} 
	
	\begin{defi}\label{invar} The score function $\eta$ is {\it translation invariant} if for all locally finite configuration $\mathscr{X}$ and $x,y\in \mathbb{R}^d$ and $m\in T$,
		{$\eta((x+y,m),\mathscr{X})=\eta((x,m),\mathscr{X}^{y})$.} 
	\end{defi}
	
	{\it A2.2.2 Restricted Case:} 
	
	As a translation may send a configuration to outside of $\Gamma_\alpha$ resulting in a completely different configuration inside $\Gamma_\alpha$, it is necessary to focus on the part that affects the score function, therefore, 
	we expect the score function to take the same value for two configurations if the parts within their stabilising radii are completely inside $\Gamma_\alpha$ and one is a translation of the other. More precisely, we have the following definition.
	
	\begin{defi}\label{traninvres0} A stabilizing score function is called {\it translation invariant} if for any $\alpha>0$, $x\in\Gamma_\alpha$  
		and $\mathscr{X}\in \bm{C}_{\mathbb{R}^d\times T}$ such that $\bar{R}(x,m,\alpha,\mathscr{X})\le d(x,\partial\Gamma_\alpha)$, {where $\partial A$ stands for the boundary of $A$}, then $\eta\(x,m,\mathscr{X},\Gamma_\alpha\)=\eta\(x',m,\mathscr{X}',\Gamma_{\alpha'}\)$ and $\bar{R}(x',m,\alpha',\mathscr{X}')=\bar{R}(x,m,\alpha,\mathscr{X})$ for all $\alpha'>0$, $x'\in\Gamma_{\alpha'}$ and $\mathscr{X}'\in \bm{C}_{\mathbb{R}^d\times T}$ such that $\bar{R}(x,m,\alpha,\mathscr{X})\le d(x',\partial\Gamma_{\alpha'})$ and $\(\mathscr{X}'_{B(x',\bar{R}(x,m,\alpha,\mathscr{X}) )}\)^{x'}=\(\mathscr{X}_{B(x,\bar{R}(x,m,\alpha,\mathscr{X}) )}\)^{x}$.
	\end{defi}
	
	Noting that there is a tacit assumption of consistency in Definition \ref{traninvres0}, which implies that if $\eta$ is translation invariant in Definition \ref{traninvres0}, there exists a $\bar{g}: \bm{C}_{\mathbb{R}^d\times T}\rightarrow \mathbb{R}$ such that $$\lim_{\alpha\rightarrow \infty}\eta\((0,m),\mathscr{X},\Gamma_\alpha\)=\bar{g}\(\mathscr{X}+\delta_{(0,m)}\)$$ for $\law_{T}$ almost sure $m\in T$ and almost all realizations $\mathscr{X}$ of the homogeneous marked Poisson point process $\Xi\sim\mathscr{P}_{\lambda,\mathscr{L}_T}$. Furthermore, we can see that for each score function $\eta$ satisfying the translation invariance in Definition \ref{traninvres0}, there exists a score function for the unrestricted case by setting $\bar{\eta}((x,m),\mathscr{X}):=\bar{g}(\mathscr{X}^{x})\mathbf{1}_{(x,m)\in \mathscr{X}}$ and writing the radii of stabilization in the sense of Definition~\ref{defi4} as $R$. From the construction, $\bar{\eta}$ is range bound (resp. exponentially stabilizing, polynomially stabilizing of order $\beta>0$) in the sense of Definition~\ref{defi4} if $\eta$ is range bound (resp. exponentially stabilizing, polynomially stabilizing of order $\beta>0$) in the sense of Definition~\ref{defi4r}. Moreover, if $B(x, R(x)) \subset \Gamma_\alpha$, then $\bar{R}(x, \alpha)= R(x)$  and if $B(x, R(x)) \not\subset \Gamma_\alpha$, then $\bar{R}(x,\alpha)>d(x,\partial \Gamma_\alpha)$, but there is no definite relationship between $\bar{R}$ and $R$.

	\noindent{\it A2.3 Moment condition} 
	
	{\it Unrestricted Case:} The score function $\eta$ is said to satisfy the $k$th moment condition if
	\begin{equation}\label{thm2.1}
		\mathbb{E}\left|\eta\left(({\bf 0},M_3),\Xi+a_1\delta_{(x_1,M_1)}+a_2\delta_{(x_2,M_2)}+\delta_{({\bf 0},M_3)}\right)\right|^k\le C
	\end{equation}
	for some positive constant $C$ for all $a_i\in \{0,1\}$, distinct $x_i\in \mathbb{R}^d$, $i\in\{1,2\}$, and  $\iid$ random elements $M_1$, $M_2$, $M_3$ that are independent of $\Xi$ following the distribution $\mathscr{L}_{T}$.

	{\it Restricted Case:} For restricted score functions, $\eta$ is said to satisfy the $k$th moment condition if there exists a positive constant $C$ independent of $\alpha$ such that
	\begin{equation}\label{thm2.1r}
		\mathbb{E}\left|\eta\left((x_3,M_3),\Xi_{\Gamma_\alpha}+a_1\delta_{(x_1,M_1)}+a_2\delta_{(x_2,M_2)}+\delta_{(x_3,M_3)}\right)\right|^k\le C
	\end{equation}
	for all $a_1,a_2\in \{0,1\}$, distinct $x_1,x_2,x_3\in \Gamma_\alpha$, and  $\iid$ random elements $M_1$, $M_2$, $M_3$ that are independent of $\Xi$ following the distribution $\mathscr{L}_{T}$. From the construction, if $\eta$ is stabilizing, then $\bar{\eta}$ satisfies the moment condition of the same order in the sense of \Ref{thm2.1}.
	
	\noindent{\it A2.4 Non-singularity} 
	
	{\it Unrestricted Case:} The score function is said to be {\it non-singular} if 
	\begin{equation}
		\law\left(\left.\sum_{(x,m)\in \Xi}\eta\((x, m),\Xi\)\mathbf{1}_{d(x, N_0)<R(x)}\right|\sigma(\Xi_{N_0^c})\right)\label{non-sin}\end{equation}  
	has a positive probability to be non-singular for some bounded set $N_0$. That is, the sum of the values of the score function that affected by the points in $N_0$ is non-singular. 
	
	{\it Restricted Case:} We define the non-singularity when the score function is stabilizing. The score function $\eta$ for restricted input is said to be  {\it non-singular} if it is stabilizing and the corresponding $\bar{\eta}$ satisfies that
	{	\begin{equation}
			\law\left(\left.\sum_{(x,m)\in \Xi}\bar{\eta}\((x, m),\Xi\)\mathbf{1}_{d(x, N_0)<R(x)}\right|\sigma(\Xi_{N_0^c})\right)
			\label{non-sinr}
	\end{equation}  }
	has a positive probability to be non-singular {for some bounded set $N_0$}.

	The main result for $W_\alpha$ (unrestricted case) is summarized below.
	
	\begin{thm}\label{thm2} Let $Z_\alpha\sim N(\mean W_{\alpha}, \var(W_{\alpha}))$. Assume that the score function $\eta$ is translation invariant in Definition~\ref{invar} and non-singular~\Ref{non-sin}.
		\begin{description} 
			\item{(i)} If $\eta$ is range-bound as in Definition~\ref{defi4} and satisfies the third moment condition~\Ref{thm2.1}, then
			$$d_{TV}(W_\alpha,Z_\alpha)\le O\(\alpha^{-\frac{1}{2}}\).$$
			\item{(ii)} If $\eta$ is exponentially stabilizing as in Definition~\ref{defi4} and satisfies the third moment condition~\Ref{thm2.1}, then
			$$d_{TV}(W_\alpha,Z_\alpha)\le O\(\alpha^{-\frac{1}{2}}\ln(\alpha)^{\frac{5d}{2}}\).$$ 
			\item{(iii)} If $\eta$ is polynomially stabilizing as in Definition~\ref{defi4} with parameter $\beta>\frac{(15k-14)d}{k-2}$ and satisfies the $k'$-th moment condition ~\Ref{thm2.1} with $k'>k\ge 3$, then 
			\begin{align*}
				d_{TV}\(W_\alpha,Z_\alpha\)&\le O\(\alpha^{-\frac{\beta(k-2)[\beta(k-2)-d(15k-14)]}{(k\beta-2\beta-dk)(5dk+2\beta k-4\beta)}}\).
			\end{align*}
		\end{description}
	\end{thm}
	
	When approximation error is measured in terms of the Kolmogorov distance, the distributions of $W_\alpha$ and $\bar{W}_\alpha$ are often close for large $\alpha$. However, in terms of the total variation distance, one can not infer the accuracy of normal approximation of $W_\alpha$ by taking limit of that for $\bar{W}_\alpha$. For this reason, we need to adapt the conditions accordingly and tackle $\bar{W}_\alpha$ separately.
	We state the main result for $\bar{W}_\alpha$ (restricted case) in the following theorem.
	
	\begin{thm}\label{thm2a} Let $\bar{Z}_\alpha\sim N(\mean\bar{W}_\alpha,\var(\bar{W}_\alpha))$. Assume that $\eta$ is translation invariant in Definition~\ref{traninvres0} and non-singular~\Ref{non-sinr}.
		\begin{description}
			\item{(i)} If $\eta$ is range-bound as in Definition~\ref{defi4r} and satisfies the third moment condition~\Ref{thm2.1r}, then
			$$d_{TV}(\bar{W}_\alpha,\bar{Z}_\alpha)\le O\(\alpha^{-\frac{1}{2}}\).$$
			\item{(ii)}	If $\eta$ is exponentially stabilizing as in Definition~\ref{defi4r} and satisfies the third moment condition~\Ref{thm2.1r}, then
			$$d_{TV}(\bar{W}_\alpha,\bar{Z}_\alpha)\le O\(\alpha^{-\frac{1}{2}}\ln(\alpha)^{\frac{5d}{2}}\).$$ 
			\item{(iii)}	If $\eta$ is polynomially stabilizing as in Definition~\ref{defi4r} with parameter $\beta>\frac{(15k-14)d}{k-2}$ and satisfies the $k'$-th moment condition~\Ref{thm2.1r} with $k'>k\ge 3$, then
			$$d_{TV}\(\bar{W}_\alpha,\bar{Z}_\alpha\)\le O\(\alpha^{-\frac{\beta(k-2)[\beta(k-2)-d(15k-14)]}{(k\beta-2\beta-dk)(5dk+2\beta k-4\beta)}}\).
			$$
		\end{description}
	\end{thm}
	
	\section{Applications}\label{Applications}
	{\rm~Our main result can be applied to a wide range of geometric probability problems, including {normal approximation of }functionals of $k$-nearest neighbors graph, Voronoi graph, sphere of influence graph, Delaunay triangulation, Gabriel graph and relative neighborhood graph. To keep our article in a reasonable size, we only show the $k$-nearest neighbors graph and the Voronoi graph in details. We can see that many functionals of the graphs such as total edge length satisfy the conditions of the main theorems naturally and the ideas for verifying these conditions are similar. For the ease of reading, we briefly introduce these graphs below, more details can be found in \cite{Devroye88, T82}. }
	
	{\rm Let $\mathscr{X}\subset {\mathbb{R}^d}$ be a locally finite point set:}
	
	\begin{description}
		
		\item{(i)} {\it $k$-nearest neighbors graph}~{\rm~The $k$-nearest neighbors graph $NG(\mathscr{X})$ is the graph obtained by including $\{x,y\}$ as an edge whenever $y$ is one of the $k$ points nearest to $x$ or $x$ is one of the $k$ points nearest to $y$. A variant of the $NG(\mathscr{X})$ which has been considered in the literature is the {\it directed} graph $NG'\(\mathscr{X}\)$, which is constructed by inserting a directed edge $(x,y)$ if $y$ is one of the $k$ nearest neighbors of $x$.}
		
		\item{(ii)} {\it Voronoi tessellation}~{\rm We enumerate the points in $\mathscr{X}$ as $\{x_1,~x_2,\dots\}$, denote the locus of points in $\mathbb{R}^d$ closer to $x_i$ than to any other points in $\mathscr{X}$ by $C(x_i):=C(x_i,\mathscr{X})$ for all $i \in \mathbb{N}$. We can see that $C(x_i)$ is the intersection of half-planes and when the point set $\mathscr{X}$ has $n<\infty$ points, $C(x_i)$'s is a convex polygonal region with at most $n-1$ sides for $i\le n$. The cells $C(x_i)$ form a partition of $\mathbb{R}^d$, the partition is called Voronoi tessellation and the points in $\mathscr{X}$ are usually called Voronoi generators. }
		
		\item{(iii)} {\it Delaunay triangulation}~{\rm~The Delaunay triangulation graph puts an edge between two points in $\mathscr{X}$ if these points are centers of adjacent Voronoi cells, which is a dual to the Voronoi tessellation. }
		
		\item{(iv)} {\it Gabriel graph}~{\rm~Gabriel graph puts an edge between two points $x$ and $y$ in $\mathscr{X}$ if the ball centered at the middle point $\frac{x+y}{2}$ with radius $\|\frac{x-y}{2}\|$ does not contain any other points in $\mathscr{X}$. We can see that Gabriel graph is a subgraph of Delaunay triangulation graph. }
		
		\item{(v)} {\it Relative neighborhood graph}~{\rm~Relative neighborhood graph puts an edge between two points $x$ and $y$ in $\mathscr{X}$ if $B(x,\|x-y\|)\cap B(y,\|x-y\|)\cap\mathscr{X} =\emptyset$, i.e., the loon between $x$ and $y$ does not contain any other points in $\mathscr{X}$. This graph is a subgraph of Gabriel graph, so is also a subgraph of Delaunay triangulation graph. }
		
		\item{(vi)} {\it Sphere of influence graph}~{\rm~The sphere of influence graph of a locally finite point set $\mathscr{X}\subset {\mathbb{R}^d}$ is the graph obtained by including $\{x,y\}$ as an edge whenever $x,y\in\mathscr{X}$, $\|x-y\|\le \|x-N(x,\mathscr{X})\|+\|y-N(y, \mathscr{X})\|$, where for $z\in\mathscr{X}$, $N(z,\mathscr{X})$ is the nearest point of $z$ in $\mathscr{X}$. That is, for every point $z\in \mathscr{X}$, we draw a circle with center $z$ and radius being the distance between $z$ and its nearest point in $\mathscr{X}$, then two points are connected if the circles centered at two points $x,y$ intersect. }
	\end{description}
	
	\subsection{The total edge length of $k$-nearest neighbors graph}\label{knear}
	
	\begin{thm}
		If $\Xi$ is a homogeneous Poisson point process, then the total edge length $\bar{W}_\alpha$ (resp. $\bar{W}_\alpha'$) of $NG\(\Xi_{\Gamma_{\alpha}}\)$ (resp. $NG'\(\Xi_{\Gamma_{\alpha}}\)$) satisfies 
		$$d_{TV}\(\bar{W}_\alpha, {\bar{Z}}_\alpha\)\le O\(\alpha^{-\frac{1}{2}}\ln(\alpha)^{\frac{5d}{2}}\)~\(\mbox{resp.~}d_{TV}\(\bar{W}_\alpha', {\bar{Z}}_\alpha'\)\le O\(\alpha^{-\frac{1}{2}}\ln(\alpha)^{\frac{5d}{2}}\)\),$$ where ${\bar{Z}}_\alpha$ (resp. ${\bar{Z}}_\alpha'$) is a normal random variable with the same mean and variance as those of $\bar{W}_\alpha$ (resp. $\bar{W}_\alpha'$).
	\end{thm}
	
	\noindent{\it Proof.~} We only show the claim for the total edge length of $NG\(\Xi_{\Gamma_{\alpha}}\)$ since $NG'\(\Xi_{\Gamma_{\alpha}}\)$ can be handled with the same idea. The score function in this case is $${\eta\(x,\mathscr{X},\Gamma_\alpha\):=\frac{1}{2}\sum_{y\in \mathscr{X}_{\Gamma_\alpha}}\|y-x\|\mathbf{1}_{\{(x,y)\in NG(\mathscr{X}_{\Gamma_\alpha})\}},}$$ 
	\begin{wrapfigure}{r}{6cm}
		\begin{tikzpicture}[scale=3]
			\filldraw (0,0) circle (1pt);
			\foreach \x in {15,75,135,195,255,315}
			\draw[] (0,0) -- (\x:1);
			\draw[] (-0.5657,0.5657) -- (135:1); 
			\node at (45:0.55) {$T_i(t)$};
			\draw[latex-latex] (135:0.045) -- (135:0.79);
			\draw (0,0) circle (0.8);
			\node[fill=white,inner sep=2pt] at (135:0.45) {{$t$}};
			\node[fill=white,inner sep=1pt] at (219:0.13) {$x$};
			\draw[dashed] (-1,0)--(1,0) node[right]{};
			\draw[dashed] (0,-1)--(0,1) node[above]{};
		\end{tikzpicture}
		\caption{$k$-nearest: stabilization}
		\label{kns}
		\vspace{-15pt}
	\end{wrapfigure}
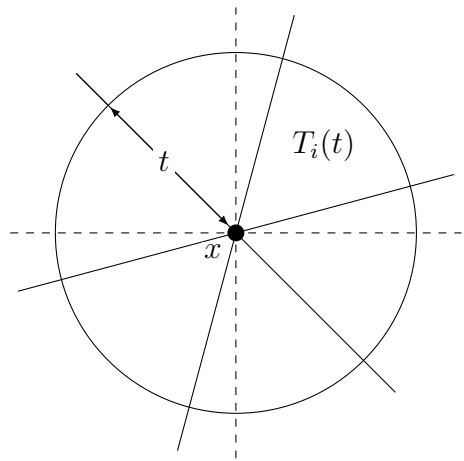
	which is clearly translation invariant.
	To apply Theorem~\ref{thm2a}, we need to check the moment condition \Ref{thm2.1r}, non-singularity \Ref{non-sinr} and stabilizing condition as in Definition~\ref{defi4r}. For simplicity, we show these conditions in two dimensional case and the argument can be easily extended to $\mathbb{R}^d$.
	
	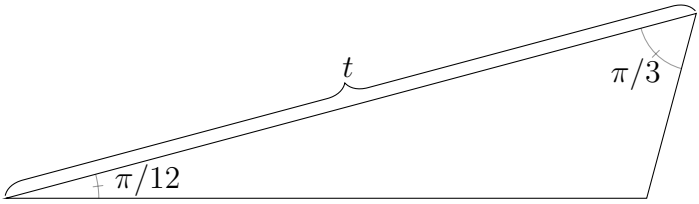
\begin{wrapfigure}{l}{9.5cm}
		\begin{tikzpicture}[scale=1.9]
			\coordinate (O) at (0,0);
			\coordinate (A) at (4.4825,0);
			\coordinate (B) at (4.8293,1.2940);
			\draw (O)--(A)--(B)--cycle;
			\draw [decorate, decoration={brace,amplitude=8pt,raise=1pt}] (0,0) -- (4.8293,1.2940) node [midway, anchor=south west, yshift=-1mm, xshift=-6mm, outer sep=10pt,font=] {$t$};
			\tkzMarkAngle[fill= orange,size=0.65cm,%
			opacity=.4](A,O,B)
			\tkzLabelAngle[pos = 1.0](A,O,B){$\pi/12$}
			\tkzMarkAngle[fill= orange,size=0.4cm,%
			opacity=.4](O,B,A)
			\tkzLabelAngle[pos = 0.6](O,B,A){$\pi/3$}
		\end{tikzpicture}
		\caption{$k$-nearest: $A_t$}
		\label{kntriangle}
		\vspace{-15pt}
	\end{wrapfigure}We start with the exponential stabilization and fix $\alpha>0$ and $x\in \Gamma_\alpha$. Referring to Figure~\ref{kns}, for each $t>0$, we construct six disjoint sectors of the same size $T_j(t)$, $1\le j\le 6$, with $x$ as the centre and angle $\frac{\pi}{3}$. In consideration of edge effects near the boundary of $\Gamma_\alpha$, the sectors are rotated around $x$ such that all straight edges of the sectors have a minimal angle $\pi/12$ with respect to the edges of $\Gamma_\alpha$. It is clear that $T_j(t)\subset T_j(u)$ for all $0<t<u$. Set $T_{j}(\infty)=\cup_{t>0}T_j(t)$ for $1\le j\le 6$, then from the property of the Poisson point process, there are infinitely many points in $\Xi\cap T_{j}(\infty)$ for all $j$ a.s. Let $|A|$ denote the cardinality of the set $A$ and define
	\begin{equation}t_{x,\alpha}=\inf\{t:\ |T_j(t)\cap \Gamma_\alpha\cap \Xi|\ge k+1\mbox { or }T_j(t)\cap \Gamma_\alpha=T_j(\infty)\cap \Gamma_\alpha,\ 1\le j\le 6\}\label{nearest1}
	\end{equation}
	and $\bar{R}\(x,\alpha\)=3t_{x,\alpha}$. We show that $\bar{R}$ is a radius of stabilization and its tail distribution can be bounded by an exponentially decaying function independent of $\alpha$ and $x$. For the the radius of stabilization, there are two cases to consider. The first case is that none of $T_j(t_{x,\alpha})\cap \Gamma_\alpha\cap\Xi$, $1\le j\le 6$, contains at least $k+1$ points, thus 
	$B(x, \bar{R}\(x,\alpha\))\supset \Gamma_\alpha$ and \Ref{defi4.1} is obvious. The second case is {at least }one of  $T_j(t_{x,\alpha})\cap \Gamma_\alpha\cap\Xi$, $1\le j\le 6$, contains at least $k+1$ points, which means that the $k$ nearest neighbors $\{x_1,\dots,x_k\}$ of $x$ are in $B\(x,t_{x,\alpha}\)$. {If a point $y\in \Gamma_\alpha\backslash B(x,t_{x,\alpha})$, then $y\in \Gamma_\alpha\cap (T_j(\infty)\backslash T_j(t_{x,\alpha}))$ for some $j$. Since $\Gamma_\alpha\cap (T_j(\infty)\backslash T_j(t_{x,\alpha}))$ is non-empty, $T_j(t_{x,\alpha})$ contains at least $k+1$ points $\{y_1,\dots, y_{k+1}\}$ and $d(x,y)>d(y_i,y)$ for all $i\le k+1$, then $y$ cannot have $x$ as one of its $k-$nearest neighbors. This ensures that all points having  $x$ as one of their $k-$nearest neighbors are in $B(x, t_{x,\alpha})$. Noting that the diameter of $B\(x,t_{x,\alpha}\)$ is $2t_{x,\alpha}$ and there are at least $k+1$ points in $B\(x,t_{x,\alpha}\)$, we can see that whether a point  $y$ in $B\(x,t_{x,\alpha}\)$ having $x$ as one of its $k-$nearest neighbors is entirely determined by $\Xi\cap B(y,2t_{x,\alpha})\subset \Xi\cap B(x,3t_{x,\alpha})$. This guarantees that $\eta\(x,\mathscr{X},\Gamma_\alpha\)$ is $\Xi_{B(x,3t_{x,\alpha})}$ measurable and $\bar{R}\(x,\alpha\)$ is a radius of stabilization. For the tail distribution of $\bar{R}$, referring to Figure~\ref{kntriangle},
		we consider the number of points of $\Xi$ falling into a triangle $A_t$ as a result of a sector being chopped off by the edge of $\Gamma_\alpha$. This is the worst situation for capturing the number of points by one sector. A routine trigonometry calculation gives that the area of $A_t$ is at least $0.116t^2$.  Define $\tau:=\inf\{t:\ |\Xi\cap A_t|\ge k+1\}$, then
		\begin{equation}\prob\(\bar{R}\(x,\alpha\)>t\)\le 6\prob\(\tau>t/3\)\le 6e^{-0.116\lambda  (t/3)^2}\sum_{i=0}^k\frac{\(0.116\lambda  (t/3)^2\)^i}{i!},\ t>0,\label{nearest0}\end{equation}
		which implies the exponential stabilization.
		
		The non-singularity \Ref{non-sinr} can be proved through the corresponding unrestricted score function $\bar{\eta}(x,\mathscr{X})=\frac{1}{2}\sum_{y\in \mathscr{X}}\|y-x\|\mathbf{1}_{\{(x,y)\in NG(\mathscr{X})\}}.$ Referring to Figure~\ref{knn}, we take $N_0:=B(0,{0.5})$, observe that $\partial B(0,6)$ can be covered by finitely many $B(x,3)$ with $\|x\|=5$ and write the centers of these balls as $x_1$, $\dots$, $x_n$ (in two dimensional case, $n=5$). Let
		$E$ be the event that $|B(x_i,1)\cap \Xi|\ge k+1$ for all $1\le i\le n$,  $|\(B(0,{1})\backslash B(0,{0.5})\)\cap\Xi|=k$ and $\Xi\cap \(B(0,6)\backslash \(\cup_{i\le n}B(x_i,1)\cup B(0,{1})\)\)$ is empty, then $E$ is $\Xi_{N_0^c}$ measurable and $\prob(E)>0$. Conditional on $E$, we we can see that $E_1:=\{|\Xi\cap B(0,{0.5})|=1\}$
		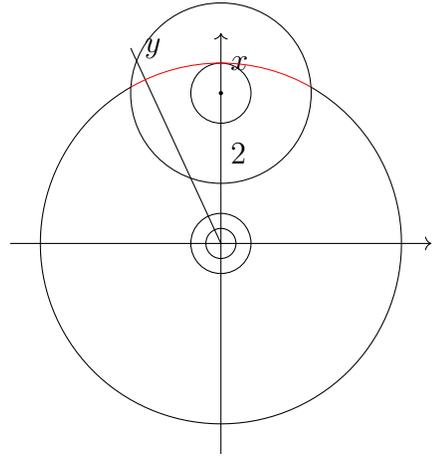
\begin{wrapfigure}{r}{6cm}
			\begin{tikzpicture}[scale=0.4]
				\draw[->](-7,0)--(7,0);
				\draw[->](0,-7)--(0,7);
				\draw[-](0,0)--(-3,6.5);
				\node[label={[black]0:$y$}] at (-3.25,6.5) {};
				\draw (0,0) circle (1);
				\draw (0,0) circle (0.5);
				\draw (0,5) circle (1);
				\fill (0,5) circle (2pt);
				\node[label={[black]0:$x$}] at (-0.4,6) {};	
				\node[label={[black]0:$2$}] at (-0.4,3) {};	
				\draw (0,5) circle (3);
				\draw[red] (2.99332590942,5.2) arc (60.07356513:119.92643487:6);
				\draw[black] (-2.99332590942,5.2) arc (119.92643487:420.07356513:6);
				
			\end{tikzpicture}
			\caption{$k$-nearest: non-singular condition}
			\label{knn}
		\end{wrapfigure}
		satisfies $\prob(E_1|E)>0$ and on $E_1$, all summands in \Ref{non-sinr} that are random are those involving the point of $\Xi_{B(0,{0.5})}$ and we now establish that these random score functions are entirely determined by $\Xi_{B(0,{1})}$. As a matter of fact, any point in $\Xi\cap B(x_i,1)$, $1\le i\le n$, has $k$ nearest points with distances no larger than $2$, so points in $\Xi\cap B(0,{1})$ cannot be $k$ nearest points to points in $\Xi\cap B(x_i,1)$. For any point $y\in \Xi\cap B(0,6)^c$, the line between $0$ and $y$ intersects $\partial B(0,6)$ at $y'$ which is in $B(x_{i_0},3)$ for some $1\le i_0\le n$, so the distances between $y$ and points in $\Xi\cap B(x_{i_0},3)$ are at most $\|y-y'\|+4=\|y\|-2$, the distances between $y$ and points in $\Xi\cap B(0,{1})$ are at least $\|y\|-1$, which ensures points of $\Xi\cap B(0,6)^c$ cannot have points in $\Xi\cap B(0,{1})$ as their $k$ nearest neighbors. {On the other hand, on $E\cap E_1$, there are $k+1$ points in $\Xi\cap B(0,1)$, for any point $x\in B(0,1)$, since the distances between $x$ and other points in $B(0,1)$ are less than $2$, the points outside $B(x,2)\subset B(0,3)$ will not be $k-$nearest points of $x$, $\bar{\eta}(x,\Xi)=\frac{1}{2}\sum_{y\in \(\Xi\cap B(0,1)\)\backslash\{x\}}\|x-y\|$.   }Hence, given $E$, all random score functions contributing in the sum of \Ref{non-sinr} are those completely determined by $\Xi_{B(0,2)}$, giving 
		$$\mathbf{1}_{E_1}\sum_{\substack{x\in\Xi\mbox{\scriptsize{ such that }}\\
				\bar{\eta}\(x,\Xi\)\mbox{\scriptsize{ is random given}~}\Xi_{N_0^c}}}\bar{\eta}\(x,\Xi\)\mathbf{1}_{d(x, N_0)<R(x)}=\mathbf{1}_{E_1}\left\{\sum_{\substack{y\in \Xi\cap B(0,{0.5}),\\ x \in  \(B(0,{1})\backslash B(0,{0.5})\)\cap\Xi}}\|x-y\|+X\right\},$$
		where $X$ is $\sigma(\Xi_{N_0^c})$ measurable. 
		Since this is a continuous function of $y\in \Xi\cap B(0,{0.5})$, the non-singularity \Ref{non-sinr} follows.
		
		For the moment condition \Ref{thm2.1r}, recalling the definition of $t_{x,\alpha}$ in \Ref{nearest1}, we replace $x$ with $x_3$ to get $t_{x_3,\alpha}$. We now establish that
		\begin{equation}\eta\(x_3,\Xi_{\Gamma_{\alpha}}+a_1\delta_{x_1}+a_2\delta_{x_2}+\delta_{x_3}\)\le 3.5kt_{x_3,\alpha}.
			\label{nearest2}
		\end{equation}
		In fact, the $k$ nearest neighbors to $x_3$ have the contribution of the total edge length $\le \frac12 kt_{x_3,\alpha}$. On the other hand, for $1\le j\le 6$, each point in $\(\Xi_{\Gamma_{\alpha}}+a_1\delta_{x_1}+a_2\delta_{x_2}\)\cap T_{j}\(t_{x_3,\alpha}\)$ may take $x_3$ as its $k$ nearest neighbor, which contributes to the total edge length $\le \frac12 t_{x_3,\alpha}$. As there are six sectors  $\(\Xi_{\Gamma_{\alpha}}+a_1\delta_{x_1}+a_2\delta_{x_2}\)\cap T_{j}\(t_{x_3,\alpha}\)$, $1\le j\le 6$, and each sector has no more than $k$ points with $x_3$ as their $k$ nearest neighbors, the contribution of the total edge length from this part is bounded by $3kt_{x_3,\alpha}$. By adding up these two bounds, we obtain \Ref{nearest2}. 
		Finally, we combine \Ref{nearest2} and \Ref{nearest0} to get
		\begin{align*}
			&\mathbb{E}\left\{\eta\(x_3,\Xi_{\Gamma_{\alpha}}+a_1\delta_{x_1}+a_2\delta_{x_2}+\delta_{x_3}\)^3\right\}
			\le 42.875k^3\mathbb{E}\left\{t_{x_3,\alpha}^3\right\}\le C <\infty,
		\end{align*}
		and the proof is completed by applying Theorem~\ref{thm2a}.
		\qed
		
		\subsection{ The total edge length of Voronoi tessellation}
		
		\begin{wrapfigure}{r}{0.35\textwidth}
			\vspace{-20pt}
			\begin{center} 
				\includegraphics[trim = 60mm 90mm 60mm 90mm, clip,width=0.35\textwidth]{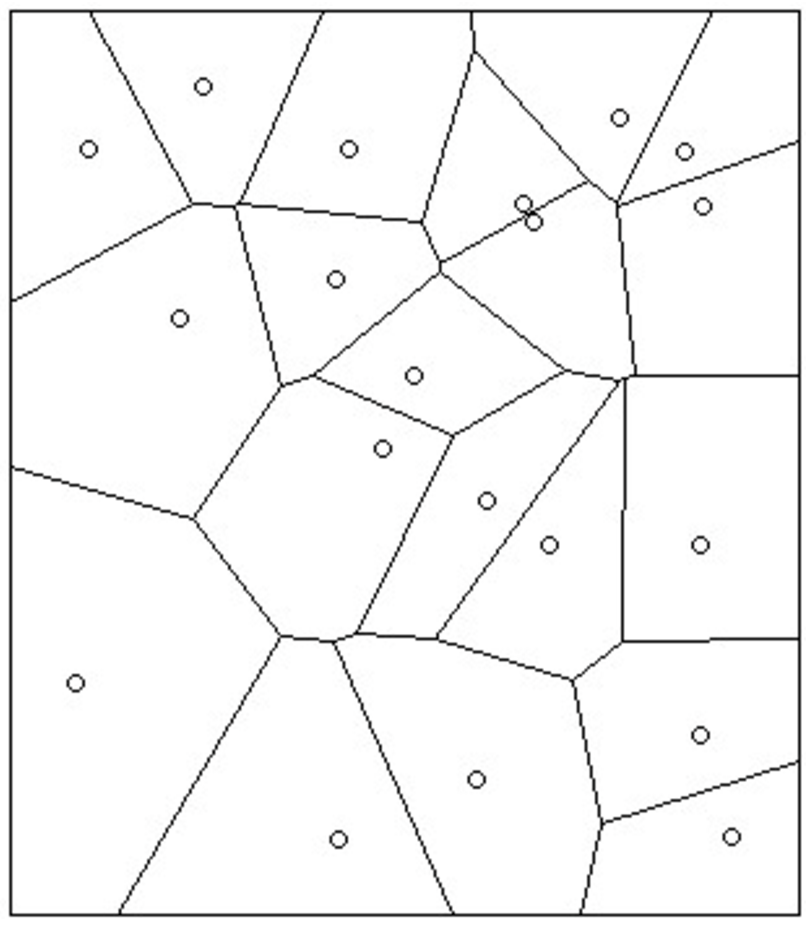}
				\vspace{-20pt}
				\caption{Voronoi tessellation} %
				\label{figure1}
			\end{center}
			\vspace{-15pt}
		\end{wrapfigure}
		
		Consider a finite point set $\mathscr{X}\subset \Gamma_\alpha$, the Voronoi tessellation in $\Gamma_\alpha$ generated by $\mathscr{X}$ is the partition formed by cells $C(x_i,\mathscr{X})\cap \Gamma_\alpha$, see Figure~\ref{figure1}. We write the graph of this tessellation as $V(\mathscr{X},\alpha)$ and the total edge length of $V(\mathscr{X},\alpha)$ as $\mathscr{V}(\mathscr{X},\alpha)$.
		
		\begin{thm}\label{exthm1}
			If $\Xi$ is a homogeneous Poisson point process, then $$d_{TV}\(\mathscr{V}(\Xi_{\Gamma_{\alpha}},\alpha), {\bar{Z}}_\alpha\)\le O\(\alpha^{-\frac{1}{2}}\ln(\alpha)^{\frac{5d}{2}}\),$$ where ${\bar{Z}}_\alpha$ is a normal random variable with the same mean and variance as those of $\mathscr{V}(\mathscr{X},\alpha)$.
		\end{thm}
		\noindent{\it Proof.} Before going into details, we observe that $$\mathscr{V}\(\Xi_{\Gamma_{\alpha}},\alpha\)=l(\partial \Gamma_\alpha)+\sum_{\{x,y\}\subset \Xi_{\Gamma_{\alpha}},x\neq y}l\(\partial C(x, \Xi_{\Gamma_{\alpha}})\cap \partial C(y, \Xi_{\Gamma_{\alpha}})\),$$ where $l(\cdot)$ is the volume of a set in dimension $d-1$. We restrict our attention to Voronoi tessellations of random point sets in $\mathbb{R}^2$ and, with notational complexity, the approach also works in $\mathbb{R}^d$.  Because $l(\partial \Gamma_\alpha)=4\alpha^{\frac{1}{2}}$ is a constant, by removing this constant, we have $\mathscr{V}'\(\Xi_{\Gamma_{\alpha}},\alpha\):=\mathscr{V}\(\Xi_{\Gamma_{\alpha}},\alpha\)-4\alpha^{\frac{1}{2}}$ and $d_{TV}\(\mathscr{V}\(\Xi_{\Gamma_{\alpha}}, \alpha\),{\bar{Z}}_\alpha\)=d_{TV}\(\mathscr{V}'\(\Xi_{\Gamma_{\alpha}},\alpha\),{\bar{Z}}_\alpha'\)$, where ${\bar{Z}}_\alpha'$ is a normal random variable having the same mean and variance as those of $\mathscr{V}'\(\Xi_{\Gamma_{\alpha}},\alpha\)$. We can set the score function corresponding to $\mathscr{V}'$ as $${\eta(x,\mathscr{X},\Gamma_\alpha)=\frac{1}{2}\sum_{y\in \mathscr{X},~y\neq x}l\(\partial C(x,\mathscr{X})\cap \partial C(y, \mathscr{X})\)=\frac{1}{2}l\(\partial (C(x, \mathscr{X})\cap \Gamma_\alpha)\backslash (\partial \Gamma_\alpha)\)}$$ for all $x\in \mathscr{X}\subset \Gamma_\alpha$, i.e., $\eta(x,\mathscr{X},\Gamma_\alpha)$ is a half of the total length of edges of $C(x, \mathscr{X})\cap \Gamma_\alpha$ excluding the boundary of $\Gamma_\alpha$. The score function $\eta$ is clearly translation invariant, thus, to apply Theorem~\ref{thm2a}, we need to verify stabilization as in Definition~\ref{defi4r},
		moment condition~\Ref{thm2.1r} and non-singularity~\Ref{non-sinr}.

		We start from showing that the score function is exponentially stabilizing. Referring to Figure~\ref{vdr}, similar to Section~\ref{knear}, we construct six disjoint equilateral triangles $T_{xj}(t)$, $1\le j\le 6$, such that $x$ is a vertex of these triangles and the triangles are rotated so that all edges with $x$ as a vertex have a minimal angle $\pi/12$ against the edges of $\Gamma_\alpha$. Let $T_{xj}(\infty)=\cup_{t\ge 0} T_{xj}(t)$, $1\le j\le 6$, then  $\cup_{1\le j\le 6}T_{xj}(\infty)=\mathbb{R}^2$.
		Define $$R_{xj}:=R_{xj}\(x,\alpha,\Xi_{\Gamma_\alpha}\):=\inf\{t:\ T_{xj}(t)\cap \Xi_{\Gamma_\alpha}\ne \emptyset\mbox{ or }T_{xj}(t)\cap\Gamma_\alpha=T_{xj}(\infty)\cap\Gamma_\alpha\}$$ 
		and $$R_{x0}:=R_{x0}\(x,\alpha,\Xi_{\Gamma_\alpha}\):=\max_{1\le j\le 6}R_{xj}\(x,\alpha,\Xi_{\Gamma_\alpha}\).$$
		We note that there is a minor issue of the counterpart of $R_{x0}$ defined in \cite{MY99} when $x$ is close to the corners of $\Gamma_\alpha$. 
		We now show that $\bar{R}(x,\alpha):=3R_{x0}(x,\alpha,\Xi_{\Gamma_\alpha})$ is a radius of stabilization. In fact, for any point $x'$ in $\Gamma_\alpha\backslash \overline{\left(\cup_{1\le j\le 6}T_{xj}(R_{x0})\right)}$, $x'$ is contained in a triangle $T_{xj_0}(\infty)\backslash\overline{T_{xj_0}(R_{x0})}$. This implies $T_{xj_0}(R_{x0})\cap \Xi_{\Gamma_\alpha}\ne\emptyset$, i.e., we can find a point $y\in T_{xj_0}(R_{x0})\cap \Xi_{\Gamma_\alpha}$ and the point $y$ satisfies $d(x',y)\le d(x,x')$, hence $x'\notin C(x,\Xi_{\Gamma_\alpha})\cap \Gamma_\alpha$, which ensures $C(x,\Xi_{\Gamma_\alpha})\subset \overline{\left(\cup_{1\le j\le 6}T_{xj}(R_0)\right)}$. Consequently, if a point $y$ in $\Xi_{\Gamma_\alpha}$ generates an edge of $C(X)\cap\Gamma_\alpha$, $d(x,y)\le 2R_0$ and $\bar{R}(x,\alpha)$ satisfies Definition~\ref{defi4r}. As in Section~\ref{knear}, we use $A_t$ in Figure~\ref{kntriangle} again to define $\tau:=\inf\{t:\ |\Xi\cap A_t|\ge 1\}$, then
		\begin{equation}\prob\(\bar{R}\(x\)>t\)\le 6\prob\(\tau>t/3\)\le 6e^{-0.116\lambda  (t/3)^2},\ t>0.\label{voronoi0}\end{equation}
		This completes the proof of the exponential stabilization of $\eta$.
		
		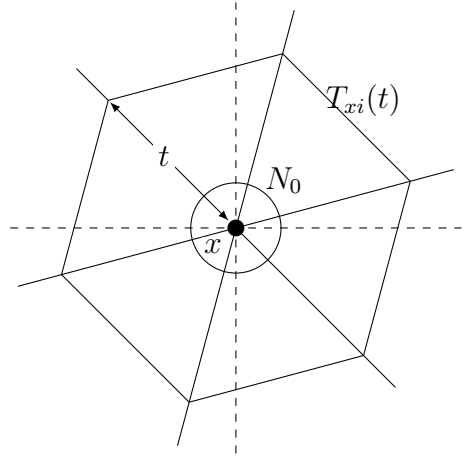
\begin{wrapfigure}{r}{6cm}
			\begin{tikzpicture}[scale=3]
				\newdimen\R
				\R=0.8cm
				\draw[yshift=0\R] (15:\R) \foreach \x in {15,75,135,...,315} {
					-- (\x:\R)
				} -- cycle (90:\R) node[above] {} ;
				\filldraw (0,0) circle (1pt);
				\foreach \x in {15,75,195,255,315}
				\draw[] (0,0) -- (\x:1);
				\draw[] (-0.5657,0.5657) -- (135:1); 
				\node at (45:0.3) {$N_0$};
				\node at (45:0.8) {$T_{xi}(t)$};
				\draw[latex-latex] (135:0.045) -- (135:0.79);
				\draw (0,0) circle (0.2);
				\node[fill=white,inner sep=2pt] at (135:0.45) {{$t$}};
				\node[fill=white,inner sep=1pt] at (219:0.13) {$x$};
				\draw[dashed] (-1,0)--(1,0) node[right]{};
				\draw[dashed] (0,-1)--(0,1) node[above]{};
			\end{tikzpicture}
			\caption{Voronoi: stabilization}
			\label{vdr}
			\vspace{-15pt}
		\end{wrapfigure}
		The non-singularity~\Ref{non-sinr} can be examined by using a non-restricted counterpart $\bar{\eta}$ of $\eta$, taking $N_0 =B(0,1)$ and filling the moat $B(0,4)\backslash B(0,3)$ with sufficiently dense points of $\Xi$ such that when $\Xi_{N_0^c}$ is fixed, the random score functions contributing to the sum of \Ref{non-sinr} are purely determined by a point in $\Xi\cap N_0$. More precisely, we cover the circle $\partial B(0,3)$ by disjoint squares with side length $\sqrt{2}/4$ and enumerate the squares as $S_i,\ 1\le i\le k$. Note that all the squares are contained in $B(0,4)\backslash B(0,2)$. Let $E=\cap_{1\le i\le k}\{|\Xi\cap(S_i)|\ge 1\}$, $E_1=\{|\Xi\cap N_0|=1\}$, then $E$ is $\sigma\(\Xi_{N_0^c}\)$ measurable, $\prob(E)>0$ and $\prob(E_1|E)>0$. Since the points in $\Xi\cap \(\cup_{i=1}^kS_i\)$ have neighbors within distance 1, for any $x\in N_0$, $T_{xj}(6)$ contains as least one point from $\Xi\cap \(B(0,4)\backslash B(0,2)\)$. As argued in the stabilization, points in $\Xi\cap \(B(0,12)^c\)$ do not affect the cell centered at $x\in N_0$, and by symmetry, $x\in N_0$ does not affect the Voronoi cells centered at points in $\Xi\cap \(B(0,12)^c\)$. This ensures that, conditional {on} $E$, all random score functions contributing in the sum of \Ref{non-sinr} are those completely determined by $\Xi_{N_0}$, giving 
		\begin{align*}&\mathbf{1}_{E_1}\sum_{\substack{x\in\Xi\mbox{\scriptsize{ such that }}\\
					\bar{\eta}\(x,\Xi\)\mbox{\scriptsize{ is random given}~}\Xi_{N_0^c}}}\bar{\eta}\(x,\Xi\)\mathbf{1}_{d(x, N_0)<R(x)}\\
			=&\mathbf{1}_{E_1}\left\{\sum_{x\in \Xi\cap B(0,12)}\bar{\eta}\(x,\Xi\)\mathbf{1}_{d(x, N_0)<R(x)}+X\right\},\end{align*}
		where $X$ is $\sigma\(\Xi_{N_0^c}\)$ measurable. As $\mathbf{1}_{E_1}\bar{\eta}\(x,\Xi\)$ is an almost surely (in terms of the volume measure in $\mathbb{R}^2$) continuous function of $x\in \Xi\cap N_0$, the proof of non-singularity \Ref{non-sinr} is completed.

		It remains to show the moment condition \Ref{thm2.1r}. In fact, as shown in {the} stabilizing property, we can see that $\bar{R}(x,\alpha)$ will not increase when adding points, so $C(x_3,\Xi_{\Gamma_\alpha}+a_1\delta_{x_1}+a_2\delta_{x_2}+\delta_{x_3})\cap \Gamma_\alpha\subset B(x, \bar{R}(x,\alpha))$, then the number of edges of $C(x_3,\Xi_{\Gamma_\alpha}+a_1\delta_{x_1}+a_2\delta_{x_2}+\delta_{x_3})\cap \Gamma_\alpha$ excluding those in the edge of $ \Gamma_\alpha$  is less than or equal to $(\Xi_{\Gamma_\alpha}+a_1\delta_{x_1}+a_2\delta_{x_2})\(B\(x,\bar{R}(x{,\alpha})\)\)\le\Xi_{\Gamma_\alpha}\(B\(x,\bar{R}(x{,\alpha})\)\)+2$ and each of them has length less than $2\bar{R}(x,\alpha)$. To this end, we observe that $\Xi$ restricted to outside of $\cup_{j=1}^6T_{xj}(R_{xj})$ is independent of $\bar{R}(x{,\alpha})$, hence 
		\begin{align*}\Xi_{\Gamma_\alpha}\(B\(x,\bar{R}(x{,\alpha})\)\)\le&\Xi_{\Gamma_\alpha}\(B\(x,\bar{R}(x{,\alpha})\)\backslash \(\cup_{j=1}^6T_{xj}(R_{xj})\)\)+6\\
			\stackrel{\mbox{\scriptsize{ST}}}{\le}&\Xi'_{\Gamma_\alpha}\(B\(x,\bar{R}(x{,\alpha})\)\)+6,
		\end{align*}
		where $\stackrel{\mbox{\scriptsize{ST}}}{\le}$ stands for stochastically less than or equal to and $\Xi'$ is an independent copy of $\Xi$. Hence, using \Ref{voronoi0}, we obtain
		\begin{align*}
			&\mathbb{E}\(\eta\left(x_3,\Xi_{\Gamma_\alpha}+a_1\delta_{x_1}+a_2\delta_{x_2}+\delta_{x_3}\right)^3\)\\
			\le &\mathbb{E}\(\(\Xi'\(B\({x_3},\bar{R}(x{,\alpha})\)\)+8\)^3\(2\bar{R}({x_3}{,\alpha})\)^3\)\\
			\le &\int_0^\infty \sum_{i=0}^\infty (i+8)^3(2r)^3\frac{e^{-\lambda \pi r^2} (\lambda \pi r^2)^i}{i!}6e^{-0.116\lambda  (t/3)^2}*\(0.116\lambda/9\)2rdr\\
			\le& C<\infty,
		\end{align*}
		which ensures \Ref{thm2.1r}. The proof of Theorem~\ref{exthm1} is completed by using Theorem~\ref{thm2a}. \qed
		
		{\subsection{Log volume estimation}}
		
		Log volume estimation is an essential research topic in forest science and forest management \cite{C80,Li15}. This example demonstrates that, with the marks, our theorem can be used to provide an error estimate of normal approximation of the log volume distribution. To this end, it is reasonable to assume that in a given range $\Gamma_\alpha$ {of} {a} natural forest, the locations of trees form a Poisson point process $\bar{\Xi}$, and for $x\in \bar{\Xi}$, we can use a random mark $M_x\in{T}:=\{1,\dots,n\}$ to denote the species of the tree at position $x$, then $\Xi:=\sum_{x\in \bar{\Xi}}\delta_{(x,M_x)}$ is a marked Poisson point process with independent marks. The timber volume of a tree at $x$ is a combined result of the location, the species of the tree, the configuration of species of trees in a finite range around $x$ and some other random factors that can't be explained by the configuration of trees in the range. We write $\eta(\left(x,m\right), \Xi_{\Gamma_\alpha},\Gamma_\alpha)$ as the timber volume determined by the location $x$, the species $m$ and the configuration of trees, and $\e_{x}$ as the adjusted timber volume at location $x$ due to unexplained random factors.
		
		\begin{thm} Assume that $\eta$ is a non-negative bounded score function such that $$\eta(\left(x,m\right), \Xi_{\Gamma_\alpha},\Gamma_\alpha)=\eta(\left(x,m\right), \Xi_{\Gamma_\alpha\cap B(x,r)},\Gamma_\alpha)$$ for some positive constant $r$, $$\eta(\left(x,m\right), \Xi_{B(x,r)},\Gamma_{\alpha_1})=\eta(\left(x,m\right), \Xi_{B(x,r)},\Gamma_{\alpha_2})$$ for all $\alpha_1$ and $\alpha_2$ with $B(x,r)\subset \Gamma_{\alpha_1\wedge\alpha_2}$, $\eta$ is translation invariant in Definition~\ref{traninvres0}, $\e_{x}$'s are $\iid$ random variables with finite third moment and the positive part $\e_x^+:=\e_x\vee 0$ is non-singular, and $\e_{x}$'s are independent of the configuration $\Xi$, then the log volume of the range $\Gamma_\alpha$ can be represented as 
			$${\bar{W}}_{\alpha}:=\sum_{x\in \bar{\Xi}_{\Gamma_\alpha}}\[\(\eta(\left(x,m\right), \Xi_{\Gamma_\alpha},\Gamma_\alpha)+\e_x\)\vee 0\]$$
			and it satisfies
			$$d_{TV}\({\bar{W}}_\alpha,{\bar{Z}}_\alpha\)\le O\(\alpha^{-\frac{1}{2}}\),$$ 
			where ${\bar{Z}}_\alpha$ is a normal random variable with the same mean and variance as those of ${\bar{W}}_\alpha$.
		\end{thm}
		
		\noindent{\it Proof.~}Before going into details, we first construct a new marked Poisson point process $\Xi':=\sum_{x\in \bar{\Xi}}\delta_{(x, (M_x, \e_x))}$ with $\iid$ marks $(M_x,\e_x)\in T\times \mathbb{R}$ independent of the ground process $\bar{\Xi}'=\bar{\Xi}$ and  incorporate $\e_x$ into a new score function on $\Xi'$ as
		\begin{align*}\eta'((x,{(m,\e)}),\Xi',\Gamma_\alpha):=& \eta'((x,(m,\e)),\Xi'_{\Gamma_\alpha},\Gamma_\alpha)\\
			:=&\[\(\eta(\left(x,m\right), \Xi_{\Gamma_\alpha},\Gamma_\alpha)+\e\)\vee 0\]\mathbf{1}_{(x,(m,\e))\in \Xi'_{\Gamma_\alpha}}.
		\end{align*}
		We can see that ${\bar{W}}_\alpha=\sum_{x\in \bar{\Xi}'_{\Gamma_\alpha}}\eta'((x,(m_x,\e_x)),\Xi'_{\Gamma_\alpha},\Gamma_\alpha)$. The score function $\eta'$ is clearly translation invariant, thus, to apply Theorem~\ref{thm2a}, it is sufficient to verify that $\eta'$ is range-bound as in Definition~\ref{defi4r}, satisfies the moment condition \Ref{thm2.1r} and non-singularity \Ref{non-sinr}.
		
		The range-bound property of the score function $\eta'$ is inherited from the range-bound property of $\eta$ with the same radius of stabilization $\bar{R}(x,\alpha):=r$, the moment condition \Ref{thm2.1r} is a direct consequence of the boundedness of $\eta$, the finite third moment of $\e_x$ and the Minkowski inequality, hence it remains to show the non-singularity. To this end, we observe that the corresponding unrestricted counterpart $\bar{\eta}$ of $\eta'$ is defined by \begin{align*}\bar{\eta}((x,(m_x,\e_x)),\mathscr{X}')&= \lim_{\alpha\rightarrow \infty}\eta'((x,(m_x,\e_x)),\mathscr{X}',\Gamma_{\alpha})\\&=\eta'((x,(m_x,\e_x)),\mathscr{X}',\Gamma_{\alpha_x})= \(\eta((x,m_x),\mathscr{X},\Gamma_{\alpha_x})+\e_x\)\vee 0
		\end{align*} where $\alpha_x=4(\|x\|+r)^2$ and $\mathscr{X}$ is the projection of $\mathscr{X}'$ on $\mathbb{R}^2\times T$. Let $N_0=B(0,1)$, $E=\left\{\left|\Xi'_{B(N_0,r)\backslash N_0}\right|=0\right\}$, $E_1=\left\{\left|\Xi'_{N_0}\right|=1\right\}$, then $E$ is $\sigma\(\Xi'_{N_0^c}\)$ measurable, $\mathbb{P}(E)>0$ and $\mathbb{P}(E_1|E)>0$. Writing $\e_x^-=-\min(\e_x,0)$, $\bar{\Xi}'_{N_0}=\{x_0\}$ in $E_1$, {given $E$} we have 
		\begin{align*}
			&\mathbf{1}_{E_1}\sum_{x\in\bar{\Xi'}}\bar{\eta}\((x,(m_x,\e_x)),\Xi'\)\mathbf{1}_{d(x, N_0)<r}\\
			=&\mathbf{1}_{E_1} \(\eta((x_0,m_{x_0}),\delta_{({x_0},m_{x_0})},\Gamma_{4(r+1)^2})+\e_{x_0}\)\vee0 \\
			=&\mathbf{1}_{E_1}\[\mathbf{1}_{\e_{x_0}>0}\(\eta(({x_0},m_{x_0}),\delta_{({x_0},m_{x_0})},\Gamma_{4(r+1)^2})+\e_x^+\)\right.\\
			&+\left.\mathbf{1}_{\e_x\le 0}\(\eta(({x_0},m_{x_0}),\delta_{({x_0},m_{x_0})},\Gamma_{4(r+1)^2})-\e_{x_0}^-\)\vee 0\].
		\end{align*} On $\{\e_{x_0}>0\}$, $\e_{x_0}^+$ has positive non-singularity part and is independent of \\ $\eta(({x_0},m_{x_0}),\delta_{({x_0},m_{x_0})},\Gamma_{4(r+1)^2})$, $\eta(({x_0},m_{x_0}),\delta_{({x_0},m_{x_0})},\Gamma_{4(r+1)^2})+\e_{x_0}^+$ is also non-singular, together with the fact that $\{\e_{x_0}>0\}$ and $\{\e_{x_0}\le 0\}$ are disjoint, the non-singularity follows. \qed
		
		\begin{re} If the timber volume of a tree is determined by its nearest neighboring trees, then we can set the score function $\eta$ as a function of weighted Voronoi cells. Using the idea of the proof of Theorem~\ref{exthm1}, we can establish the bound of error of normal approximation to the distribution of the log timber volume ${\bar{W}}_\alpha$ as $d_{TV}\({\bar{W}}_\alpha,{\bar{Z}}_\alpha\)\le O(\alpha^{-\frac{1}{2}}\ln(\alpha)^{\frac{5d}{2}})$.  \end{re}

\subsection{Maximal layers}

Maximal layers of points have been of considerable interest since \cite{R62,K75} and have a wide range of applications, see \cite{CHT03} for a brief review of their applications. One of the applications is \textit{the smallest color-spanning interval}~\cite{K17} which is a linear function of the distances between maximal points and the edge. In this subsection, we demonstrate  that Theorem~\ref{thm2a} with marks can be easily applied to estimate the error of normal approximation to the distribution of the sum of distances between different maximal layers if the points are from a Poisson point process.

 For $x\in \mathbb{R}^d$, we define $A_x=([0,\infty)^d+x)\cap\Gamma_\alpha$.
Given a locally finite point set $\mathscr{X}\subset \mathbb{R}^d$, a point $x$ is called {\it maximal} in $\mathscr{X}$ if $x\in\mathscr{X}$ and there is no other point $(y_1,\dots,y_d)\in \mathscr{X}$ satisfying $y_i\ge x_i$ for all $1\le i\le d$ (see Figure~\ref{mp}). 
 Mathematically, $x$ is maximal in $\mathscr{X}$ if $\mathscr{X}\cap A_x=\{x\}$. This enables us to write different maximal layers as follows: the $k$th maximal layer of points can be recursively defined as
 $$\mathscr{X}_k:=\sum_{x\in \mathscr{X}}\delta_x\bone_{[A_x\cap (\mathscr{X}\backslash (\cup_{i=1}^{k-1}\mathscr{X}_{i}))=\{x\}]},\ \ \ k\ge 1,$$ with the convention $\cup_{i=1}^0\mathscr{X}_i=\emptyset$.

For simplicity, we consider the restriction of the Poisson point process to a region
in $\mathbb{R}^d$ between two parallel $d-1$ dimensional planes for
$d\ge2$. More precisely, the region of interest is 
$$\Gamma_{\alpha,r}=\left\{(x_1,x_2,\dots,x_d);~x_i\in
[0,\alpha^{\frac{1}{d-1}}], i\le
d-1,~x_d+\sum_{i=1}^{d-1}x_i \cot(\theta_i)\in[0,r]\right\}$$ for fixed $\theta_i\in
(0,\frac{\pi}{2})$, $1\le i\le d-1$, and $\Xi_{\Gamma_{\alpha,r}}$ is a
homogeneous Poisson point process with rate $\lambda$ on
$\Gamma_{\alpha,r}$. Define $\Xi_{k,r,\alpha}$ as the $k$th maximal layer of
$\Xi_{\Gamma_{\alpha,r}}$, then the total distance between the points in $\Xi_{k,r,\alpha}$ and the upper plane
$$P:=\left\{(x_1,x_2,\dots,x_d);~x_i\in [0,\alpha^{\frac{1}{d-1}}],\ i\le d-1,x_d=-\sum_{i=1}^{d-1}x_i \cot(\theta_i)+r\right\}$$ can be represented as
${\bar{W}}_{k,r,\alpha}:=\sum_{x\in \Xi_{k,r,\alpha}}d(x,P)$.

\begin{thm}\label{exthm5}
With the above setup, when $r\in \mathbb{R}_+$ is fixed,
$$d_{TV}\({\bar{W}}_{k,r,\alpha},{\bar{Z}}_{k,r,\alpha}\)\le O\(\alpha^{-\frac{1}{2}}
\),$$ where ${\bar{Z}}_{k,r,\alpha}\sim N\(\mean {(\bar{W}}_{k,r,\alpha}{)},\var\({\bar{W}}_{k,r,\alpha}\)\)$.
\end{thm}

\begin{re}{\rm It remains a challenge to consider maximal layers induced by a homogeneous Poisson point process on $$\left\{(v_1,v_2):\ v_1\in [0,\alpha^{1/(d-1)} ]^{d-1},0\le v_2\le F(v_1)\right\},$$
where $F:\ [0,\alpha^{1/(d-1)} ]^{d-1}\to [0,\infty)$ has continuous negative partial derivatives in all coordinates, the partial derivatives are bounded away from 0 and $-\infty$, and $|F|\le {O(\alpha^{1/(d-1)})}$. We conjecture that normal approximation in total variation for the total distance between the points in a maximal layer and the upper edge surface is still valid.
}
\end{re}

\noindent{\it Proof of Theorem~\ref{exthm5}.~} As the score function $d(\cdot,P)$ is not translation invariant in
 the sense of Definition~\ref{traninvres0}, we first turn the problem to that of a marked Poisson point
process with independent marks. The idea is to project the points of $\Xi_{\Gamma_{\alpha,r}}$ on their first $d-1$ coordinates to obtain the ground Poisson point process and send the last coordinate to marks with $T=[0,r]$. To this end, define a mapping
$h':\Gamma_{\infty,r}:=\cup_{\alpha>0}\Gamma_{\alpha,r}\rightarrow
[0,\infty)^{d-1}\times [0,r]$ such that
$$h'(x_1,\dots,x_d)=(x_1,\dots,x_d)+\(0,\dots,0,\sum_{i=1}^{d-1}x_i
\cot(\theta_i)\)$$
and $h(\mathscr{X}):=\{h'(x):\ x\in \mathscr{X}\}$. Then $h'$ is a one-to-one mapping and $\Xi':=h\(\Xi_{\Gamma_{\infty,r}}\)$
can be regarded as a marked Poisson point process on
$[0,\infty)^{d-1}\times [0,r]$ with rate $r\lambda$ and independent
marks following the uniform distribution on $[0,r]$. Write the mark of $x\in \Xi'$
as $m_x$, then \begin{equation}\label{ex5.1}{\bar{W}}_{k,r,\alpha}=C(\theta_1,\dots,\theta_{d-1})\sum_{x\in
h(\Xi_{k,r,\alpha})}(r-m_x),
\end{equation} where $C(\theta_1,\dots,\theta_{d-1})$ is a constant
determined by $\theta_1,\dots,\theta_{d-1}$. Let
$\Gamma_\alpha':=[0,\alpha^{\frac{1}{d-1}}]^{d-1}$, then
$h(\Xi_{\Gamma_{\alpha,r}})=\Xi'_{\Gamma_\alpha'}$. For a point
$(x,m)\in [0,\infty)^{d-1}\times[0,r]$, we write
$A_{x,m,r,\alpha}' =h'(((h')^{-1}(x,m)+[0,\infty)^d)\cap \Gamma_{\alpha,r})$ (see Figure~\ref{mp}) and
\begin{equation}\label{ex5.2}\Xi'_{k,r,\alpha}:=h(\Xi_{k,r,\alpha})=\sum_{x\in
		\bar{\Xi}'_{\Gamma_\alpha'}}\delta_{(x,m_x)}\mathbf{1}_{A_{x,m_x,r,\alpha}'\cap\(\Xi'_{\Gamma_\alpha'}\backslash
		\cup_{i=1}^{k-1}\Xi'_{i,r,\alpha}\)=\{(x,m_x)\}}.
\end{equation} Combining \Ref{ex5.1} and \Ref{ex5.2}, ${\bar{W}}_{k,r,\alpha}$ can be represented as the sum of values of the score function
$$\eta((x,m_x),\Xi',\Gamma_\alpha'):=C(\theta_1,\dots,\theta_{d-1})(r-m_x)\mathbf{1}_{(x,m_x)\in
\Xi'_{k,\alpha}}$$ over the range $\Gamma_\alpha'$. To apply
Theorem~\ref{thm2a}, we need to check that $\eta$ is range-bound as in
Definition~\ref{defi4r}, satisfies the moment condition
\Ref{thm2.1r} and non-singularity \Ref{non-sinr}.

\begin{figure}
	\begin{subfigure}{.5\textwidth}
		\begin{tikzpicture}[scale=0.6]
			\draw (0,8) -- (8,0);
			\draw (0,8) -- (0,9);
			\draw[dashed](0,8)--(0,0);
			\draw[dashed](8,0)--(0,0);
			\draw(0,0.3)--(0.3,0.3);
			\draw(0.3,0)--(0.3,0.3);
			\draw (0,9) -- (8,1);
			\draw (8,1) -- (8,0);
			\draw[black] (8,0.2) arc (90:135:0.2);
			\node[label={[black]270:$\alpha$}] at (4,0) {};
			\node[label={[black]135:$\theta_1$}] at (8.55,-0.15) {};
			\node[label={[black]45:\textcolor{blue}{$(0,0)$}}] at (-2.16,7.12) {};
			\draw [decorate,decoration={brace,amplitude=3.5pt,mirror,raise=0pt},yshift=0pt]
			(8.05,0.0) -- (8.05,1) node [black,midway,xshift=0.25cm] {\footnotesize
				$r$};
			\draw (0.5100,7.9556) circle (4pt);
			\draw (2.5100,5.9556) circle (4pt);
			\draw (1.5100,7.2556) circle (4pt);
			\draw (3.5100,5.2556) circle (4pt);
			\draw (4.0525,4.7427) circle (4pt);
			\draw (5.4299,2.8470) circle (4pt);
			\draw (6.3619,1.9843) circle (4pt);
			\draw (5.9451,2.1521) circle (4pt);
			\draw (3.1378,5.6856) circle (4pt);
			\draw (1.3695,6.9476) circle (4pt);
			\draw (0.2547,7.7798) circle (4pt);
			\draw (5.6484,3.3019) circle (4pt);
			\draw[blue,fill=blue] (0,8) circle (.6ex);
			at (3.5100,5.2556)   () {};
			at (4.0525,4.7427)   () {};
			at (5.4299,2.8470)   () {};
			at (6.0619,1.9843)   () {};
			at (5.9451,2.1521)   () {};
			at (3.1378,5.6856)   () {};
			at (5.2438,3.4510)   () {};
			at (1.3695,6.9476)   () {};
			at (5.6484,3.3019)   () {};
			at (0.2547,7.7798)   () {};
			\node[label={[black]45:$A_x\cap \Gamma_{\alpha,r}$}] at (5,3.05) {};
			\draw (5.06,3.94) -- (4.95,3.83);
			\draw (5.17,3.83) -- (4.95,3.61);
			\draw (5.28,3.72) -- (5.06,3.5);
			\draw (5.39,3.61) -- (5.28,3.5);
			\draw (4.5,3.5) -- (4.85,3.5);
			\draw (5.05,3.5) -- (5.5,3.5);
			\draw (4.95,3.6) -- (4.95,4.05);
			\draw (4.95,3.05) -- (4.95,3.4);
			\node[label={[black]135:$x$}] at (5.3,3.25) {};
			\node[circle,draw,inner sep=2pt] at (4.95,3.5)   () {};
			\node[label={[black]45:{$A'_{x,m,r,\alpha}$}}] at (5.2,7.8) {};
			\draw (5.06,3.94) -- (4.95,3.83);
			\draw (5.05,8.55) -- (5.5,9);
			\draw (4.95,9) -- (5.5,9);
			\draw (4.95,8.55) -- (4.95,9);
			\draw (5.06,8.59) -- (5.06,9);
			\draw (5.17,8.70) -- (5.17,9);
			\draw (5.28,8.81) -- (5.28,9);
			\draw (5.39,8.92) -- (5.39,9);
			\node[label={[black]135:$h'(x)$}] at (5.2,7.8) {};
			\node[circle,draw,inner sep=2pt] at (4.95,8.54)   () {};
		\end{tikzpicture}
		\caption{maximal point}
		\label{mp}
	\end{subfigure}
	\begin{subfigure}{.5\textwidth}
		\begin{tikzpicture}[scale=0.9]
			\draw (0,0.56) 
			-- (0.7,.56) 
			-- (0.7,1.12)
			-- cycle;
			\draw (0.7,1.68) 
			-- (1.4,1.68) 
			--  (1.4,2.24)
			-- cycle;
			\draw   [thick,dash dot] (-3.5,-2.8) -- (-3.5,2.8)
			--  (3.5,2.8) -- (3.5,-2.8)
			-- (0,-2.8);
			\draw   [thick,dashed] (0,0) -- (3.5,2.8);
			\node[label={[black]45:\scriptsize{$B_1$}}] at (0,0.28) {};
			\node[label={[black]45:\scriptsize{$B_2$}}] at (0.7,1.39) {};
			\node[label={[black]45:$\left(\frac{r\tan(\theta_1)}2,\frac{r}2\right)$}] at (-.18,-0.7) {};
			\node[label={[black]45:$B_0$}] at (-1.5,-2.3) {};
			\draw (-3.5,-2.8)
			-- (0,0)
			-- (0,-2.8)
			-- cycle;
			\draw[blue,fill=blue] (0,0) circle (.4ex);   
			\definecolor{c1}{RGB}{0,129,188}
			\definecolor{c2}{RGB}{252,177,49}
			\draw[fill=c1,opacity=0.15]   (0,0.56) -- (0.7,.56) -- (0.7,1.12);
			\draw[fill=c1,opacity=0.15]   (0.7,1.68) -- (1.4,1.68) --  (1.4,2.24);
			\draw[fill=c1,opacity=0.15]   (-3.5,-2.8) -- (0,0) -- (0,-2.8);
		\end{tikzpicture}
		\vskip0.8cm
		\caption{singularity}
		\label{mp1}
	\end{subfigure}
	\caption{maximal layers}
\end{figure}
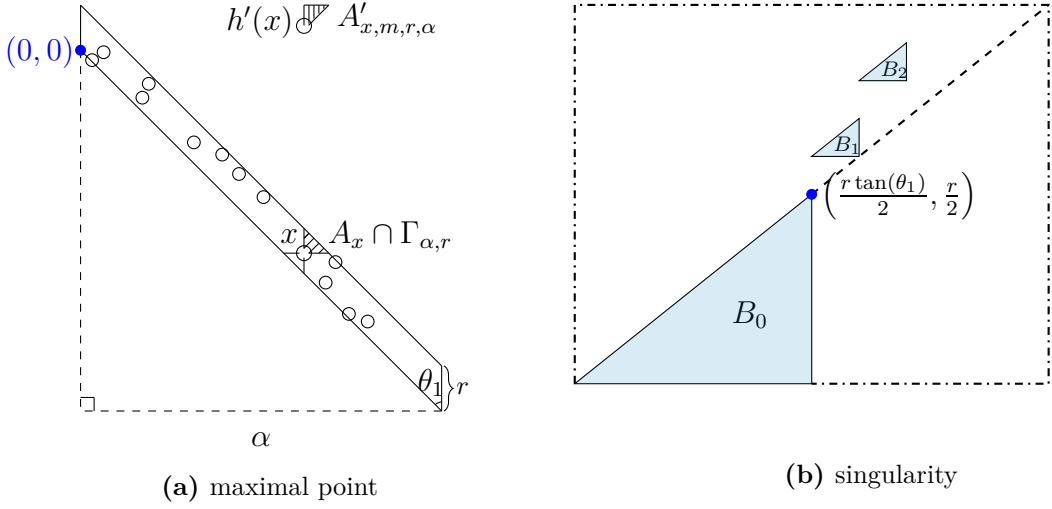
For simplicity, we only show the claim in two dimensional case and
the argument for $d>2$ is the same except notational complexity. When $d=2$,
$\Gamma_{\alpha,r}$ is a parallelogram with angle $\theta_1$ as in
Figure~\ref{mp}, $P$ and $C(\theta_1,\dots,\theta_{d-1})$ reduce to an
edge in $\mathbb{R}^2$ and $\sin(\theta_1)$ respectively. Since
$\sin(\theta_1)(r-m_x)$ is given by the mark of $x$, to show $\eta$ is
range-bound, it is sufficient to show that $\mathbf{1}_{\{(x,m_x)\in
\Xi'_{k,r,\alpha}\}}$ is completely determined by $\Xi'\cap
A_{x,m_x,r,\alpha}'$. In fact, we can accomplish this by observing that $(x,m_x)\in \Xi'_{k,r,\alpha}$ iff there is a sequence $\{(x_j,m_{x_j}),\ 1\le j\le k\}\subset \Xi'\cap
A_{x,m_x,r,\alpha}'$ (which ensures $A_{x_j,m_{x_j},r,\alpha}'\subset A_{x,m_x,r,\alpha}'$) such that $(x_k,m_{x_k})=(x,m)$ and $A_{x_j,m_{x_j},r,\alpha}'\cap(\Xi'\backslash
\cup_{i=1}^{j-1}\Xi'_{i,r,\alpha})=\{(x_j,m_{x_j})\}$ for $1\le j\le k$. Since $\Xi'\cap
A_{x,m_x,\alpha}'\subset \Xi'_{[x,x+r\tan(\theta_1)]}$, we can see that $\eta$ is
range-bound in Definition~\ref{defi4r} with
$\bar{R}(x,\alpha):=r\tan(\theta_1){+1}$. The moment condition follows from the fact
that $\eta$ is bounded above by $r$. For the non-singular condition, {we extend $\Xi'$ to $\mathbb{R}^{d-1}\times[0,r]$, write $\Gamma_{\infty,r}^e:=\{x\in\mathbb{R}^d: \mbox{there exists } y\in P~\mbox{such that }y-r(0,\dots,0,1)\le x\le y\}$ and let $(\Xi_{\Gamma_{\infty,r}^e})_j$ be the $j$th maximal layer of $\Xi_{\Gamma_{\infty,r}^e}$ and $\Xi'_j=h((\Xi_{\Gamma_{\infty,r}^e})_j)$},
we can see that the corresponding unrestricted score function is
$\bar{\eta}(x,\Xi')=\sin(\theta_1)(r-m_x)\mathbf{1}_{(x,m_x)\in
	\Xi'_k}$ and the corresponding stabilizing radii
$R(x)=r\tan(\theta_1){+1}$. Referring to Figure~\ref{mp1}, we set 
$N_0:=\(0,\frac{r\tan(\theta_1)}{2}\)$, $B_0=\{(x,m);~x\in N_0,0\le m\le
x\cot(\theta_1)\}$, $B_i$ as the triangle region with vertices $\(r\tan(\theta_1)\left(\frac12+\frac{i-1}{4(k-1)}\right), r\left(\frac12+\frac{2i-1}{4(k-1)}\right)\),$ $\(r\tan(\theta_1)\left(\frac12+\frac{i}{4(k-1)}\right),r\left(\frac12+\frac{2i-1}{4(k-1)}\right)\)$ and $\(r\tan(\theta_1)\left(\frac12+\frac{i}{4(k-1)}\right), r\left(\frac12+\frac{2i}{4(k-1)}\right)\)$, $1\le i\le k-1$, and 
$${C}=\({\(\[-r\tan(\theta_1){-1},\frac{3r\tan(\theta_1)}{2}{+1}\]\backslash N_0\)}\times[0,r]\)\backslash\(\cup_{i=1}^{k-1}B_i\),$$
define 
$E:=\left\{\Xi'\cap {C}=\emptyset, \left|\Xi'\cap B_i\right|=1,1\le i\le k-1\right\}$  and $E_0:=\{|\Xi'_{N_0}|=\left|\Xi'\cap B_0\right|=1\}$, then {$E\in\sigma\(\Xi'_{N_0^c}\)$, }$\mathbb{P}(E)>0$ and $\mathbb{P}(E_0|E)>0$. {We can see that given $E$, the point in $\Xi'\cap B_i$ is in $\Xi'_{k-i}$ for all $1\le i\le k-1$}.
Moreover, on $E\cap E_0$, the point $({x_0},m_{x_0})$ in $\Xi'_{N_0}$ is in
$\Xi'_k.$ Hence, given $E$, 
$$\mathbf{1}_{E_0}\sum_{x\in\bar{\Xi'}}\bar{\eta}\(x,\Xi'\)\mathbf{1}_{d(x,
	N_0)<R(x)}=\mathbf{1}_{E_0}\sin(\theta_1)(r-m_{x_0})$$
is non-singular. \qed
 
As a final remark of the section, we mention that the unrestricted version of all examples considered here can be proved because it is trivial to show that the unrestricted version of the score function $\bar{\eta}$ satisfies the stabilization condition as in Definition~\ref{defi4} and moment condition \Ref{thm2.1} using the same method, and non-singularity \Ref{non-sin} condition is the same as \Ref{non-sinr} given the score function is $\bar{\eta}$.

\section{Preliminaries and auxiliary results}\label{Preliminaries}

We start with a few technical lemmas.

\begin{lma} \label{lma1} Assume $\xi_1,\dots,\xi_n$ are $\iid$ random variables having the triangular density function
\begin{equation}\kappa_a(x)=\left\{\begin{array}{ll}
\frac1a\left(1-\frac {|x|}a\right),&\mbox{ for }|x|\le a,\\
0,&\mbox{ for } |x|>a,
\end{array}\right.\label{lma1.1}\end{equation}
where $a>0$. Let $T_n=\sum_{i=1}^n\xi_i$. Then for any $\gamma>0$,
\begin{equation}d_{TV}(T_n,T_n+\gamma)\le \frac{\gamma}a\left\{\sqrt{\frac{3}{\pi n}}+\frac{2}{(2n-1)\pi^{2n}}\right\}.\label{lma1.2}\end{equation}
\end{lma}

The following lemma says that if the distribution of a random variable is non-singular, then the distributions of random variables which are not far away from it are also non-singular.

\begin{lma}\label{non-singular1} Let $F$ be a non-singular distribution on $\mathbb{R}$ with $\alpha_F>0$ in the decomposition \Ref{decom1}. Then for any distribution $G$ such that $d_{TV}\(F, G\)<\alpha_F$, $G$ is non-singular with $\alpha_G\ge \alpha_F-d_{TV}\(F, G\)$ in its representation
$$G=(1-\alpha_G) G_s+\alpha_GG_a,$$
where $G_a$ is absolutely continuous with respect to the Lebesgue measure and $G_s$ is singular.
\end{lma}

We denote the convolution by $\ast$ and write $F^{k\ast}$ as the $k$-fold convolution of the function $F$ with itself.

\begin{lma}\label{lma2}
		For any two non-singular distributions $F_1$ and $F_2$, there exist positive constants $a>0$, $u\in\real$, $\theta\in (0,1]$ and a distribution function $H$ such that 
		\begin{equation} \label{lma2.1}
		F_1\ast F_2=(1-\theta)H+\theta K_a\ast\d_{u},
		\end{equation}
		where $K_a$ is the distribution of the triangle density $\kappa_a$ in \Ref{lma1.1} and $\d_u$ is the Dirac measure at $u$.
	\end{lma}
	
Lemma~\ref{lma2} says that $F_1\ast F_2$ is the distribution function of $(X_1+u)X_3+X_2(1-X_3)$, where $X_1\sim K_a$, $X_2\sim H$, $X_3\sim{\rm Bernoulli}(\theta)$ are independent random variables.

\begin{re} \label{remark1}{\rm~From the definition of triangular density function, if $a$, $u$, $\theta$ satisfy \eqref{lma2.1} with a distribution $H$, then for arbitrary $p,\ q$ such that $0<q\le p\le 1$, we can find an $H'$ satisfying the equation with $a'=pa$, $u'=u$ and $\theta'=q\theta$.}
	\end{re}

Using the properties of the triangular distributions, we can derive that the sum of the score functions restricted by their radii of stabilization has a similar property as shown in Lemma~\ref{lma1} when the score function is range-bound, exponentially stabilizing or polynomially stabilizing with suitable $\beta$.

	\begin{lma}\label{lma3} Let $\Xi$ be a marked homogeneous Poisson point process on $(\mathbb{R}^d\times T,\mathscr{B}(\real^d)\times \mathscr{T})$  with intensity $\lambda$ and $\iid$ marks in $(T,\mathscr{T})$ following $\mathscr{L}_{T}$. 
	\begin{description}
	\item{(a)} (unrestricted case) Assume that the score function $\eta$ is non-singular~\Ref{non-sin}. If $\eta$ is polynomially stabilizing in Definition~\ref{defi4} with order $\beta>d+1$ and the radius of stabilization $R$, define $W_{\alpha,r}:=\sum_{(x,m){\in \Xi_{\Gamma_\alpha}}}\eta(\left(x,m\right), \Xi)\mathbf{1}_{R(x)\le r}$, then
		\begin{align}
		d_{TV}(W_{\alpha,r},W_{\alpha,r}+\gamma)&\le C(|\gamma|\vee 1)\left(\alpha^{-\frac{1}{2}}r^{\frac{d}{2}}\right)	\label{statement1}	\end{align}
	 for any $\gamma\in \mathbb{R}$ and $r>R_0$,
		where $C$ and $R_0$ are positive constants independent of $\gamma$. If $\eta$ is range-bound in Definition~\ref{defi4}, then
		\begin{align}d_{TV}(W_\alpha, W_\alpha+\gamma)\le C(|\gamma|\vee 1)\alpha^{-\frac{1}{2}}\label{statement2}\end{align} for some positive constant $C$ independent of $\gamma$.
		
		\item{(b)} (restricted case) Assume that the score function $\eta$ is non-singular~\Ref{non-sinr}. If $\eta$ is polynomially stabilizing in Definition~\ref{defi4r} with order $\beta>d+1$ and the radius of stabilization $R$, define $\bar{W}_{\alpha,r}:={\sum_{(x,m)\in \Xi_{\Gamma_\alpha}}\eta(\left(x,m\right), \Xi,\Gamma_\alpha)\mathbf{1}_{\bar{R}(x,\alpha)\le r}}$, then
		\begin{align}
		d_{TV}(\bar{W}_{\alpha,r},\bar{W}_{\alpha,r}+\gamma)&\le C(|\gamma|\vee 1)\left(\alpha^{-\frac{1}{2}}r^{\frac{d}{2}}\right)\label{statement3}		\end{align}
	 for any $\gamma\in \mathbb{R}$ and $r>R_0$,
		where $C$ and $R_0$ are positive constants independent of $\gamma$. If $\eta$ is range-bound in Definition~\ref{defi4r}, then
		\begin{align}d_{TV}(\bar{W}_\alpha,\bar{W}_\alpha+\gamma)\le C(|\gamma|\vee 1)\alpha^{-\frac{1}{2}}\label{statement4}\end{align} for some positive constant $C$ independent of $\gamma$.		
		\end{description}
	\end{lma} 

\begin{re}{\rm Since exponential stabilization implies polynomial stabilization, the statements {\Ref{statement1} and \Ref{statement3}} \ignore{of Lemma~\ref{lma3} }also hold under corresponding exponential stabilization conditions.}
\end{re}

Now, a more general version of Lemma~\ref{lma3} with $\gamma$ replaced by a function of $\Xi_N$ for some Borel set $N$ and the expectation replaced by a conditional expectation can be easily established.

	\begin{cor}\label{cor1} For $\alpha,r>0$, let $\{N_{\alpha,r}^{(k)}\}_{k\in\{1,2,3\}}\subset\mathscr{B}(\mathbb{R}^d)$  such that $\(N_{\alpha,r}^{(1)}\cup N_{\alpha,r}^{(2)}\cup N_{\alpha,r}^{(3)}\)\cap \Gamma_\alpha\in B(x, C\alpha^{\frac{1}{d}})$ for a point $x\in \mathbb{R}^d$ and a positive constant $C\in (0,\frac{1}{2})$, $\mathscr{F}_{0,\alpha,r}$ be a sub $\sigma$-algebra of $\sigma\(\Xi_{N_{\alpha,r}^{(1)}}\)$, and $h_{\alpha,r}$ be a measurable function mapping configurations on $N_{\alpha, r}^{(2)}\times T$ to the real space.
	\begin{description}
	\item{(a)} (unrestricted case)  Define $W_{\alpha,r}':=\sum_{(x,m)\in \Xi_{\Gamma_\alpha\backslash N_{\alpha,r}^{(3)}}}\eta(\left(x,m\right), \Xi)\mathbf{1}_{R(x)\le r}$. If the conditions of Lemma~\ref{lma3}~(a) hold, then  
		\begin{align}
        &d_{TV}\left(W_{\alpha,r}',W_{\alpha,r}'+h_{\alpha,r}\left(\Xi_{N_{\alpha, r}^{(2)}}\right)\middle|\mathscr{F}_{0,\alpha,r}\right)
        \nonumber\\
        &\le \mathbb{E}\left(\left|h_{\alpha,r}\left(\Xi_{N_{\alpha, r}^{(2)}}\right)\right|\vee 1\middle|\mathscr{F}_{0,\alpha,r} \right)O\left(\alpha^{-\frac{1}{2}}r^{\frac{d}{2}}\right)~a.s.,\label{lma3coro01}
		\end{align}
		where $O\left(\alpha^{-\frac{1}{2}}r^{\frac{d}{2}}\right)$ is independent of sets $\{N_{\alpha,r}^{(k)}\}_{\alpha,r\in \mathbb{R}_+,k\in\{1,2,3\}}$, functions $\{h_{\alpha,r}\}_{\alpha,r\in \mathbb{R}_+}$ and $\sigma$-algebras $\{\mathscr{F}_{0,\alpha,r}\}_{\alpha,r\in \mathbb{R}_+}$.
		
	\item{(b)}	(restricted case) Define $\bar{W}_{\alpha,r}':=\sum_{(x,m)\in \Xi_{\Gamma_\alpha\backslash N_{\alpha,r}^{(3)}}}\eta(\left(x,m\right), \Xi,{\Gamma_\alpha})\mathbf{1}_{\bar{R}(x,\alpha)<r}\overline{\Xi}(dx)$. If the conditions of Lemma~\ref{lma3}~(b) hold, then
		\begin{align}
&d_{TV}\left(\bar{W}_{\alpha,r}',\bar{W}_{\alpha,r}'+h_{\alpha,r}\left(\Xi_{N_{\alpha, r}^{(2)} }\right)\middle|\mathscr{F}_{0,\alpha,r}\right)\nonumber\\
&\le \mathbb{E}\left(\left|h_{\alpha,r}\left(\Xi_{N_{\alpha, r}^{(2)} }\right)\right|\vee 1\middle|\mathscr{F}_{0,\alpha,r} \right)O\left(\alpha^{-\frac{1}{2}}r^{\frac{d}{2}}\right)~a.s.,\label{lma3coro02}
		\end{align}
		where $O\left(\alpha^{-\frac{1}{2}}r^{\frac{d}{2}}\right)$ is independent of sets $\{N_{\alpha,r}^{(k)}\}_{\alpha,r\in \mathbb{R}_+,k\in\{1,2,3\}}$, functions $\{h_{\alpha,r}\}_{\alpha,~r\in \mathbb{R}_+}$ and $\sigma$-algebras $\{\mathscr{F}_{0,\alpha,r}\}_{\alpha,r\in \mathbb{R}_+}$.
		\end{description}
\end{cor}

As discussed in the inspiring example, the order of $\var\(W_{\alpha}\)$ and $\var\(\bar{W}_{\alpha}\)$ plays the crucial role in the accuracy of normal approximation. The next lemma says that under exponential stabilization, optimal order of the variances can be achieved.

\begin{lma}\label{lma4}
\begin{description}
\item{(a)} (unrestricted case) If the score function $\eta$ satisfies the third moment condition~\Ref{thm2.1}, non-singularity~\Ref{non-sin} and exponential stabilization~in Definition~\ref{defi4}, then
	  $\var\(W_{\alpha}\)=\Omega\(\alpha\)$.
\item{(b)} (restricted case) If the score function $\eta$ satisfies the third moment condition~\Ref{thm2.1r}, non-singularity~\Ref{non-sinr} and exponential stabilization~in Definition~\ref{defi4r}, then
$\var\(\bar{W}_{\alpha}\)$ $=\Omega\(\alpha\)$.
	 \end{description}
\end{lma}

For polynomially stabilizing score functions, we do not know the optimal order of the variance, but we can get a lower bound as shown in the next lemma.

\begin{lma}\label{lma5}
\begin{description}
\item{(a)} (unrestricted case) If the score function $\eta$ satisfies the $k'$-th moment condition~\Ref{thm2.1} with $k'>k\ge3$, non-singularity~\Ref{non-sin} and is polynomially stabilizing in Definition~\ref{defi4} with parameter $\beta>(3k-2)d/(k-2)$, 
then
$$\var\(W_{\alpha}\)\ge O\(\alpha^{\frac{k\beta-2\beta-3dk+2d}{k\beta-2\beta-dk}}\).$$

\item{(b)} (restricted case) If the score function $\eta$ satisfies the $k'$-th moment condition~\Ref{thm2.1r} with $k'>k\ge3$, non-singularity~\Ref{non-sinr} and polynomially stabilizing in Definition~\ref{defi4r} with parameter $\beta>(3k-2)d/(k-2)$, then
$$\var\(\bar{W}_{\alpha}\)\ge O\(\alpha^{\frac{k\beta-2\beta-3dk+2d}{k\beta-2\beta-dk}}\).$$
	 \end{description}
\end{lma} 

\section{The proofs of the auxiliary and main results}\label{Theproofs}

We need Palm processes and reduced Palm processes as the tools in our proofs, and for the ease of reading, we briefly recall their definitions. Let $H$ be a Polish space with Borel $\sigma-$algebra $\mathscr{B}\(H\)$ and configuration space $\(\bm{C}_H, \mathscr{C}_H\)$, let $\Psi$ be a point process on $\(\bm{C}_H, \mathscr{C}_H\)$, write the mean measure of $\Psi$ as $\psi (d x) := \mathbb{E} \Psi (d x)$, the family of point processes $\{ \Psi_x : x \in H\}$ are said to be the {\it Palm processes} associated with $\Psi$ if for any measurable function $f : H\times \bm{C}_H \rightarrow [0,\infty)$,
\begin{equation}\label{palm1}
	\mathbb{E} \[ \int_{H} f(x,\Psi)\Psi(dx) \] = \int_{H} \mathbb{E} f(x,\Psi_x)  \psi(dx) ,
\end{equation}
\cite[\S~10.1]{Kallenberg83}. A Palm process $\Psi_x$ contains a point at $x$ and it is often more convenient to consider the {\it reduced Palm process} $\Psi_x-\delta_x$ at $x$ by removing the point $x$ from $\Psi_x$. Furthermore, suppose that the {\it factorial moments} $\psi^{[2]}(dx,dy) : = \mathbb{E} [\Psi(dx)(\Psi-\d_x)(dy)]$ and $\psi^{[3]}(dx,dy,dz) : = \mathbb{E} [\Psi(dx)(\Psi-\d_x)(dy)(\Psi-\d_x-\delta_y)(dz)]$ are finite, then we can respectively define the {\it second order Palm processes} $\{ \Psi_{xy}: x,y \in H \}$ and {\it third order Palm processes} $\{ \Psi_{xyz}: x,y,z \in H \}$ associated with $\Psi$ by
\begin{align}
	&\mathbb{E} \[ \iint_{H^2}
	f(x,y;\Psi)\Psi(dx)(\Psi-\d_x)(dy) \]= \iint_{H^2} \mathbb{E}
	f(x,y;\Psi_{xy}) \psi^{[2]} (dx,dy) ,\label{palm2}\\
	&\mathbb{E} \[ \iiint_{H^3}
	f(x,y,z;\Psi)\Psi(dx)(\Psi-\d_x)(dy)(\Psi-\d_x-\d_y)(dz) \] \nonumber\\& = \iiint_{H^3} \mathbb{E}
	f(x,y,z;\Psi_{xyz}) \psi^{[3]} (dx,dy,dz) ,\label{palm3}
\end{align}
for all measurable functions $f : H^2 \times \bm{C}_H \rightarrow [0,\infty)$ in \Ref{palm2} and $f: H^3 \times \bm{C}_H \rightarrow [0,\infty)$ in \Ref{palm3} \cite[\S~12.3]{Kallenberg83}. Using reduced Palm processes, Slivnyak-Mecke theorem \cite{Mecke63} states that a point process such that the distributions of its reduced Palm processes are the same as that of the point process if and only if it is a Poisson point process. Then we can see that for homogeneous Poisson point process with rate $\lambda$, its mean measure can be written as $\Lambda(dx)=\lambda dx$ and its Palm processes satisfy $\Psi_{x}\overset{d}{=}\Psi+\d_x$, $\Psi_{xy}\overset{d}{=}\Psi+\d_x+\d_y$ and $\Psi_{xyz}\overset{d}{=}\Psi+\d_x+\d_y+\d_z$, the factorial moments $\psi^{[2]}(dx,dy)=\lambda^2dxdy$ and $\psi^{[3]}(dx,dy,dz)=\lambda^3dxdydz$ for all distinct $x$, $y$, $z\in H$.  

We can adapt \Ref{palm1}, \Ref{palm2} and \Ref{palm3} to the marked homogeneous Poisson point process case that we are dealing with, assume the rate of $\overline{\Xi}$ is $\lambda$, $\{M_i\}_{i\in \mathbb{N}}$ is a sequence of $\iid$ random elements on $(T,\mathscr{T})$ following the distribution $\mathscr{L}_{T}$ which is independent of $\Xi$, then because of the independence of marks, we can see that
\begin{align}
	\label{palm4}
	&\mathbb{E} \[ \int_{\bm{S}} f(x,\Xi)\overline{\Xi}(dx) \] = \int_{\bm{S}}  \mathbb{E} f(x,\Xi+\d_{(x,M_1)})  \lambda dx, \\
	&\mathbb{E} \[ \iint_{\bm{S}^2}
	f(x,y;\Xi)\overline{\Xi}(dx)(\overline{\Xi}-\d_x)(dy) \] = \iint_{\bm{S}^2} \mathbb{E}
	f(x,y;\Xi+\d_{(x,M_1)}+\d_{(y,M_2)}) \lambda^2dxdy ,\label{palm5}\\
	&\mathbb{E} \[ \iiint_{\bm{S}^3}
	f(x,y,z;\Xi)\overline{\Xi}(dx)(\overline{\Xi}-\d_x)(dy)(\overline{\Xi}-\d_x-\d_y)(dz) \]\nonumber
	\\ =& \iiint_{\bm{S}^3} \mathbb{E}
	f(x,y,z;\Xi+\d_{(x,M_1)}+\d_{(y,M_2)}+\d_{(z,M_3)}) \lambda^3dxdydz,\label{palm6}
\end{align}
for all measurable functions $f :\bm{S} \times (\bm{C}_{\bm{S}}\times T) \rightarrow [0,\infty)$ in \Ref{palm4}, $f:\bm{S}^2 \times (\bm{C}_{\bm{S}}\times T) \rightarrow [0,\infty)$ in \Ref{palm5}  and $f:\bm{S}^3 \times (\bm{C}_{\bm{S}}\times T) \rightarrow [0,\infty)$ in \Ref{palm6}.

Recalling the shift operator defined in Section~\ref{Generalresults}, we can write $g(\mathscr{X}^x):=\eta(\left(x,m\right), \mathscr{X})$ (resp. $g_\alpha(x,\mathscr{X}):= \eta\(\(x,m\),\Xi,\mathscr{X},\Gamma_\alpha\)$) for all configuration $\mathscr{X}$, $(x,m)\in \mathscr{X}$ and $\alpha>0$ so that notations can be simplified significantly, e.g, . 
\begin{eqnarray*}
		W_\alpha&=&\sum_{(x,m)\in\Xi_{\Gamma_\alpha}}\eta(\left(x,m\right), \Xi)=\int_{\Gamma_\alpha} g(\Xi^x)\overline{\Xi}(dx)=\sum_{x\in \bar{\Xi}}g(\Xi^x),\\
		\bar{W}_\alpha&=&\sum_{(x,m)\in\Xi_{\Gamma_\alpha}}\eta(\left(x,m\right), \Xi,\Gamma_\alpha)=\int_{\Gamma_\alpha}{g_\alpha(x, {\Xi})\overline{\Xi}(dx)},
\end{eqnarray*}
where $\overline{\Xi}$ is the projection of $\Xi$ on $\mathbb{R}^d$, $R$ and $\bar{R}(x,\alpha)$ are the corresponding radii of stabilization.

We now proceed to establish a few lemmas needed in the proofs.

\begin{lma}\label{lma6} (Conditional Total Variance Formula)	Let $X$ be a random variable on probability space $(\Omega, \mathscr{G},\mathbb{P})$ with finite second moment, $\mathscr{G}_1$ and $\mathscr{G}_2$ be two sub $\sigma$-algebras of $\mathscr{G}$ such that $\mathscr{G}_1\subset \mathscr{G}_2$, then 
	$$\var\(X|\mathscr{G}_1\)=\mathbb{E}\left(\var\left(X\middle|\mathscr{G}_2\right)\middle|\mathscr{G}_1\right)+\var\left(\mathbb{E}\left(X\middle|\mathscr{G}_2\right)\middle|\mathscr{G}_1\right).$$
\end{lma}	
\noindent{\it Proof.} From the definition of the conditional variance, we can see that 
\begin{align*}
	\var\(X|\mathscr{G}_1\)=&\mathbb{E}\left(X^2\middle|\mathscr{G}_1\right)-\mathbb{E}\left(X\middle|\mathscr{G}_1\right)^2
	\\=&\mathbb{E}\left(\mathbb{E}\left(X^2\middle|\mathscr{G}_2\right)\middle|\mathscr{G}_1\right)-\mathbb{E}\left(\mathbb{E}\left(X\middle|\mathscr{G}_2\right)\middle|\mathscr{G}_1\right)^2
	\\=&\mathbb{E}\left(\mathbb{E}\left(X^2\middle|\mathscr{G}_2\right)\middle|\mathscr{G}_1\right)-\mathbb{E}\left(\mathbb{E}\left(X\middle|\mathscr{G}_2\right)^2\middle|\mathscr{G}_1\right)+\mathbb{E}\left(\mathbb{E}\left(X\middle|\mathscr{G}_2\right)^2\middle|\mathscr{G}_1\right)\\
	&-\mathbb{E}\left(\mathbb{E}\left(X\middle|\mathscr{G}_2\right)\middle|\mathscr{G}_1\right)^2
	\\=&\mathbb{E}\left(\var\left(X\middle|\mathscr{G}_2\right)\middle|\mathscr{G}_1\right)+\var\left(\mathbb{E}\left(X\middle|\mathscr{G}_2\right)\middle|\mathscr{G}_1\right),
\end{align*} so the statement holds. \qed

Also, given the value of a random variable in a certain event, we can find a lower bound for the conditional variance.
\begin{lma}\label{lma7}
	Let $X$ be a random variable on probability space $(\Omega, \mathscr{G},\mathbb{P})$ with finite second moment, for any event $A$ and $\sigma$-algebra such that $\mathscr{F}\subset\mathscr{G}$, 
	\begin{equation}\label{lma7.1}
		\var\left(X\middle|\mathscr{F}\right)\ge \var\left(X\mathbf{1}_A+\frac{\mathbb{E}\left(X\mathbf{1}_A\middle|\mathscr{F}\right)}{\mathbb{P}\left(A\middle|\mathscr{F}\right)}\mathbf{1}_{A^c}\middle|\mathscr{F}\right),
	\end{equation}
	where $\frac{0}{0}=0$ by convention.
\end{lma}
\noindent{\it Proof of Lemma~\ref{lma7}.} For the case that event $A$ has probability $0$, the statement is trivially true, so we focus on the case that $\mathbb{P}(A)>0$. Let $A\cap\mathscr{F}=\{B\cap A;~B\in \mathscr{F}\}$, which is a $\sigma-$algebra on $A$, and $\mathbb{P}_A$ as a probability measure on $(A, A\cap\mathscr{F})$ such that $\mathbb{P}_A\(B\cap A\)=\mathbb{P}\left(B\middle|A\right)$ for all $B\in \mathscr{F}$, then we have the corresponding conditional expectation $$\mathbb{E}_A\left(X\middle|A\cap\mathscr{F}\right)=\mathbb{E}_A\left(X\mathbf{1}_A\middle|A\cap\mathscr{F}\right)=\mathbf{1}_A\mathbb{E}_A\left(X\middle|A\cap\mathscr{F}\right),$$ which equals to $0$ on $\mathbf{1}_{A^c}$. 
The proof relies on the following observations:\\
\begin{lma} \label{lma8}
	For any random variable $Y$ with $\mathbb{E}\vert Y\vert < \infty$ and $A\in \mathscr{G}$ such that $\mathbb{P}(A)>0$,
	$$
		\mathbb{E}\left(\mathbb{E}_A\left(Y\middle|A\cap\mathscr{F}\right)\middle|\mathscr{F}\right)=\mathbb{E}\left(\mathbf{1}_A Y\middle|\mathscr{F}\right).
	$$
\end{lma}
\noindent{\it Proof.} Both sides are $\mathscr{F}$ measurable, together with the fact that for any $B\in \mathscr{F}$, 
\begin{align}
	\mathbb{E}\left(\mathbb{E}\left(\mathbb{E}_A\left(Y\middle|A\cap\mathscr{F}\right)\middle|\mathscr{F}\right)\mathbf{1}_B\right)=&\mathbb{E}\left(\mathbb{E}_A\left(Y\middle|A\cap\mathscr{F}\right)\mathbf{1}_B\right)\nonumber 
	\\=&\mathbb{E}\left(\mathbb{E}_A\left(Y\middle|A\cap\mathscr{F}\right)\mathbf{1}_{B\cap A}\right)\nonumber
	\\=&\mathbb{E}\left(\mathbb{E}_A\left(\mathbf{1}_{B\cap A}Y\middle|A\cap\mathscr{F}\right)\right)\nonumber
	\\=&\mathbb{P}(A)\mathbb{E}_A\left(\mathbb{E}_A\left(\mathbf{1}_{B\cap A}Y\middle|A\cap\mathscr{F}\right)\right)\nonumber
	\\=&\mathbb{P}(A)\mathbb{E}_A\left(\mathbf{1}_{B\cap A}Y\right)=\mathbb{E}\left(\mathbf{1}_{B\cap A}Y\right)=\mathbb{E}\left(\mathbb{E}\left(\mathbf{1}_A Y\middle|\mathscr{F}\right)\mathbf{1}_B\right),\nonumber
\end{align} as claimed. \qed
\begin{lma} \label{lma9}
	For any random variable $Y$ with $\mathbb{E}\vert Y\vert <\infty$ and $A\in \mathscr{G}$ such that $\mathbb{P}(A)>0$,
	\begin{equation}\label{lma9.1}
		\frac{\mathbb{E}\left(\mathbf{1}_A Y \middle|\mathscr{F}\right)}{\mathbb{P}\left(A\middle|\mathscr{F} \right)}\mathbf{1}_A=\mathbb{E}_A\left(Y\middle|\mathscr{F}\right).
	\end{equation}
\end{lma}
\noindent{\it Proof.} Both sides equal $0$ on $A^c$, and are measurable on $A\cap\mathscr{F}$ when restricted to $A$. From the construction of $ A\cap\mathscr{F}$, any set $B'\in  A\cap\mathscr{F}$ is of the form $B\cap A$ for some $B\in\mathscr{F}$, hence \Ref{lma9.1} is equivalent to
$$\mathbb{E}_A\left(\frac{\mathbb{E}\left(\mathbf{1}_A Y \middle|\mathscr{F}\right)}{\mathbb{P}\left(A\middle|\mathscr{F} \right)}\mathbf{1}_A\mathbf{1}_{A\cap B}\right)=\mathbb{E}_A\left(\mathbb{E}_A\left(Y\middle|\mathscr{F}\right)\mathbf{1}_{A\cap B}\right)$$ for all $B\in \mathscr{F}$. 
Now we have 
\begin{align}
	\mathbb{E}_A\left(\frac{\mathbb{E}\left(\mathbf{1}_A Y \middle|\mathscr{F}\right)}{\mathbb{E}\left(\mathbf{1}_A\middle|\mathscr{F} \right)}\mathbf{1}_{A\cap B}\right)=&\frac{1}{\mathbb{P}(A)}\mathbb{E}\left(\frac{\mathbb{E}\left(\mathbf{1}_A Y \middle|\mathscr{F}\right)}{\mathbb{E}\left(\mathbf{1}_A\middle|\mathscr{F} \right)}\mathbf{1}_A\mathbf{1}_B\right)\nonumber
	\\=&\frac{1}{\mathbb{P}(A)}\mathbb{E}\left(\mathbb{E}\left(\frac{\mathbb{E}\left(\mathbf{1}_A Y \middle|\mathscr{F}\right)}{\mathbb{E}\left(\mathbf{1}_A\middle|\mathscr{F} \right)}\mathbf{1}_A\mathbf{1}_B\middle|\mathscr{F}\right)\right)\nonumber
	\\=&\frac{1}{\mathbb{P}(A)}\mathbb{E}\left(\mathbb{E}\left(\mathbf{1}_A Y \middle|\mathscr{F}\right)\mathbf{1}_B\right)=\frac{1}{\mathbb{P}(A)}\mathbb{E}\left(\mathbb{E}\left(\mathbf{1}_A\mathbf{1}_B Y \middle|\mathscr{F}\right)\right)\nonumber
	\\=&\frac{1}{\mathbb{P}(A)}\mathbb{E}\left(\mathbf{1}_A\mathbf{1}_B Y\right)=\mathbb{E}_A(\mathbb{E}_A(Y\mathbf{1}_A\mathbf{1}_B|\mathscr{F}))=\mathbb{E}_A(\mathbb{E}_A(Y|\mathscr{F})\mathbf{1}_{A\cap B})\nonumber
\end{align} completing the proof. \qed\\
\noindent{\it Proof of Lemma~\ref{lma7} (continued).} We start from the left hand side of \Ref{lma7.1},
\begin{align}
	&\var\left(X\middle|\mathscr{F}\right)\nonumber
	\\=&\mathbb{E}\left(X^2\middle|\mathscr{F}\right)-\mathbb{E}\left(X\middle|\mathscr{F}\right)^2\nonumber
	\\=&\mathbb{E}\left(X^2\middle|\mathscr{F}\right)-\left(\mathbb{E}\left(X\mathbf{1}_A\middle|\mathscr{F}\right)+\mathbb{E}\left(X\mathbf{1}_{A^c}\middle|\mathscr{F}\right)\right)^2\nonumber
	\\=&\mathbb{E}\left(X^2\middle|\mathscr{F}\right)-\mathbb{E}\left(\mathbb{E}_A\left(X\mathbf{1}_A\middle|A\cap\mathscr{F}\right)+\mathbb{E}_{A^c}\left(X\middle|A^c\cap\mathscr{F}\right)\middle|\mathscr{F}\right)^2\nonumber
	\\\ge&\mathbb{E}\left(X^2\middle|\mathscr{F}\right)-\mathbb{E}\left(\left(\mathbb{E}_A\left(X\mathbf{1}_A\middle|A\cap\mathscr{F}\right)+\mathbb{E}_{A^c}\left(X\middle|A^c\cap\mathscr{F}\right)\right)^2\middle|\mathscr{F}\right)\nonumber
	\\=&\mathbb{E}\left(X^2\(\mathbf{1}_A+\mathbf{1}_{A^c}\)\middle|\mathscr{F}\right)-\mathbb{E}\left(\mathbb{E}_A\left(X\mathbf{1}_A\middle|A\cap\mathscr{F}\right)^2+\mathbb{E}_{A^c}\left(X\middle|A^c\cap\mathscr{F}\right)^2\middle|\mathscr{F}\right)\nonumber
	\\=&\mathbb{E}\left(\mathbb{E}_A\left(X^2\middle| A\cap\mathscr{F}\right)+\mathbb{E}_{A^c}\left(X^2\middle| A^c\cap\mathscr{F}\right)\middle|\mathscr{F}\right)\nonumber\\
	&-\mathbb{E}\left(\mathbb{E}_A\left(X\mathbf{1}_A\middle|A\cap\mathscr{F}\right)^2+\mathbb{E}_{A^c}\left(X\middle|A^c\cap\mathscr{F}\right)^2\middle|\mathscr{F}\right)\nonumber
	\\=&\mathbb{E}\left(\var_A\left(X\middle| A\cap\mathscr{F}\right)+\var_{A^c}\left(X\middle| A^c\cap\mathscr{F}\right)\middle|\mathscr{F}\right)\label{lma7.2},
\end{align}
where the inequality follows from Jensen's inequality, the third and the second last equalities are from Lemma~\ref{lma8}.
On the other hand, the right hand side of \Ref{lma7.1} can be written as 
\begin{align}
	&\mathbb{E}\left(X^2\mathbf{1}_A+\left(\frac{\mathbb{E}\left(X\mathbf{1}_A\middle|\mathscr{F}\right)}{\mathbb{P}\left(A\middle|\mathscr{F}\right)}\right)^2\mathbf{1}_{A^c}\middle|\mathscr{F}\right)-\left(\mathbb{E}\left(X\mathbf{1}_A\middle|\mathscr{F}\right)+\mathbb{E}\left(\frac{\mathbb{E}\left(X\mathbf{1}_A\middle|\mathscr{F}\right)}{\mathbb{P}\left(A\middle|\mathscr{F}\right)}\mathbf{1}_{A^c}\middle|\mathscr{F}\right)\right)^2\nonumber
	\\=&\mathbb{E}\left(X^2\mathbf{1}_A\middle|\mathscr{F}\right)+\left(\frac{\mathbb{E}\left(X\mathbf{1}_A\middle|\mathscr{F}\right)}{\mathbb{P}\left(A\middle|\mathscr{F}\right)}\right)^2\mathbb{P}\left(A^c\middle|\mathscr{F}\right)-\mathbb{E}\left(X\mathbf{1}_A\middle|\mathscr{F}\right)^2\frac{1}{\mathbb{P}\left(A\middle|\mathscr{F}\right)^2}\nonumber
	\\=&\mathbb{E}\left(X^2\mathbf{1}_A\middle|\mathscr{F}\right)-\left(\frac{\mathbb{E}\left(X\mathbf{1}_A\middle|\mathscr{F}\right)}{\mathbb{P}\left(A\middle|\mathscr{F}\right)}\right)^2\mathbb{P}\left(A\middle|\mathscr{F}\right)\nonumber
	\\=&\mathbb{E}\left(\mathbb{E}_A\left(X^2\middle|A\cap \mathscr{F}\right)-\left(\frac{\mathbb{E}\left(X\mathbf{1}_A\middle|\mathscr{F}\right)}{\mathbb{P}\left(A\middle|\mathscr{F}\right)}\right)^2\mathbf{1}_A\middle|\mathscr{F}\right)\nonumber
	\\=&\mathbb{E}\left(\mathbb{E}_A\left(X^2\middle|A\cap \mathscr{F}\right)-\mathbb{E}_A\left(X\middle|A\cap \mathscr{F}\right)^2\middle|\mathscr{F}\right)\nonumber
	\\=&\mathbb{E}\left(\var_A\left(X\middle| A\cap\mathscr{F}\right)\middle|\mathscr{F}\right)\label{lma7.3},
\end{align}where the third equality follows from Lemma~\ref{lma8}, and the second last equality is from Lemma~\ref{lma9}. Combining \Ref{lma7.2}, \Ref{lma7.3} and Lemma~\ref{lma6} completes the proof. \qed

The following lemma bounds the difference between two normal distribution under the total variation distance.
\begin{lma}\label{lma10}
	Let $F_{\mu,\sigma}$ be the distribution of $N(\mu,\sigma^2)$, the normal distribution with mean $\mu$ and variance $\sigma^2$, then 
	$$d_{TV}(F_{\mu_1,\sigma_1},F_{\mu_2,\sigma_2})\le \sqrt{\frac{2}{\pi}}\(\frac{|\mu_1-\mu_2|}{2\max(\sigma_1,\sigma_2)}+\frac{\max(\sigma_1,\sigma_2)}{\min(\sigma_1,\sigma_2)}-1\).$$
\end{lma}
\noindent{\it Proof.} Without loss of generality, we assume $\sigma_2>\sigma_1$. Writing the probability density function of $F_{\mu,\sigma}$ as $f_{\mu,\sigma}$, we have
\begin{equation}\label{9.1}
	\begin{aligned}
		d_{TV}(F_{\mu_1,\sigma_1},F_{\mu_2,\sigma_2})&\le  d_{TV}(F_{\mu_1,\sigma_1},F_{\mu_1,\sigma_2})+d_{TV}(F_{\mu_1,\sigma_2},F_{\mu_2,\sigma_2})\\&=d_{TV}\(F_{0,1},F_{0,\frac{\sigma_2}{\sigma_1}}\)+d_{TV}\(F_{0,1},F_{\frac{\mu_1-\mu_2}{\sigma_2},1}\).
	\end{aligned}
\end{equation} 
Then the problem turns to bound the differences between the distributions of $N(0,1)$ and $N(0,\sigma^2)$ and between the distributions of $N(0,1)$ and $N(\mu,1)$ for $\mu>0$ and $\sigma>1$.  For $\sigma>1$, we can see that the probability density functions $f_{0,1}$ and $f_{0,\sigma}$ meet at $\pm x_\sigma:=\pm\sqrt{\frac{2\ln(\sigma)\sigma^2}{\sigma^2-1}}$, and $f_{0,1}>f_{0,\sigma}$ on $(-x_\sigma,x_\sigma)$ and the inequality sign is reversed outside the interval. We can see that $1<x_\sigma<\sigma$, so we have
\begin{equation}
	\begin{aligned}\label{9.3}
		d_{TV}(F_{0,1},F_{0,\sigma})
		&=F_{0,1}(x_\sigma)-F_{0,1}(-x_\sigma)-(F_{0,\sigma}(x_\sigma)-F_{0,\sigma}(-x_\sigma))\\&=F_{0,1}(x_\sigma)-F_{0,1}\(\frac{x_\sigma}{\sigma}\)-F_{0,1}(-x_\sigma)+F_{0,1}\(-\frac{x_\sigma}{\sigma}\)\le \sqrt{\frac{2}{\pi}}(\sigma-1),
	\end{aligned}
\end{equation} where the inequality follows from the fact that the probability density function $f_{0,1}$ is bounded by $\frac{1}{\sqrt{2\pi}}$. Similarly, we can see that  
\begin{equation}
	\begin{aligned}\label{9.2}
		d_{TV}(F_{0,1},F_{\mu,1})
		&=F_{0,1}\(\frac{\mu}{2}\)-F_{\mu,1}\(\frac{\mu}{2}\)=F_{0,1}\(\frac{\mu}{2}\)-F_{0,1}\(-\frac{\mu}{2}\)\le \frac{\mu}{\sqrt{2\pi}},
	\end{aligned}
\end{equation} for $\mu>0$, where again we use the fact that $f_{0,1}$ is bounded by $\frac{1}{\sqrt{2\pi}}$ in the inequality. Substituting \Ref{9.2} and \Ref{9.3} into \Ref{9.1} yields the claim. \qed

The following lemma says that under stabilizing conditions, the cost of throwing away the terms with large radii of stabilization is negligible.

\begin{lma}\label{lma105} 
	\begin{description}
		\item{(a)} (unrestricted case) If the score function is exponentially stabilizing in Definition~\ref{defi4}, then we have
		$$d_{TV}(W_\alpha, W_{\alpha,r})\le C_1\alpha e^{-C_2r}$$
		for some positive constants $C_1$, $C_2$.
		If the score function is polynomially stabilizing with parameter $\beta$ in Definition~\ref{defi4}, then we have 
		$$d_{TV}(W_\alpha, W_{\alpha,r})\le C\alpha r^{-\beta}$$ for some positive constant $C$.
		
		\item{(b)} (restricted case) If the score function is exponentially stabilizing in Definition~\ref{defi4r}, then we have
		$$d_{TV}(\bar{W}_\alpha, \bar{W}_{\alpha,r})\le C_1\alpha e^{-C_2r}$$
		for some positive constants $C_1$, $C_2$.
		If the score function is polynomially stabilizing with parameter $\beta$ in Definition~\ref{defi4r}, then we have 
		$$d_{TV}(\bar{W}_\alpha, \bar{W}_{\alpha,r})\le C\alpha r^{-\beta}$$ for some positive constant $C$.
	\end{description}
\end{lma}
\noindent{\it Proof.} We first show the statement is true for $\bar{W}_\alpha$ and $\bar{W}_{\alpha,r}$. For convenience of writing, we define $M_x$ as random elements following the law $\mathscr{L}_T$ which are independent of $\Xi$ for all $x\in \mathbb{R}^d$. From the construction of $\bar{W}_\alpha$ and $\bar{W}_{\alpha,r}$, we can see that the event $\{\bar{W}_\alpha\neq \bar{W}_{\alpha,r}\}\subset\{\mbox{at least one }x\in \overline{\Xi}\cap\Gamma_\alpha~\mbox{with }\bar{R}(x, \alpha)>r\}$, so from \Ref{palm4}, we have 
\begin{align*}
	d_{TV}(\bar{W}_\alpha, \bar{W}_{\alpha,r})\le &\mathbb{P}\(\{\bar{W}_\alpha\neq \bar{W}_{\alpha,r}\}\)
	\\\le & \mathbb{P}\(\{\mbox{at least one }x\in \overline{\Xi}\cap\Gamma_\alpha~{\mbox{such that }}\bar{R}(x, \alpha)>r\}\)
	\\ \le& \mathbb{E}\int_{\Gamma_\alpha}\mathbf{1}_{\bar{R}(x,\alpha)>r} \overline{\Xi}(dx)
	\\ =&\int_{\Gamma_\alpha}\mathbb{E}\(\mathbf{1}_{\bar{R}(x,M_x,\alpha, \Xi+\delta_{(x, M_x)})>r}\)\lambda dx
	\\ =&\int_{\Gamma_\alpha}\mathbb{P}\(\bar{R}(x,M_x,\alpha, \Xi+\delta_{(x, M_x)})>r\) \lambda dx
	\\ \le&\alpha\lambda\bar{\tau}(r),
\end{align*} which, together with the stabilization conditions, gives the claim for $\bar{W}_{\alpha}$.

The statement is also true for $W_\alpha$, which can be proved by replacing corresponding counterparts $\bar{W}_\alpha$ with $W_\alpha$; $\bar{W}_{\alpha,r}$ with $W_{\alpha,r}$; $\bar{R}(x,\alpha)$ with $R(x)$; $\bar{R}(x,M_x,\alpha,\Xi+\delta_{(x, M_x)})$ with $R(x,M_x,\Xi+\delta_{(x, M_x)})$; $\bar{\tau}$ with $\tau.$\qed

\noindent{\it Proof of Lemma~\ref{lma1}.} For convenience, we write $G_n$, $g_n$ and $\psi_n$ as the distribution, density and characteristic functions of $T_n$ respectively. It is well-known that the triangular density $\kappa_a$ has the characteristic function $\psi_1(s)=\frac{2(1-\cos(as))}{(as)^2}$, which gives $\psi_n(s)=\left(\frac{2(1-\cos(as))}{(as)^2}\right)^n$. Using the fact that the convolution of two symmetric unimodal distributions on $\real$ is unimodal \cite{Wintner38}, we can conclude that the distribution of $T_n$ is unimodal and symmetric. This ensures that 
\begin{equation}d_{TV}(T_n,T_n+\gamma)=\sup_{x\in\real}|G_n(x)-G_n(x-\gamma)|=\int_{-\gamma/2}^{\gamma/2}g_n(x)dx.\label{lma1.3}\end{equation}
Applying the inversion formula, we have
\begin{eqnarray*}
g_n(x)&=&\frac1{2\pi}\int_\real e^{-\ime sx}\psi_n(s)ds=\frac1{2\pi}\int_\real \cos(sx)\psi_n(s)ds\nonumber\\
&=&\frac1{a\pi}\int_0^\infty \cos(sx/a)\left(\frac{2(1-\cos s)}{s^2}\right)^nds,
\end{eqnarray*}
where $\im=\sqrt{-1}$ and the second equality is due to the fact that $\sin(sx)\psi_n(s)$ is an odd function. Obviously, $g_n(x)\le g_n(0)$ so we need to establish an upper bound for $g_n(0)$. A direct verification gives
$$0\le \frac{2(1-\cos s)}{s^2}\le e^{-\frac{s^2}{12}}\mbox{ for }0\le s\le 2\pi,$$
which implies
\begin{eqnarray}
g_n(0)&\le&\frac1{a\pi} \left\{\int_0^{2\pi}e^{-\frac{ns^2}{12}}ds+\int_{2\pi}^\infty \left(\frac 4{s^2}\right)^nds\right\}\nonumber\\
&\le& \frac1{a\pi\sqrt{n}}\int_0^\infty e^{-\frac{s^2}{12}}ds+\frac{2}{a(2n-1)\pi^{2n}}\nonumber\\
&=&\frac1a\sqrt{\frac{3}{\pi n}}+\frac{2}{a(2n-1)\pi^{2n}}.\label{lma1.5}
\end{eqnarray}
Now, combining \Ref{lma1.5} with \Ref{lma1.3} gives \Ref{lma1.2}.\qed 

\noindent{\it Proof of Lemma~\ref{non-singular1}.} We construct a maximal coupling \cite[p.~254]{BHJ} $(X,Y)$ such that $X\sim F$, $Y\sim G$ and $d_{TV}(F,G)=\prob(X\ne Y)$. The Lebesgue decomposition \Ref{decom1} ensures that there exists an $A\in \mathscr{B}\(\mathbb{R}\)$ such that $F_a(A)=1$ and $F_s(A)=0$. Define $\mu_G(B)=\mathbb{P}\(X\in B\cap A,X= Y\)\le \alpha_FF_a(B)$ for $B\in \mathscr{B}\(\mathbb{R}\)$, so $\mu_G$ is absolutely continuous with respect to the Lebesgue measure. On the other hand, 
$$G(B)\ge G(B\cap A)\ge\prob(Y\in B\cap A,X=Y)
=\mu_G(B),\ \mbox{for }B\in\mathscr{B}\(\mathbb{R}\),$$ 
hence $\alpha_G\ge \mu_G(\mathbb{R})= \alpha_F-\mathbb{P}\(X\neq Y\)=\alpha_F-d_{TV}\(F, G\)>0$.
 \qed
 	
\begin{wrapfigure}{r}{0.5\textwidth}
  \vspace{-20pt}
  \begin{center} 
   \includegraphics[trim = 20mm 100mm 30mm 100mm, clip,width=0.5\textwidth]{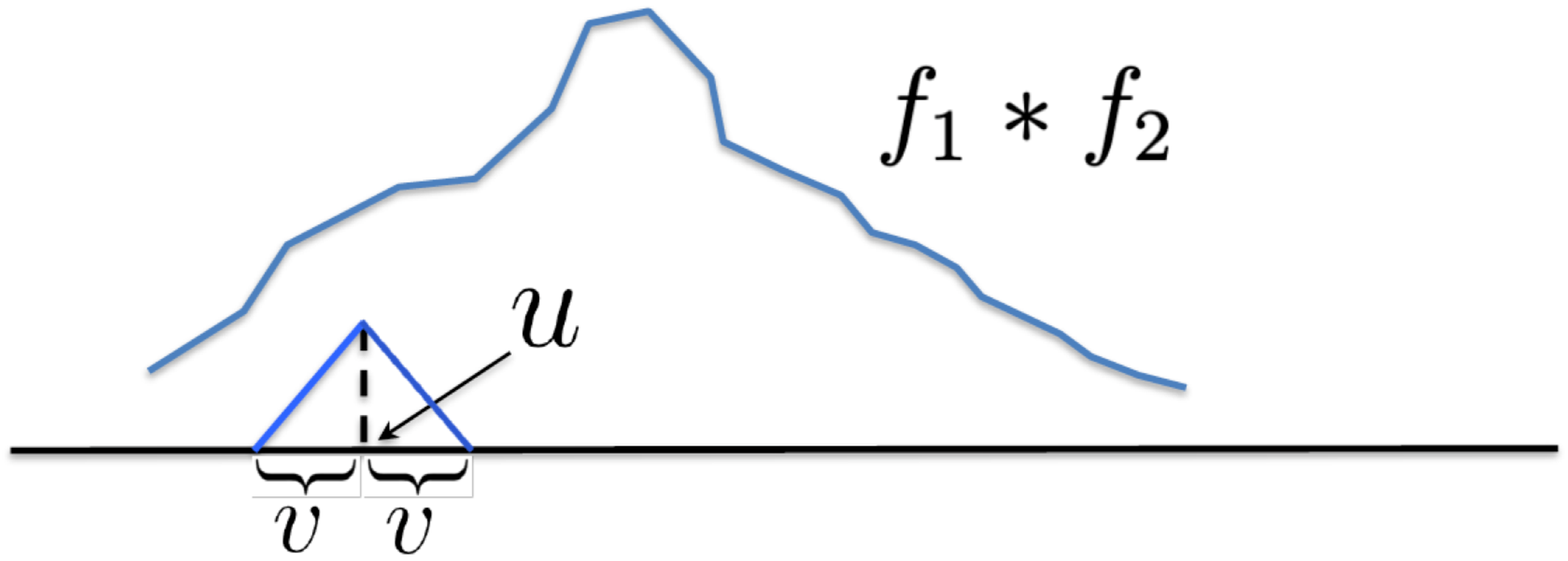}
 \vspace{-20pt}
  \caption{Existence of $u$ and $v$} %
  \label{figureone}
  \end{center}
  \vspace{-10pt}
 \end{wrapfigure}
\noindent{\it Proof of Lemma~\ref{lma2}.} Since $F_i$ is non-singular, there exists a non-zero sub-probability measure $\mu_{i}$ with a density $f_{i}$ such that $\mu_i(dx)=f_i(x)dx\le dF_i(x)$ for $x\in\mathbb{R}$. 
Without loss of generality, we can assume that both $f_{1}$ and $f_{2}$ are bounded with bounded supports, which ensures that $f_{1}\ast f_{2}$ is continuous (for the case of $f_{1}= f_{2}$, see \cite[p.~79]{Lindvall92}). In fact, as $f_{1}$ is a density, one can find a sequence of continuous functions $\{f_{1n}: \, n\ge 1\}$ satisfying $|f_{1n}-f_{1}|_1\rightarrow0$ as $n\to \infty$, where $|\cdot|_1$ is the $l_1$ norm. Now, with $|\cdot|_\infty$ denoting the supremum norm, $|f_{1n}* f_{2}-f_{1}* f_{2}|_\infty\le |f_{2}|_\infty|f_{1n}-f_{1}|_1\to 0$ as $n\to\infty$. However, the continuity is preserved under the supremum norm, the continuity of $f_{1}\ast f_{2}$ follows.

Referring to Figure~\ref{figureone}, since $f_{1}\ast f_{2}\not\equiv0$, we can find  $u\in\real$ and $v>0$ such that $f_{1}\ast f_{2}(u)>0$ and $\min_{x\in[u-v,u+v]}f_{1}\ast f_{2}(x)\ge \frac12f_{1}\ast f_{2}(u)=:b$. Let $\theta=vb$ and $a=v$, $H=\frac1{1-\theta}(F_1\ast F_2-\theta K_a\ast\delta_u)$, the claim follows. \qed
 
\noindent\noindent{\it Proof of Lemma~\ref{lma3}.} The idea of the proof is to use the radius of stabilization to limit the effect of dependence, establish that the non-singularity \Ref{non-sin} passes to the trimmed score function $\eta(\left(x,m\right), \Xi)\mathbf{1}_{R(x)\le r}$ {$\(\mbox{resp.}~\eta(\left(x,m\right), \Xi,\Gamma_\alpha)\mathbf{1}_{\bar{R}(x,\alpha)\le r}\)$} and then divide the carrier space $\Gamma_\alpha$ into maximal number of cubes so that sums of the trimmed score function on these cubes are independent. The order of the bound is then determined by the reciprocal of the number of the cubes, as in the Berry-Esseen bound. Except slightly complicated notation, the proof of the restricted case 
is the same so we first focus on the unrestricted case. From the A2.2, for the restricted case, we can find $\bar{g}$, $\bar{\eta}$ and $R$ corresponding to $\eta$ such that the stabilization radii $R$ of $\bar{\eta}$ satisfies the same stabilization property as $\eta$ in the sense of Definition~\ref{defi4}. Because $N_0$ is a bounded set, there exists an $r_1\in \mathbb{R}_+$ such that $N_0\subset B(0,r_1)$. For convenience, we write the random variables $Y:=\sum_{x\in \overline{\Xi}}\bar{g}(\Xi^x)\mathbf{1}_{d(x, N_0)<R(x)}$, $Y_r:=\sum_{x\in \overline{\Xi}}\bar{g}(\Xi^x)\mathbf{1}_{d(x, N_0)<R(x)<r}$, and write the event $\{Y\neq Y_r\}$ as $E_r$ for $r\in \mathbb{R}_+$. We can see that $$\mathbb{P}(E_r)\le \mathbb{P}(\{\mbox{there is at least one point }x\in \overline{\Xi}\mbox{ such that }d(x, N_0)\vee r<R(x)\})=:\mathbb{P}(E_r'),$$ and the right hand side is a decreasing function of $r$. We show that any one of the stabilization conditions implies that  $\mathbb{P}\(E_r\)\to0$ as $r\to\infty$, that is, $Y_r$ converges to $Y$ almost surely. In fact, 
\begin{equation}\label{lma3.2}
	\begin{aligned}
		&\mathbb{P}(E_r')
		\\	\le&\mathbb{P}\(\{\mbox{there is at least one point }x\in \overline{\Xi}\cap B(0,r_1+r)\mbox{ such that } r\le R(x)\}\)
		\\&+\mathbb{P}\(\{\mbox{there is at least one point }x\in \overline{\Xi}\cap B(0,r_1+r)^c\mbox{ such that } |x|-r_1\le R(x)\}\).
	\end{aligned}
\end{equation} Using the property of Palm process, we can see that the first term of \Ref{lma3.2} satisfies 
\begin{align}
	&\mathbb{P}\(\{\mbox{there is at least one point }x\in \overline{\Xi}\cap B(0,r_1+r)\mbox{ such that } r\le R(x)\}\)\nonumber
	\\\le & \mathbb{E}\int_{B(0,r_1+r)}\mathbf{1}_{R(x)\ge r}\overline{\Xi}(dx)\nonumber
	=\int_{B(0,r_1+r)}\mathbb{E}\mathbf{1}_{R(x, M_x,\Xi+\delta_{(x,M_x)})\ge r}\lambda dx	\\= & \int_{B(0,r_1+r)}\mathbb{P}\(R(x, M_x,\Xi+\delta_{(x,M_x)})\ge r\)\lambda dx \nonumber
	\\ \le & \int_{B(0,r_1+r)}\tau(r)\lambda dx= \frac{\lambda(r_1+r)^d\pi^{d/2}\tau(r)}{\Gamma(\frac{d}{2}+1)}\label{lma3.3}
\end{align} and the second term is bounded by 
\begin{align}
	&\mathbb{P}\(\{\mbox{there is at least one point }x\in \overline{\Xi}\cap B(0,r_1+r)^c\mbox{ such that } |x|-r_1\le R(x)\}\)\nonumber
	\\\le &\mathbb{E}\int_{B(0,r_1+r)^c}\mathbf{1}_{R(x)\ge |x|-r_1}\overline{\Xi}(dx)\nonumber
	= \int_{B(0,r_1+r)^c}\mathbb{E}\mathbf{1}_{R(x, M_x,\Xi+\delta_{(x,M_x)})\ge |x|-r_1}\lambda dx\\= & \int_{B(0,r_1+r)^c}\mathbb{P}\(R(x, M_x,\Xi+\delta_{(x,M_x)})\ge |x|-r_1\)\lambda dx \nonumber
	\\ \le & \int_{B(0,r_1+r)^c}\tau(|x|-r_1)\lambda dx= \int_{r_1+r}^\infty\frac{d\lambda t^{d-1}\pi^{d/2}\tau(t-r_1)}{\Gamma(\frac{d}{2}+1)}dt.\label{lma3.4}
\end{align}
When the score function satisfies one of the stabilization conditions, both bounds in \Ref{lma3.3} and \Ref{lma3.4} converge to $0$ as $r\rightarrow \infty$, so $\mathbb{P}(E_r)\le \mathbb{P}(E_r')\rightarrow 0$ as $r\rightarrow \infty$.

Recall that we say two measures $\mu_1\le\mu_2$ if $\mu_1(A)\le\mu_2(A)$ for all measurable sets $A$. The non-singularity~\Ref{non-sin} ensures that, with a positive probability, the conditional distribution $\law\left(Y\middle|\sigma\(\Xi_{N_0^c}\)\right)$ is non-singular. Which means that we can find a $\sigma\(\Xi_{N_0^c}\)$ measurable random measure $\xi$ on $\mathbb{R}$ which is absolutely continuous, $\xi\le \law\left(Y\middle|\sigma\(\Xi_{N_0^c}\)\right)$ $a.s.$  and $\mathbb{P}\(\xi\(\mathbb{R}\)>0\)>0$. Since $\lim_{u\downarrow 0}\mathbb{P}\(\xi(\mathbb{R})>u\)=\mathbb{P}\(\xi\(\mathbb{R}\)>0\)>0$, we can find a $p>0$ such that $\mathbb{P}\(\xi(\mathbb{R})>p\)>4p$. Because $E_r'$ is decreasing in the sense of inclusion in $r$, and $\mathbb{P}\(E_r'\)\rightarrow 0$ as $r\rightarrow \infty,$ we can find an $R_0\in \mathbb{R}_+$ such that $\mathbb{P}(E_{R_0}')\le p^2$, which ensures
\begin{equation}\label{difference}
	\mathbb{P}\(\mathbb{P}\left(E_{R_0}'\middle|\sigma\(\Xi_{N_0^c}\)\right)>\frac{p}{2}\)\le 2p.
\end{equation} Writing $\tilde{Y}:=Y\mathbf{1}_{E_{R_0}'^c}$, $A_1:=\left\{d_{TV}\(\left.Y,\tilde{Y}\right|\sigma\(\Xi_{N_0^c}\)\)>p/2\right\}$, $A_2:=\{\xi(\mathbb{R})>p\}$, then $A_1$ and $A_2$ are both $\sigma\(\Xi_{N_0^c}\)$ measurable and
$$
\mathbb{P}\(\mathbb{P}\left(d_{TV}\(\left.Y,\tilde{Y}\right|\sigma\(\Xi_{N_0^c}\)\)\right)>\frac{p}{2}\)\le 
\mathbb{P}\(\mathbb{P}\left(E_{R_0}'\middle|\sigma\(\Xi_{N_0^c}\)\right)>\frac{p}{2}\)\le 2p,
$$
giving
$\prob(A_2\cap A_1^c)>2p$. For $\omega \in A_1^c\cap A_2$, $d_{TV}\(\left.Y,\tilde{Y}\right|\sigma\(\Xi_{N_0^c}\)\)(\omega)<p/2$ and $\xi(\omega)(\mathbb{R})>p$. By Lemma~\ref{non-singular1} and \Ref{difference}, there exists an absolutely continuous $\sigma\(\Xi_{N_0^c}\)$ measurable random measure $\tilde{\xi}$ such that $\tilde{\xi}\le \law\left(\tilde{Y}\middle|\sigma\(\Xi_{N_0^c}\)\right)$ $a.s.$  and $\mathbb{P}\(\tilde{\xi}\(\mathbb{R}\)>\frac{p}{2}\)>2p$.
We write $\Xi'$ as an independent copy of $\Xi$ and the corresponding $Y_r$ and $\tilde{\xi}$ as $Y_r'$ and $\tilde{\xi}'$ respectively. 
Using Lemma~\ref{lma2}, we can find $\sigma\(\Xi_{N_0^c}, \Xi'_{N_0^c}\)$ measurable random variables
$\Theta_1\ge 0,\ \Theta_2\ge 0$ and $U\in \mathbb{R}$ such that
$\mathbb{P}(\Theta_1>0,\Theta_2>0)=4p^2$,
\begin{align*}
	\tilde{\xi}\star\tilde{\xi}'\ge \Theta_1 K_{\Theta_2}\star \delta_U.
\end{align*}
However,
$$\lim_{\epsilon\downarrow 0}\mathbb{P}\left\{\Theta_1/\Theta_2\ge
\epsilon, \Theta_2\ge\epsilon\right\}=\mathbb{P}\left\{\Theta_1>0,\Theta_2>0\right\},
$$
from Remark~\ref{remark1},  we can find an $\epsilon>0$ such that
\begin{align*}& \mathbb{P}\left\{\tilde{\xi}\star\tilde{\xi}'\ge \epsilon^2K_\epsilon\star \delta_U\right\} \ge 2p^2   
\end{align*}
for a $\sigma\(\Xi_{N_0^c}, \Xi'_{N_0^c}\)$ measurable $U$. 
From the fact that we can write
$Y_r=\tilde Y+Y_r\mathbf{1}_{E_{R_0}'}$, we have for any $B\in \mathscr{B}(\mathbb{R}\backslash\{0\})$, 
\begin{eqnarray*}
	\prob\(\left.Y_r\in B\right|\sigma\(\Xi_{N_0^c}\)\)&=&\prob\(\left.\tilde Y\in B, E_{R_0}'^c\right|\sigma\(\Xi_{N_0^c}\)\)+\prob\(\left.Y_r\in B, E_{R_0}'\right|\sigma\(\Xi_{N_0^c}\)\)\\
	&\ge& \prob\(\left.\tilde Y\in B,E_{R_0}'^c\right|\sigma\(\Xi_{N_0^c}\)\)\\
	&=&\prob\(\left.\tilde Y\in B\right|\sigma\(\Xi_{N_0^c}\)\).
\end{eqnarray*}
Hence  
\begin{align}
	&\law\left(Y_r\middle|\sigma\(\Xi_{N_0^c}\)\right)(\cdot)\ge \law\left(\tilde Y\middle|\sigma\(\Xi_{N_0^c}\)\right)(\cdot\backslash\{0\})\nonumber\\
	&\ge \tilde\xi (\cdot\backslash\{0\})=\tilde\xi (\cdot) \mbox{ a.s. for all }r\ge R_0.\label{forlma5.9}
\end{align}
Therefore, using $U\in \mathscr{A}$ to stand for $U$ being $\mathscr{A}$ measurable, we have
\begin{align*}
	&\sup_{U\in \sigma\(\Xi_{B(N_0,2r)\backslash N_0},\Xi'_{B(N_0,2r)\backslash N_0}\)} \prob\left\{\law\left(Y_r\middle|\sigma\(\Xi_{B(N_0,2r)\backslash N_0}\)\right)\star\law\left(Y_r'\middle|\sigma\(\Xi'_{B(N_0,2r)\backslash N_0}\)\right)\ge \epsilon^2K_\epsilon\star \delta_U\right\}
	\\
	&=\sup_{U\in \sigma\(\Xi_{B(N_0^c)},\Xi'_{B(N_0^c)}\)} \prob\left\{\law\left(Y_r\middle|\sigma\(\Xi_{B(N_0,2r)\backslash N_0}\)\right)\star\law\left(Y_r'\middle|\sigma\(\Xi'_{B(N_0,2r)\backslash N_0}\)\right)\ge \epsilon^2K_\epsilon\star \delta_U\right\}
	\\
	&= \sup_{U\in \sigma\(\Xi_{B(N_0^c)},\Xi'_{B(N_0^c)}\)}\prob\left\{\law\left(Y_r\middle|\sigma\(\Xi_{N_0^c}\)\right)\star\law\left(Y_r'\middle|\sigma\(\Xi'_{N_0^c}\)\right)\ge \epsilon^2K_\epsilon\star \delta_U\right\}
	\\ 
	&\ge  \sup_{U\in \sigma\(\Xi_{B(N_0^c)},\Xi'_{B(N_0^c)}\)}\prob\left\{\tilde{\xi}\star\tilde{\xi}'\ge \epsilon^2K_\epsilon\star \delta_U \right\}\\
	&\ge 2p^2,
\end{align*}
which ensures that, for any $r>R_0$, we can find a $\sigma\(\Xi_{B(N_0,2r)\backslash N_0},\Xi'_{B(N_0,2r)\backslash N_0}\)$ measurable $U$ such that \begin{equation}
	\mathbb{P}\left\{\law\left(Y_r\middle|\sigma\(\Xi_{B(N_0,2r)\backslash N_0}\)\right)\star\law\left(Y_r'\middle|\sigma\(\Xi'_{B(N_0,2r)\backslash N_0}\)\right)\ge \epsilon^2K_\epsilon\star \delta_U\right\}\ge p^2.\label{lemma4.4proofa1}
\end{equation}

If $\alpha\le (2(4r+2r_1))^d$, \Ref{statement1} is trivial with $C=\left\{2(4+\frac{2r_1}{R_0})\right\}^{d/2}$, so we now assume $\alpha> \{2(4r+2r_1)\}^d$. From the structure of $\Xi$, we can see that $\Xi(A,D)\overset{d}{=}\Xi(x+A,D)$ and $\Xi(A,D)$ is independent of $\Xi(B,D)$ for all disjoint $A$, $B\in\mathscr{B}(\real^d)$, $D\in \mathscr{T}$ and $x\in \real^d$. For a fixed $r>R_0$, we can divide $\Gamma_\alpha$ into disjoint cubes $\C_1,\cdots,\C_{m_{\alpha,r}}$ with edge length $4r+2r_1$ and centers $c_1$, $\cdots$, $c_{m_{\alpha,r}}$, aiming to maximize the number of cubes, so $m_{\alpha,r}\sim \alpha(4r+2r_1)^{-d}$, which has order $O\(\alpha r^{-d}\)$. Without loss of generality, we can assume that $m_{\alpha,r}$ is even or we simply delete one from them and the above properties still holds. For $i\le m_{\alpha,r}$, we define $A_i=c_i+N_0$, $B_i=B(A_i,r)$, $C_i=B(B_i,r)$, $D_i=C_i\backslash A_i$, $\cn_{0,\alpha,r}:=\cup_{1\le i\le m_{\alpha,r}}A_i$, $\cn_{1,\alpha,r}:=\cup_{1\le i\le m_{\alpha,r}}B_i$, $\cn_{2,\alpha,r}:=\cup_{1\le i\le m_{\alpha,r}} D_i$, $\scrF_{1,\alpha,r}:=\sigma(\Xi_{\mathbb{R}^d\backslash \cn_{0,\alpha,r}})$, $\scrF_{2,\alpha,r}:=\sigma(\Xi_{ \cn_{2,\alpha,r}})$, $W_{\alpha,r}^0=\int_{\cn_{1,\alpha,r}}\bar{g}(\Xi^{x})\mathbf{1}_{R(x)<r}\overline{\Xi}(dx)$ and $W_{\alpha,r}^1=W_{\alpha,r}-W_{\alpha,r}^0$. Note that for all $x$ such that $d(x,\partial \Gamma_\alpha)\ge r$, $\eta((x,m),\Xi,\Gamma_\alpha)\mathbf{1}_{\bar{R}(x,\alpha)<r}= \bar\eta((x,m),\Xi)\mathbf{1}_{R(x)<r}$ for all $(x,m)\in \Xi$ a.s.
From the definition of total variation distance, $d_{TV}(W_{\alpha,r},W_{\alpha,r}+\gamma)=\sup_{A\in \mathscr{B}(R)}\(\mathbb{P}(W_{\alpha,r}\in A)-\mathbb{P}(W_{\alpha,r}\in A-\gamma)\)$, hence the tower property ensures
\begin{align}
	&d_{TV}\left(W_{\alpha,r},W_{\alpha,r}+\gamma\right)\nonumber
	\\=&\sup_{A\in \mathscr{B}(\mathbb{R})}\E\left(\mathbf{1}_{W_{\alpha,r}\in A}-\mathbf{1}_{W_{\alpha,r}\in A-\gamma}\right)\nonumber
	\\=&\sup_{A\in \mathscr{B}(\mathbb{R})}\E\left(\E\left(\mathbf{1}_{W_{\alpha,r}\in A}-\mathbf{1}_{W_{\alpha,r}\in A-\gamma}|\mathscr{F}_{1,\alpha,r}\right)\right)\nonumber
	\\=&\sup_{A\in \mathscr{B}(\mathbb{R})}\E\left(\E\left(\mathbf{1}_{W_{\alpha,r}^0\in A-W_{\alpha,r}^1}-\mathbf{1}_{W_{\alpha,r}^0\in A-\gamma-W_{\alpha,r}^1}|\mathscr{F}_{1,\alpha,r}\right)\right)\nonumber
	\\\le&\E\left(\sup_{A\in \mathscr{B}(\mathbb{R})}\left[\E\left(\mathbf{1}_{W_{\alpha,r}^0\in A}-\mathbf{1}_{W_{\alpha,r}^0\in A-\gamma}|\mathscr{F}_{1,\alpha,r}\right)\right]\right)\nonumber
	\\=&\E\left(\sup_{A\in \mathscr{B}(\mathbb{R})}\left[\E\left(\mathbf{1}_{W_{\alpha,r}^0\in A}-\mathbf{1}_{W_{\alpha,r}^0\in A-\gamma}|\mathscr{F}_{2,\alpha,r}\right)\right]\right),\label{lma3.7}
\end{align}where the last equality follows from the fact that $W_{\alpha,r}^0$ depends on $\scrF_{2,\alpha,r}$ in $\scrF_{1,\alpha,r}$.

From \Ref{lma3.7}, to show \Ref{statement1}, it is sufficient to show that $$\E\left(\sup_{A\in \mathscr{B}(\mathbb{R})}\left[\E\left(\mathbf{1}_{W_{\alpha,r}^0\in A}-\mathbf{1}_{W_{\alpha,r}^0\in A-\gamma}|\scrF_{2,\alpha,r}\right)\right]\right)\le (|\gamma|\vee 1)O\left(\alpha^{-\frac{1}{2}}r^{\frac{d}{2}}\right).$$ 
Using the fact that $\int_{B_i} \bar{g}(\Xi^{x})\mathbf{1}_{R(x)<r}\overline{\Xi}(dx)$ depends only on $\sigma(\Xi_{D_i})$ in $\scrF_{2,\alpha,r}$ for $i\le m_{\alpha,r} $, and from the independence of $\sigma(\Xi_{D_i})$ for different $i$, we can see that 
\begin{align}
	\law\(W_{\alpha,r}^0|\scrF_{2,\alpha,r}\)&=\law\(\sum_{i=1}^{m_{\alpha,r}}\left.\int_{B_i} \bar{g}(\Xi^{x})\mathbf{1}_{R(x)<r}\overline{\Xi}(dx)\right|\scrF_{2,\alpha,r}\)\nonumber
	\\&=\law\(\left.\sum_{i=1}^{m_{\alpha,r}}\int_{B_i} \bar{g}(\Xi^{x})\mathbf{1}_{R(x)<r}\overline{\Xi}(dx)\right|\sigma\(\Xi_{D_i}, i\le m_{\alpha,r}\)\).\label{lemma4.4proofa2}
\end{align} 
Using \Ref{lemma4.4proofa1}, we obtain 
\begin{align*}
	&\law\(\left.\sum_{i=2j-1}^{2j}\int_{B_i} \bar{g}(\Xi^{x})\mathbf{1}_{R(x)<r}\overline{\Xi}(dx)\right|\sigma\(\Xi_{D_i}, i\le m_{\alpha,r}\)\)\\
	&=\law\(\left.X_{1,j}\(1-J_{1,j}\)+X_{2,j}J_{1,j}\(1-J_{2,j}\)+(X_{3,j}+U_j)J_{1,j}J_{2,j}\right|\sigma\(\Xi_{D_{2j-1}}, \Xi_{D_{2j}}\)\),
\end{align*} 
where $J_{1,j}$, $J_{2,j}$ and $U_j$ are $\sigma\(\Xi_{D_{2j-1}}, \Xi_{D_{2j}}\)$ measurable with $\prob(J_{1,j}=1)=1-\prob(J_{1,j}=0)=p^2$, $\prob(J_{2,j}=1)=1-\prob(J_{2,j}=0)=\e^2$, $J_{1,j} \indep J_{2,j}$, $X_{1,j} $ and $X_{2,j} $ are $\sigma\(\Xi_{B_{2j-1}}, \Xi_{B_{2j}}\)$ measurable, and $X_{3,j}\sim K_\e,\ 1\le j\le m_{\alpha,r}/2,$ are i.i.d. and independent of $\sigma\(\Xi_{D_i}, i\le m_{\alpha,r}\)$.
Hence, define $\Sigma_1:= \sum_{j=1}^{m_{\alpha,r}/2}(X_{1,j}\(1-J_{1,j}\)+X_{2,j}J_{1,j}\(1-J_{2,j}\)+(X_{3,j}+U_j)J_{1,j}J_{2,j})$, $\Sigma_2:= \sum_{j=1}^{m_{\alpha,r}/2}X_{3,j}J_{1,j}J_{2,j}$, $\Sigma_{3,l}:= \sum_{j=1}^lX_{3,j}$ and $I\sim${\rm Binomial}$(m_{\alpha,r}/2,\e^2p^2)$ which is independent of $\{X_{3,j}:\ j\le m_{\alpha,r}/2\}$, it follows from \Ref{lma3.7} and \Ref{lemma4.4proofa2} that
\begin{align}
	&d_{TV}\left(W_{\alpha,r},W_{\alpha,r}+\gamma\right)\nonumber
	\\ \le&\E\left(\sup_{A\in \mathscr{B}(\mathbb{R})}\left[\E\left(\mathbf{1}_{\Sigma_1\in A}-\mathbf{1}_{\Sigma_1\in A-\gamma}|\sigma\(\Xi_{D_i}, i\le m_{\alpha,r}\)\right)\right]\right)\nonumber
	\\ \le&\E\left(\sup_{A\in \mathscr{B}(\mathbb{R})}\left[\E\left(\mathbf{1}_{\Sigma_2\in A}-\mathbf{1}_{\Sigma_2\in A-\gamma}|\sigma\(\Xi_{D_i}, i\le m_{\alpha,r}\)\right)\right]\right)\nonumber
	\\ =&\E\left(\sup_{A\in \mathscr{B}(\mathbb{R})}\left[\E\left(\mathbf{1}_{\Sigma_{3,I}\in A}-\mathbf{1}_{\Sigma_{3,I}\in A-\gamma}|I\right)\right]\right)\nonumber
	\\\le&\prob(I\le (\E I)/2)+\sum_{(\E I)/2<j\le m_{\alpha,r}/2}\sup_{A\in \mathscr{B}(\mathbb{R})}\left[\E\left(\mathbf{1}_{\Sigma_{3,j}\in A}-\mathbf{1}_{\Sigma_{3,j}\in A-\gamma}\right)\right]\prob(I=j)\nonumber
	\\\le& O\left(\alpha^{-1}r^d\right)+O\left(\alpha^{-\frac{1}{2}}r^{\frac{d}{2}}\right)|\gamma|= \left(|\gamma|\vee 1\right)O\left(\alpha^{-\frac{1}{2}}r^{\frac{d}{2}}\right), \label{lma3.10}
\end{align} where the first term of \Ref{lma3.10} is from Chebyshev's inequality and the second terms is due to Lemma~\ref{lma1}. This completes the proof of \Ref{statement1}. 

In terms of \Ref{statement2}, since range-bound implies polynomially stabilizing with arbitrary order $\beta$, \Ref{statement1} still holds for all $r>R_0$. On the other hand, $\bar{W}_\alpha=\bar{W}_r$ a.s. when $r>t$ for some positive constant $t$, \Ref{statement2} follows by taking $r=R_0\vee t+1$. 

The claim \Ref{statement3} can be proved by replacing $W_\alpha$ with $\bar{W}_\alpha$; $W_{\alpha,r}$ with $\bar{W}_{\alpha,r}$; {$\bar{g}$ with $g$,} $W^0_{\alpha,r}$ with $\bar{W}^0_{\alpha,r}$; $W^1_{\alpha,r}$ with $\bar{W}^1_{\alpha,r}$; $\Xi^{x}$ by $\Xi^{\Gamma_\alpha,x}$; $R(x)$ with $\bar{R}(x,\alpha)$; $R(x,M_x,\Xi+\delta_{(x, M_x)})$ with $\bar{R}(x,M_x,\alpha,\Xi+\delta_{(x, M_x)})$ and redefining $\scrF_{1,\alpha,r}:=\sigma(\Xi_{\Gamma_\alpha\backslash \cn_{0,\alpha,r}})$. The bound \Ref{statement4} can be argued in the same way as that for \Ref{statement2}. \qed

\noindent{\it Proof of Corollary~\ref{cor1}.}  The proof can be easily adapted from the second half of the proof of Lemma~\ref{lma3} and we start with \Ref{lma3coro02}. If $\alpha^{-
	\frac{1}{d}}(1-2C)^{-1}(4r+2r_1)>\frac{1}{3}$, \Ref{lma3coro02} is obvious because the total variation distance is bounded above by $1$. Now we assume $\alpha^{-
	\frac{1}{d}}(1-2C)^{-1}(4r+2r_1)\le\frac{1}{3}$. Similar to the proof of Lemma~\ref{lma3}, we embed disjoint cubes with edge length $4r+2r_1$ into $\Gamma_\alpha\backslash \(N_{\alpha,r}^{(1)}\cup N_{\alpha,r}^{(2)}\cup N_{\alpha,r}^{(3)}\)$, aiming to maximize the number $m_{\alpha,r}$ of the cubes. Without loss, we assume that $m_{\alpha,r}$ is even. Then, we have $$\alpha(1-2C)^{d}(12r+6r_1)^{-d}{-1}\le m_{\alpha,r}\le \alpha(1-2C)^{d}(4r+2r_1)^{-d},$$ giving $m_{\alpha,r}=O(\alpha r^{-d})$. 

We use the same notations as in the proof of Lemma~\ref{lma3} but with $\Gamma_\alpha$ replaced by $\Gamma_\alpha\backslash(N_{\alpha,r}^{(1)}\cup N_{\alpha,r}^{(2)}\cup N_{\alpha,r}^{(3)})$ and define $\mathscr{F}_{2,\alpha,r}':=\sigma\(\Xi_{\cn_{2,\alpha,r}\cup N_{\alpha,r}^{(2)}}\)$. Bearing in mind that $N_{\alpha,r}^{(1)}\cup N_{\alpha,r}^{(2)}\cup N_{\alpha,r}^{(3)}$ is excluded in the $m_{\alpha,r}$ cubes, we have $\mathscr{F}_{0,\alpha,r}\subset \mathscr{F}_{1,\alpha,r}$, giving the following analogous result of \Ref{lma3.7}:
\begin{align*}
	& d_{TV}\left(\bar{W}_{\alpha,r}',\bar{W}_{\alpha,r}'+h_{\alpha,r}\(\Xi_{N_{\alpha, r}^{(2)}}\)\middle|\mathscr{F}_{0,\alpha,r}\right)
	\\=&\sup_{A\in \mathscr{B}(\mathbb{R})}\E\left(\mathbf{1}_{\bar{W}_{\alpha,r}'\in A}-\mathbf{1}_{\bar{W}_{\alpha,r}\in A-h_{\alpha,r}\(\Xi_{N_{\alpha, r}^{(2)}}\)}\middle|\mathscr{F}_{0,\alpha,r}\right)
	\\=&\sup_{A\in \mathscr{B}(\mathbb{R})}\E\left(\E\left(\mathbf{1}_{\bar{W}_{\alpha,r}\in A}-\mathbf{1}_{\bar{W}_{\alpha,r}\in A-h_{\alpha,r}\(\Xi_{N_{\alpha, r}^{(2)}}\)}\middle|\mathscr{F}_{1,\alpha,r}\right)\middle|\mathscr{F}_{0,\alpha,r}\right)
	\\=&\sup_{A\in \mathscr{B}(\mathbb{R})}\E\left(\E\left(\mathbf{1}_{\bar{W}_{\alpha,r}^0\in A}-\mathbf{1}_{\bar{W}_{\alpha,r}^0\in A-h_{\alpha,r}\(\Xi_{N_{\alpha, r}^{(2)}}\)}\middle|\mathscr{F}_{1,\alpha,r}\right)\middle|\mathscr{F}_{0,\alpha,r}\right)
	\\\le &\E\left(\sup_{A\in \mathscr{B}(\mathbb{R})}\[\E\left(\mathbf{1}_{\bar{W}_{\alpha,r}^0\in A}-\mathbf{1}_{\bar{W}_{\alpha,r}^0\in A-h_{\alpha,r}\(\Xi_{N_{\alpha, r}^{(2)}}\)}\middle|\mathscr{F}_{1,\alpha,r}\right)\]\middle|\mathscr{F}_{0,\alpha,r}\right)
	\\ =&\E\left(\sup_{A\in \mathscr{B}(\mathbb{R})}\[\E\left(\mathbf{1}_{\bar{W}_{\alpha,r}^0\in A}-\mathbf{1}_{\bar{W}_{\alpha,r}^0\in A-h_{\alpha,r}\(\Xi_{N_{\alpha, r}^{(2)}}\)}\middle|\mathscr{F}_{2,\alpha,r}'\right)\]\middle|\mathscr{F}_{0,\alpha,r}\right).
\end{align*} 
The remaining part is a line-by-line repetition of the proof of Lemma~\ref{lma3} with $\gamma$ replaced by $h_{\alpha,r}\(\Xi_{N_{\alpha, r}^{(2)}}\)$ and expectation replaced by the conditional expectation given $\mathscr{F}_{0,\alpha,r}$, leading to 
\begin{align}
	& d_{TV}\left(\bar{W}_{\alpha,r}',\bar{W}_{\alpha,r}'+h_{\alpha,r}\(\Xi_{N_{\alpha, r}^{(2)}}\)\middle|\mathscr{F}_{0,\alpha,r}\right)\nonumber
	\\\le&\E\left(\sup_{A\in \mathscr{B}(\mathbb{R})}\[\E\left(\mathbf{1}_{\bar{W}_{\alpha,r}^0\in A}-\mathbf{1}_{\bar{W}_{\alpha,r}^0\in A-h_{\alpha,r}\(\Xi_{N_{\alpha, r}^{(2)}}\)}\middle|\mathscr{F}_{2,\alpha,r}'\right)\]\middle|\mathscr{F}_{0,\alpha,r}\right)\nonumber
	\\\le&\E\left(\E\left[\sup_{A\in \mathscr{B}(\mathbb{R})}\left(\mathbf{1}_{{\Sigma_{3,I}}\in A}-\mathbf{1}_{{\Sigma_{3,I}}\in A-h_{\alpha,r}\(\Xi_{N_{\alpha, r}^{(2)}}\)}\middle|I,\sigma\(\Xi_{N_{\alpha,r}^{(2)}}\)\right)\right] \middle|\mathscr{F}_{0,\alpha,r}\right)\label{cor1.1}
	\\\le&\sum_{(\E I)/2<j\le m_{\alpha,r}/2}\prob(I=j)\mathbb{E}\left(\E\left[\sup_{A\in \mathscr{B}(\mathbb{R})}\left(\mathbf{1}_{{\Sigma_{3,j}}\in A}-\mathbf{1}_{{\Sigma_{3,j}}\in A-h_{\alpha,r}\(\Xi_{N_{\alpha, r}^{(2)}}\)}\middle| \sigma\(\Xi_{N_{\alpha,r}^{(2)}}\)\right)\right]\middle|\mathscr{F}_{0,\alpha,r}\right)\nonumber
	\\&+\prob(I\le (\E I)/2)\nonumber
	\\\le& O\left(\alpha^{-1}r^d\right)+O\left(\alpha^{-\frac{1}{2}}r^{\frac{d}{2}}\right)\mathbb{E}\left(\left|h_{\alpha,r}\(\Xi_{N_{\alpha, r}^{(2)}}\)\right|\middle| \mathscr{F}_{0,\alpha,r}\right)\label{cor1.5}
	\\=& \mathbb{E}\left(\left|h_{\alpha,r}\(\Xi_{N_{\alpha, r}^{(2)}}\)\right|\vee 1\middle| \mathscr{F}_{0,\alpha,r}\right)O\left(\alpha^{-\frac{1}{2}}r^{\frac{d}{2}}\right),\nonumber
\end{align}
where \Ref{cor1.1} follows from the fact that $\sup_{A\in \mathscr{B}(\mathbb{R})}\left[\mathbf{1}_{{\Sigma_{3,I}}\in A}-\mathbf{1}_{{\Sigma_{3,I}}\in A-h_{\alpha,r}\(\Xi_{N_{\alpha, r}^{(2)}}\)}\right]$ is a function of $I$, $\Xi_{N_{\alpha, r}^{(2)}}$, the first term of \Ref{cor1.5} is from Chebyshev's inequality and the second term is due to Lemma~\ref{lma1}. This completes the proof for the statement of $\bar{W}_{\alpha,r}$.

The claim \Ref{lma3coro01} can be proved by replacing corresponding counterparts  $\bar{W}_{\alpha,r}$ with $W_{\alpha,r}$; $\bar{W}'_{\alpha,r}$ with $W'_{\alpha,r}$; $\bar{W}^0_{\alpha,r}$ with $W^0_{\alpha,r}$. \qed

The moments of $W_{\alpha,r}$ and $W_\alpha$ (resp. $\bar{W}_{\alpha,r}$ and $\bar{W}_\alpha$) can be established using the ideas in \cite[Section~4]{XY15}. Let $\|X\|_p:=\mathbb{E}\(|X|^p\)^{\frac{1}{p}}$ be the $L_p$ norm of $X$ provided it is finite.

\begin{lma}\label{lma11} 
	\begin{description}
		\item{(a)} (unrestricted case) If the score function $\eta$ satisfies $k'$-th moment condition \Ref{thm2.1} with $k'>k\ge 1$, then
		$\max_{0< l\le k}\left\{\|W_{\alpha}\|_l,\|W_{\alpha,r}\|_l\right\}\le C\alpha.$
		\item{(b)} (restricted case) If the score function $\eta$ satisfies $k'$-th moment condition \Ref{thm2.1r} with $k'>k\ge 1$, then $\max_{0< l\le k}\left\{\|\bar{W}_{\alpha}\|_l,\|\bar{W}_{\alpha,r}\|_l\right\}\le C\alpha.$
	\end{description}
\end{lma}
\noindent{\it Proof.} The proof is adapted from that of \cite[Lemma~4.1]{XY15}. We use the notations as in the proof of Lemma~\ref{lma3} and start with the restricted case. To this end, it suffices to show $\|\bar{W}_{\alpha}\|_k{\vee \|W_{\alpha,r}\|_k}\le C\alpha$ and the claim follows from H\"older's inequality. Let $N_\alpha:=\left|\bar\Xi_{\Gamma_{\alpha}}\right|$, then $N_\alpha$ follows Poisson distribution with parameter $\alpha \lambda$. Using Minkowski's inequality, we obtain

\begin{align}
	\|\bar{W}_{\alpha}\|_k&\le\left\|\sum_{x\in \bar{\Xi}_{\Gamma_\alpha}}\left|g_\alpha\(x,\Xi\)\right|\right\|_k\nonumber 
	\\=&\left\|\sum_{x\in \bar{\Xi}_{\Gamma_\alpha}}\left|g_\alpha\(x,\Xi\)\right|\(\mathbf{1}_{N_\alpha\le \alpha\lambda}+\sum_{j=0}^\infty\mathbf{1}_{\alpha\lambda 2^{j}<N_\alpha\le \alpha\lambda 2^{j+1}}\)\right\|_k\nonumber 
	\\\le& \left\|\sum_{x\in \bar{\Xi}_{\Gamma_\alpha}}\left|g_\alpha\(x,\Xi\)\right|\mathbf{1}_{N_\alpha\le \alpha\lambda}\right\|_k+\sum_{j=0}^\infty \left\|\sum_{x\in \bar{\Xi}_{\Gamma_\alpha}}\left|g_\alpha\(x,\Xi\)\right|\mathbf{1}_{\alpha\lambda 2^{j}<N_\alpha\le \alpha\lambda 2^{j+1}}\right\|_k.\label{lma11.1}
\end{align} 
Let $s=\frac{k'}{k}>1$ and $t$ be its conjugate, i.e., $\frac{1}{s}+\frac{1}{t}=1$, using H\"older's inequality and Minkowski's inequality, for any $j\in \mathbb{N}$, we have
\begin{align}
	&\left\|\sum_{x\in \bar{\Xi}_{\Gamma_\alpha}}\left|g_\alpha\(x,\Xi\)\right|\mathbf{1}_{\alpha\lambda 2^{j}<N_\alpha\le \alpha\lambda 2^{j+1}}\right\|_k\nonumber
	\\=&\left\|\sum_{x\in \bar{\Xi}_{\Gamma_\alpha}}\left|g_\alpha\(x,\Xi\)\right|\mathbf{1}_{N_\alpha\le \alpha\lambda 2^{j+1}}\mathbf{1}_{\alpha\lambda 2^{j}<N_\alpha }\right\|_k\nonumber
	\\\le&\left\|\sum_{x\in \bar{\Xi}_{\Gamma_\alpha}}\left|g_\alpha\(x,\Xi\)\right|\mathbf{1}_{N_\alpha\le \alpha\lambda 2^{j+1}}\right\|_{k'}\(\mathbb{P}\(N_\alpha> \alpha\lambda 2^{j}\)\)^{\frac{1}{kt}}\nonumber
	\\ =&   \left\|\sum_{x\in \bar{\Xi}_{\Gamma_\alpha}}\left|g_\alpha\(x,\Xi\)\right|\mathbf{1}_{N_\alpha\le \alpha\lambda 2^{j+1}}\right\|_{k'}\mathbb{P}\(N_\alpha-\alpha\lambda>\alpha\lambda \(2^{j}-1\)\)^{\frac{1}{kt}}.\label{lma11.2}
\end{align} 
For the term $\|\cdot\|_{k'}$ in \Ref{lma11.2}, we have
\begin{align}
	&\left\|\sum_{x\in \bar{\Xi}_{\Gamma_\alpha}}\left|g_\alpha\(x,\Xi\)\right|\mathbf{1}_{N_\alpha\le n}\right\|_{k'}\nonumber
	\\=&\left\{\mathbb{E}\[\(\sum_{j=1}^n\sum_{x\in \bar{\Xi}_{\Gamma_\alpha}}\left|g_\alpha\(x,\Xi\)\right|\mathbf{1}_{N_\alpha=j}\)^{k'}\]\right\}^{\frac{1}{k'}}\nonumber
	\\=& \left\{\sum_{j=1}^n\mathbb{E}\[\(\sum_{x\in \bar{\Xi}_{\Gamma_\alpha}}\left|g_\alpha\(x,\Xi\)\right|\mathbf{1}_{N_\alpha=j}\)^{k'}\]\right\}^{\frac{1}{k'}},\label{lma11.4}
\end{align} where the first equality holds because $\sum_{x\in \bar{\Xi}_{\Gamma_\alpha}}\left|g_\alpha\(x,\Xi\)\right|^{k'}\mathbf{1}_{N_\alpha=0}=0$ and the last equality follows from the fact that $\{N_\alpha=j\}$, $1\le j\le n$, are disjoint events. On $\{N_\alpha=j\}$ for some fixed $j\in \mathbb{N}$, if we write $j$ points in $\Xi\cap \Gamma_\alpha$ as $\{(x_1,m_1),\dots,(x_j,m_j)\}$ and let $\{\(U_{\alpha,i},M_i\)\}_{i\in \mathbb{N}}$ be a sequence of \iid\ random elements having distribution $U\(\Gamma_\alpha\)\times \mathscr{L}_T$ and be independent of $\Xi$, where $U\(\Gamma_\alpha\)$ is the uniform distribution on $\Gamma_\alpha$, then 
\begin{align}
	&\mathbb{E}\[\(\sum_{x\in \bar{\Xi}_{\Gamma_\alpha}}\left|g_\alpha\(x,\Xi\)\right|\mathbf{1}_{N_\alpha=j}\)^{k'}\]\nonumber
	\\ \le& \left\{\sum_{i=1}^j\mathbb{E}\[\left|g_\alpha\((x_i,m_i),\Xi\)\right|^{k'}\mathbf{1}_{N_\alpha=j}\]^{\frac{1}{k'}}\right\}^{k'}\nonumber
	\\= &j^{k'}\mathbb{E}\[\left|g_\alpha\(\(U_{\alpha,i},M_i\), \(\sum_{i=1}^j \delta_{\(U_{\alpha,i},M_i\)}\)\)\right|^{k'}\mathbf{1}_{N_\alpha=j}\],\label{lma11.5}
\end{align} where the inequality follows from Minkovski's inequality and the equality follows from the fact that when $\left|\bar{\Xi}\cap \Gamma_\alpha\right|$ is fixed, points in $\bar{\Xi}\cap \Gamma_\alpha$ are independent and follow uniform distribution on $\Gamma_\alpha$. Combining \Ref{lma11.4} and \Ref{lma11.5}, we have 
\begin{align}
	&\left\|\sum_{x\in \bar{\Xi}_{\Gamma_\alpha}}\left|g_\alpha\(x,\Xi\)\right|\mathbf{1}_{N_\alpha\le n}\right\|_{k'}\nonumber\nonumber
	\\ \le &\left\{\sum_{j=1}^nj^{k'}\mathbb{E}\[\left|g_\alpha\(\(U_{\alpha,1},M_1\), \(\sum_{i=1}^j \delta_{\(U_{\alpha,i},M_i\)}\)\)\right|^{k'}\mathbf{1}_{N_\alpha=j}\]\right\}^{\frac{1}{k'}}\nonumber
	\\ =&\left\{\sum_{j=1}^n\lambda \alpha j^{k'-1}\mathbb{E}\[\left|g_\alpha\(\(U_{\alpha,1},M_1\), \(\sum_{i=1}^j \delta_{\(U_{\alpha,i},M_i\)}\)\)\right|^{k'}\mathbf{1}_{N_\alpha=j-1}\]\right\}^{\frac{1}{k'}}\nonumber
	\\ \le & (\lambda\alpha)^{\frac{1}{k'}}n^{\frac{k'-1}{k'}}\left\{\mathbb{E}\[\sum_{j=0}^{n-1}\left|g_\alpha\(\(U_{\alpha,1},M_1\), \(\sum_{i=1}^{j+1} \delta_{\(U_{\alpha,i},M_i\)}\)\)\right|^{k'}\mathbf{1}_{N_\alpha=j}\]\right\}^{\frac{1}{k'}}\nonumber
	\\ \le & (\lambda\alpha)^{\frac{1}{k'}}n^{\frac{k'-1}{k'}}\left\{\mathbb{E}\[\sum_{j=0}^{\infty}\left|g_\alpha\(\(U_{\alpha,1},M_1\), \(\sum_{i=1}^{j+1} \delta_{\(U_{\alpha,i},M_i\)}\)\)\right|^{k'}\mathbf{1}_{N_\alpha=j}\]\right\}^{\frac{1}{k'}}\nonumber
	\\ =& (\lambda\alpha)^{\frac{1}{k'}}n^{\frac{k'-1}{k'}}\left\{\int_{\Gamma_\alpha}\mathbb{E}\[\left|g_\alpha\((x,M),{\Xi_{\Gamma_\alpha}+\delta_{(x,M)}}\)\right|^{k'}\frac{1}{\alpha}dx\]\right\}^{\frac{1}{k'}}\nonumber
	\\\le& (\lambda\alpha)^{\frac{1}{k'}}n^{\frac{k'-1}{k'}} C_0^{\frac{1}{k'}},\label{lma11.6}
\end{align}
where the first equality follows from the fact that $N_\alpha$ is independent of $\{\(U_{\alpha,i},M_i\)\}_{i\in \mathbb{N}}$ and $\mathbb{P}(N_\alpha=j)=\frac{\lambda \alpha}{j}\mathbb{P}(N_\alpha=j-1)$, the last equality follows from the construction of marked Poisson point process.  

Combining \Ref{lma11.2} and \Ref{lma11.6}, we have
\begin{align}
	\left\|\sum_{x\in \bar{\Xi}_{\Gamma_\alpha}}\left|g_\alpha\(x,\Xi\)\right|\mathbf{1}_{\alpha\lambda 2^{j}<N_\alpha\le \alpha\lambda 2^{j+1}}\right\|_k\le  \alpha\lambda 2^{\frac{(k'-1)(j+1)}{k'}}C_0^{\frac{1}{k'}} \mathbb{P}\(N_\alpha-\alpha\lambda>\alpha\lambda \(2^{j}-1\)\)^{\frac{1}{kt}}.\label{lma11.7}
\end{align}
Using \Ref{lma11.6} and H\"older's inequality, we have 
\begin{align}
	&\left\|\sum_{x\in \bar{\Xi}_{\Gamma_\alpha}}\left|g_\alpha\(x,\Xi\)\right|\mathbf{1}_{N_\alpha\le \alpha\lambda}\right\|_k\le\left\|\sum_{x\in \bar{\Xi}_{\Gamma_\alpha}}\left|g_\alpha\(x,\Xi\)\right|\mathbf{1}_{N_\alpha\le \alpha\lambda}\right\|_{k'}\le \alpha\lambda C_0^{\frac{1}{k'}}.\label{lma11.3}
\end{align} Combining\Ref{lma11.7} and \Ref{lma11.3}, together with the fact that $\mathbb{P}\(N_\alpha-\alpha\lambda>\alpha\lambda k\)$ decrease exponentially fast with respective to $k$, we have from \Ref{lma11.1} that
$$\|\bar{W}_{\alpha}\|_k\le\left\|\sum_{x\in \bar{\Xi}_{\Gamma_\alpha}}\left|g_\alpha\(x,\Xi\)\right|\right\|_k\le C\alpha.$$ 
The proof of (b) is completed by observing that, for arbitrary $r\in \mathbb{R}_+$, 
$$\|\bar{W}_{\alpha,r}\|_k=\left\|\sum_{x\in \bar{\Xi}_{\Gamma_\alpha}}g_\alpha\(x,\Xi\)\mathbf{1}_{R(x)\le r}\right\|_k\le\left\|\sum_{x\in \bar{\Xi}_{\Gamma_\alpha}}\left|g_\alpha\(x,\Xi\)\right|\right\|_k\le C\alpha.$$ 

The claim (a) can be established by replacing $\bar{W}_{\alpha}$ with $W_\alpha$, $\bar{W}_{\alpha,r}$ with $W_{\alpha,r}$;  $\sum_{i=1}^j \delta_{\(U_{\alpha,i},M_i\)}$ with $\sum_{i=1}^j \delta_{\(U_{\alpha,i},M_i\)}+\Xi_{\Gamma_\alpha^c}$; {$g_\alpha(x,\mathscr{X})$ as $g(\mathscr{X}^{x})$}. \qed

\vskip10pt
\begin{re} {\rm~The proof of Lemma~\ref{lma11} does not depend on the shape of $\Gamma_\alpha$, so the claims still hold if we replace $\Gamma_\alpha$ with a set $A\in\mathscr{B}(\mathbb{R}^d)$ and $\alpha$ in the upper bound with the volume of $A$.}
\end{re}

With these preparations, we are ready to bound the differences $\left|\var\(W_{\alpha}\)-\var\(W_{\alpha,r}\)\right|$ and $\left|\var\(\bar{W}_{\alpha}\)-\var\(\bar{W}_{\alpha,r}\)\right|$.

\begin{lma}\label{lma12}
	\begin{description}
		\item{(a)} (unrestricted case) Assume the score function $\eta$ satisfies $k'$th moment condition~\Ref{thm2.1} for some $k'>2$. If $\eta$ is exponentially stabilizing in Definition~\ref{defi4}, then there exist positive constants $\alpha_0$ and $C$ such that 
		\begin{equation}\left|\var\(W_{\alpha}\)-\var\(W_{\alpha,r}\)\right|\le \frac{1}{\alpha}\label{lma12s1}\end{equation} for all $\alpha\ge\alpha_0$ and $r\ge C\ln (\alpha)$. If $\eta$ is polynomially stabilizing in Definition~\ref{defi4} with parameter $\beta$, then for any $k\in (2,k') $, then there exists a positive constant $C$ such that
		\begin{align}
			\left|\var\(W_{\alpha}\)-\var\(W_{\alpha,r}\)\right|&\le C\(\alpha^{\frac{3k-2}{k}}r^{-\beta\frac{k-2}{k}}\)\vee \(\alpha^{\frac{3k-1}{k}}r^{-\beta\frac{k-1}{k}}\)\label{lma12s2}
		\end{align} for all $r\le \alpha^{\frac{1}{d}}$.
		\item{(b)} (restricted case) Assume the score function $\eta$ satisfies $k'$th moment condition~\Ref{thm2.1r} for some $k'>2$. If $\eta$ is exponentially stabilizing in Definition~\ref{defi4r}, then there exist positive constants $\alpha_0$ and $C$ such that 
		\begin{equation}\left|\var\(\bar{W}_{\alpha}\)-\var\(\bar{W}_{\alpha,r}\)\right|\le \frac{1}{\alpha}\label{lma12s3}\end{equation} for all $\alpha\ge\alpha_0$ and $r\ge C\ln (\alpha)$. If $\eta$ is polynomially stabilizing in Definition~\ref{defi4r} with parameter $\beta$, then for any $k\in (2,k') $, then there exists a positive constant $C$ such that
		\begin{align}
			\left|\var\(\bar{W}_{\alpha}\)-\var\(\bar{W}_{\alpha,r}\)\right|&\le C\(\alpha^{\frac{3k-2}{k}}r^{-\beta\frac{k-2}{k}}\)\vee \(\alpha^{\frac{3k-1}{k}}r^{-\beta\frac{k-1}{k}}\)\label{lma12s4}
		\end{align} for all $r\le \alpha^{\frac{1}{d}}$.
	\end{description}
\end{lma}

\noindent{\it Proof.} We start with \Ref{lma12s3}. From Lemma~\ref{lma11}~(b), for fixed $k\in (2,k') $, we have 
\begin{equation}\max_{0< l\le k}\left\{\|\bar{W}_{\alpha}\|_l,\|\bar{W}_{\alpha,r}\|_l\right\}\le C_0\alpha\label{lmapr1}
\end{equation} for some positive constant $C_0$. Without loss, we assume $\alpha_0>1$.
Since \begin{equation}\left|\var\(\bar{W}_{\alpha}\)-\var\(\bar{W}_{\alpha,r}\)\right|\le \left|\mathbb{E}\(\bar{W}_\alpha^2-\bar{W}_{\alpha,r}^2\)\right|+\left|\(\mathbb{E}\bar{W}_\alpha\)^2-\(\mathbb{E}\bar{W}_{\alpha,r}\)^2\right|, \label{lma12.0}
\end{equation} 
assuming that the score function is exponentially stabilizing \Ref{defi4r}, we show that each of the terms at the right hand side of \Ref{lma12.0} is bounded by $\frac{1}{2\alpha}$ for $\alpha$ and $r$ sufficient large. Clearly, the definition of $\bar{W}_{\alpha,r}$ implies that $\bar{W}_\alpha^2-\bar{W}_{\alpha,r}^2= 0$ if $\bar{R}(x,\alpha)\le r$ for all $x\in\bar\Xi_{\Gamma_\alpha}$, hence it remains to tackle $E_{r,\alpha}:=\{\bar{R}(x,\alpha)\le r~\mbox{for all }x\in\bar\Xi_{\Gamma_\alpha}\}^c$. As shown in the proof of Lemma~\ref{lma105}, $\mathbb{P}\(E_{r,\alpha}\)\le\alpha C_1e^{-C_2r}$, which, together with H\"older's inequality,
ensures 
\begin{align}
	\left|\mathbb{E}\(\bar{W}_\alpha^2-\bar{W}_{\alpha,r}^2\)\right|&=\left|\mathbb{E}\[\(\bar{W}_\alpha^2-\bar{W}_{\alpha,r}^2\)\mathbf{1}_{E_{r,\alpha}}\]\right|\nonumber
	\\&\le \|\bar{W}_\alpha^2-\bar{W}_{\alpha,r}^2\|_{\frac{k}{2}}\|\mathbf{1}_{E_{r,\alpha}}\|_{\frac{k}{k-2}}\nonumber
	\\&\le \(\|\bar{W}_\alpha^2\|_{\frac{k}{2}}+\|\bar{W}_{\alpha,r}^2\|_{\frac{k}{2}}\)\mathbb{P}(E_{r,\alpha})^{\frac{k-2}{k}}\nonumber
	\\&= \(\|\bar{W}_\alpha\|_{k}^2+\|\bar{W}_{\alpha,r}\|_{k}^2\)\mathbb{P}(E_{r,\alpha})^{\frac{k-2}{k}}\le 2\(C_0\alpha\)^2 \(\alpha C_1e^{-C_2r}\)^{\frac{k-2}{k}}.\label{lma12.1}
\end{align} 

For the remaining term of \Ref{lma12.0}, we have $$\left|\(\mathbb{E}\bar{W}_\alpha\)^2-\(\mathbb{E}\bar{W}_{\alpha,r}\)^2\right|=\left|\mathbb{E}\bar{W}_\alpha-\mathbb{E}\bar{W}_{\alpha,r}\right|\left|\mathbb{E}\bar{W}_\alpha+\mathbb{E}\bar{W}_{\alpha,r}\right|.$$ The bound \Ref{lmapr1} implies $\left|\mathbb{E}\bar{W}_\alpha+\mathbb{E}\bar{W}_{\alpha,r}\right|\le 2C_0\alpha $. However, using H\"older's inequality, Minkowski's inequality and \Ref{lmapr1} again, we have 
\begin{align}
	\left|\mathbb{E}\bar{W}_\alpha-\mathbb{E}\bar{W}_{\alpha,r}\right|&=\left|\mathbb{E}\[\(\bar{W}_\alpha-\bar{W}_{\alpha,r}\)\mathbf{1}_{E_{r,\alpha}}\]\right|\nonumber
	\\&\le \|\bar{W}_\alpha-\bar{W}_{\alpha,r}\|_k\|\mathbf{1}_{E_{r,\alpha}}\|_{\frac{k}{k-1}}\nonumber
	\\&\le \(\|\bar{W}_\alpha\|_{k}+\|\bar{W}_{\alpha,r}\|_{k}\)\mathbb{P}(E_{r,\alpha})^{\frac{k-1}{k}}\nonumber
	\\&\le 2C_0\alpha \(\alpha C_1e^{-C_2r}\)^{\frac{k-1}{k}},\label{lma12.2}
\end{align} 
giving \begin{equation}
	\left|\(\mathbb{E}\bar{W}_\alpha\)^2-\(\mathbb{E}\bar{W}_{\alpha,r}\)^2\right|\le 4\(C_0\alpha\)^2\(\alpha C_1e^{-C_2r}\)^{\frac{k-1}{k}}\label{lma12.20}.
\end{equation} 
We set $r=C\ln(\alpha)$ in the upper bounds of \Ref{lma12.1} and \Ref{lma12.20} and find $C$ such that both bounds are bounded by $1/(2\alpha)$, completing 
the proof of \Ref{lma12s3}.

The same proof can be adapted for \Ref{lma12s4}. With \Ref{lma12.0} in mind, recalling the fact established in the proof of Lemma~\ref{lma105} that $\mathbb{P}(E_{r,\alpha})\le C_1\alpha r^{-\beta}$, we replace the last inequalities of  \Ref{lma12.1}, \Ref{lma12.2} and \Ref{lma12.20} with the corresponding bound of $ \mathbb{P}(E_{r,\alpha})$ to obtain
\begin{eqnarray}
	&&\left|\mathbb{E}\(\bar{W}_\alpha^2-\bar{W}_{\alpha,r}^2\)\right|\le 2(C_0\alpha)^2 \(C_1\alpha r^{-\beta}\)^{\frac{k-2}{k}},\label{lma12.3}\\
	&&\left|\mathbb{E}\bar{W}_\alpha-\mathbb{E}\bar{W}_{\alpha,r}\right|\le 2C_0\alpha \(C_1\alpha r^{-\beta}\)^{\frac{k-1}{k}},\label{lma12.4}
	\\
	&&\left|\(\mathbb{E}\bar{W}_\alpha\)^2-\(\mathbb{E}\bar{W}_{\alpha,r}\)^2\right|\le 4\(C_0\alpha\)^2\(C_1\alpha r^{-\beta}\)^{\frac{k-1}{k}}\label{lma12.5}.
\end{eqnarray}
The claim \Ref{lma12s4} follows by combining \Ref{lma12.3} and \Ref{lma12.5}, extracting $\alpha$ and $r$, and then taking $C$ as the sum of the remaining constant terms.

A line-by-line repetition of the above proof with $\bar{W}_\alpha$ and $\bar{W}_{\alpha,r}$ replaced by $W_\alpha$ and $W_{\alpha,r}$ gives \Ref{lma12s1} and \Ref{lma12s2} respectively. \qed

Next, we apply Lemma~\ref{lma6} and Lemma~\ref{lma7} to establish lower bounds for $\var\(W_{\alpha,r}\)$ and $\var\(\bar{W}_{\alpha,r}\)$.

\begin{lma}\label{lma13}
	\begin{description}
		\item{(a)} (unrestricted case) If the score function $\eta$ satisfies non-singularity \Ref{non-sin}, then $\var\(W_{\alpha,r}\)\ge C\alpha r^{-d}$ for $R_0\le r\le \alpha^{1/d}/6$, where $C,R_0>0$ are independent of $\alpha$.
		\item{(b)} (restricted case) If the score function $\eta$ satisfies non-singularity \Ref{non-sinr}, then $\var\(\bar{W}_{\alpha,r}\)\ge C\alpha r^{-d}$ for $R_0\le r\le \alpha^{1/d}/6$, where $C,R_0>0$ are independent of $\alpha$.
	\end{description}
\end{lma}

\noindent{\it Proof.} For (b), recalling the notations in the paragraph after \Ref{lemma4.4proofa1}, we obtain from the total variance formula that
\begin{align}
	\var\(\bar{W}_{\alpha,r}\)=&\mathbb{E}\(\var\left(\bar{W}_{\alpha,r}\middle| \mathscr{F}_{2,\alpha,r} \right)\)+\var\(\mathbb{E}\left(\bar{W}_{\alpha,r}\middle| \mathscr{F}_{2,\alpha,r} \right)\)\nonumber
	\\\ge& \mathbb{E}\(\var\left(\bar{W}_{\alpha,r}\middle| \mathscr{F}_{2,\alpha,r} \right)\)\nonumber
	\\=&\sum_{i=1}^{m_{\alpha,r}}\mathbb{E}\(\var\left(\sum_{x\in \overline{\Xi}\cap B_i}{\bar{g}}(\Xi^{\Gamma_\alpha,x})\mathbf{1}_{\bar{R}(x,\alpha)\le r }\middle|\Xi_{D_i}\right)\)\nonumber
	\\=&m_{\alpha,r}\mathbb{E}\(\var\left(Y_r\middle|\Xi_{N_0^c}\)\)\label{lma13.1},
\end{align} 
where $m_{\alpha,r}$ is the number of disjoint cubes with length $4r+2r_1$ embedded into $\Gamma_\alpha$. 

Using \Ref{forlma5.9}, there exists an $R_0\ge r_1>0$ such that for all $r>R_0$,
$$\law\left(Y_r\mathbf{1}_{(E_{R_0}')^c}\middle|\Xi_{N_0^c}\right)=\law\left(\tilde{Y}\middle|\Xi_{N_0^c}\right)\ge \tilde{\xi} \mbox{ a.s.},$$
where $\tilde{\xi}$ is an absolutely continuous $\sigma(\Xi_{N_0^c})$ measurable random measure satisfying $\mathbb{P}\(\tilde{\xi}(\mathbb{R})>\frac{p}{2}\)>2p$. Hence, for $r> R_0$, we apply Lemma~\ref{lma7} with $X:=Y_r$, $A:=(E_{R_0}')^c$ and use the fact that $Y_r\mathbf{1}_{(E_{R_0}')^c}=\tilde{Y}$ for all $r\ge R_0$ to obtain 
\begin{equation}
	\mean\var\left(Y_r\middle|\Xi_{N_0^c}\right)\ge \mean\var\left({\tilde{Y}}+\frac{\mathbb{E}\left({\tilde{Y}}\middle|\Xi_{N_0^c}\right)}{\mathbb{P}\left((E_{R_0}')^c\middle|\Xi_{N_0^c}\right)}\mathbf{1}_{E_{R_0}'}\middle|\Xi_{N_0^c}\right)=:b>0.\label{lma13.2}
\end{equation} The proof of claim (b) is completed by combining \Ref{lma13.1} and \Ref{lma13.2} with the observation that
{$R_0\le r\le \alpha^{1/d}/6$}  ensures $m_{\alpha,r}\ge 12^{-d}\alpha r^{-d}$.

The claim (a) can be proved by replacing $\bar{W}_{\alpha,r}$ with $W_{\alpha,r}${; $\bar{g}$ with $g$} throughout the above argument. \qed

Finally, we make use of \cite[Lemma~4.6]{XY15}, Lemma~\ref{lma12} and Lemma~\ref{lma13} to establish Lemma~\ref{lma4}.

\noindent{\it Proof of Lemma~\ref{lma4}.} To begin with, we combine \Ref{lma12s1}, \Ref{lma12s3} and Lemma~\ref{lma13}~(a) to find  an $r:=C_1\ln(\alpha)$ such that 
\begin{eqnarray}
	&&\left|\var\(W_{\alpha}\)-\var\(W_{\alpha,r}\)\right|\le \frac{1}{\alpha},\label{lma4proof1}\\
	&&\left|\var\(\bar{W}_{\alpha}\)-\var\(\bar{W}_{\alpha,r}\)\right|\le \frac{1}{\alpha},\label{lma4proof2}\\
	&&\var\(W_{\alpha,r}\)\ge C_2 \alpha \ln(\alpha)^{-d},\label{lma4proof3}\end{eqnarray}
for positive constants $C_1,C_2$. The inequalities \Ref{lma4proof1} and \Ref{lma4proof3} imply $\var\(W_{\alpha}\)\ge O\(\alpha \ln(\alpha)^{-d}\),$ hence the claim (a) follows from the dichotomy established in \cite[Lemma~4.6]{XY15} saying either $\var\(W_{\alpha}\)=\Omega\(\alpha\)$ or $\var\(W_{\alpha}\)=O\(\alpha^{\frac{d-1}{d}}\)$.

In terms of (b), it suffices to show $\var\(\bar{W}_{\alpha}\)-\var\(W_{\alpha}\)=o(\alpha)$ if we take $\bar{\eta}$ as the score function in the unrestricted case. To this end, noting that \Ref{lma4proof1} and \Ref{lma4proof2}, it remains to show $\var\(W_{\alpha,r}\)-\var\(\bar{W}_{\alpha,r}\)=o(\alpha)$. However, by the Cauchy-Schwarz inequality, we have 
\begin{align*}
		&\left|\var\(W_{\alpha,r}\)-\var\(\bar{W}_{\alpha,r}\)\right|\\
		=&\left|\var\(W_{\alpha,r}-\bar{W}_{\alpha,r}\)-2\cov\(W_{\alpha,r}-\bar{W}_{\alpha,r},W_{\alpha,r}\)\right|
		\\ \le& \var\(W_{\alpha,r}-\bar{W}_{\alpha,r}\)+2\sqrt{\var\(W_{\alpha,r}-\bar{W}_{\alpha,r}\)\var\(W_{\alpha,r}\)},
\end{align*}
and it follows from $\var\(W_{\alpha}\)=\Omega\(\alpha\)$ and \Ref{lma4proof1} that $\var\(W_{\alpha,r}\)=\Omega\(\alpha\)$, hence the proof is reduced to showing $\var\(W_{\alpha,r}-\bar{W}_{\alpha,r}\)=o(\alpha)$. 

Since ${g_\alpha\(x,\Xi\)}\mathbf{1}_{\bar{R}(x,\alpha)<r}={\bar{g}}\(\Xi^{x}\)\mathbf{1}_{R(x)<r}$ if $d(x, \partial\Gamma_{\alpha})>r$, we have $W_{\alpha,r}-\bar{W}_{\alpha,r}=W_{1,\alpha,r}-W_{2,\alpha,r}$ where $$W_{1,\alpha,r}:=\sum_{x\in \overline{\Xi}_{B\(\partial\Gamma_{\alpha},r\)\cap\Gamma_\alpha}}{\bar{g}}\(\Xi^{x}\)\mathbf{1}_{R(x)<r},~W_{2,\alpha,r}:=\sum_{x\in \bar{\Xi}_{B\(\partial\Gamma_{\alpha},r\)\cap\Gamma_\alpha}}{g_\alpha}\(x,\Xi\)\mathbf{1}_{\bar{R}(x,\alpha)<r}.$$ As the summands of $W_{1,\alpha,r}$ and $W_{2,\alpha,r}$ are in the moat within distance $r$ from the boundary of $\Gamma_\alpha$, both $\var\(W_{1,\alpha,r_{\alpha}'}\)$ and $\var\(W_{2,\alpha,r_{\alpha}'}\)$ are of order $o(\alpha)$, as detailed below. In fact, it follows from \Ref{palm4} that 
{\begin{align}
		&\mathbb{E}\( g_\alpha\(x,\Xi\)\mathbf{1}_{\bar{R}(x,\alpha)<r}\overline{\Xi}(dx)\)\nonumber\\
		=&\mathbb{E}\(g_\alpha\(x,\Xi+\delta_{(x,M_x)}\)\mathbf{1}_{\bar{R}(x,M_x,\alpha,\Xi+\delta_{(x,M_x)})<r}\)\lambda dx\nonumber
		\\=&:P_{x,\alpha,r} dx,\label{lma4.7}
\end{align} }
if we set
\begin{equation}\overline{\Xi}_\alpha^\ast(dx):={g_\alpha\(x,\Xi\)}\mathbf{1}_{\bar{R}(x,\alpha)<r}\overline{\Xi}(dx)-P_{x,\alpha,r} dx,\label{lma4.7p}\end{equation}
then $\mean\(\overline{\Xi}_\alpha^\ast(dx)\overline{\Xi}_\alpha^\ast(dy)\){=\mean\(\overline{\Xi}_\alpha^\ast(dx)\)\mean\(\overline{\Xi}_\alpha^\ast(dy)\)}=0$ if $d(x,y)>2r$. Therefore, 
\begin{align}
	&\var\({W_{2,\alpha,r}}\)\nonumber\\
	=&\int_{x,y\in B\(\partial\Gamma_{\alpha},r\)\cap\Gamma_\alpha}\mean\(\overline{\Xi}_\alpha^\ast(dx)\overline{\Xi}_\alpha^\ast(dy)\)\nonumber
	\\=&\int_{x,y\in B\(\partial\Gamma_{\alpha},r\)\cap\Gamma_\alpha,d(x,y)\le 2r}\mean\(\overline{\Xi}_\alpha^\ast(dx)\overline{\Xi}_\alpha^\ast(dy)\)\nonumber\\
	=&\int_{x,y\in B\(\partial\Gamma_{\alpha},r\)\cap\Gamma_\alpha,d(x,y)\le 2r}
	\left\{\mean\[{g_\alpha\(x,\Xi\)}\mathbf{1}_{\bar{R}(x,\alpha)<r}{g_\alpha\(y,\Xi\)}\mathbf{1}_{\bar{R}(y,\alpha)<r}\overline{\Xi}(dx)\overline{\Xi}(dy)\]-P_{x,\alpha,r}P_{y,\alpha,r} dxdy\right\}.\label{lma4.5}
\end{align} 
Recalling the second order Palm distribution in \Ref{palm5}, we can use the moment condition~\Ref{thm2.1r} together with H\"older's inequality to obtain 
\begin{align}
	&\mean\[\left|{g_\alpha\(x,\Xi\)}\right|\mathbf{1}_{\bar{R}(x,\alpha)<r}\left|{g_\alpha\(y,\Xi\)}\right|\mathbf{1}_{\bar{R}(y,\alpha)<r}\overline{\Xi}(dx)\overline{\Xi}(dy)\]\le C^2(\lambda^2dxdy+\lambda dx),\label{lma4.5ad1}\\
	&\left|P_{x,\alpha,r}\right|\left|P_{y,\alpha,r}\right| dxdy\le C^2 \lambda^2dxdy,\label{lma4.5ad2}
\end{align} 
where $C\ge 1$. Combining these estimates with \Ref{lma4.5} gives
$$\var\({W_{2,\alpha,r}}\)=O\(\alpha^{\frac{d-1}{d}}r^{d+1}\)=o(\alpha).$$
The proof of $\var\({W_{1,\alpha,r}}\)=o(\alpha)$ is similar except we replace \Ref{thm2.1r} with \Ref{thm2.1}. {Consequently, $$\var\(W_{\alpha,r}-\bar{W}_{\alpha,r}\)=\var\(W_{1,\alpha,r}-W_{2,\alpha,r}\)\le 2\(\var\({W_{1,\alpha,r}}\)+\var\({W_{2,\alpha,r}}\)\)=o(\alpha)$$ and the statement follows.}\qed

As the lower bounds in Lemma~\ref{lma5} are very conservative, their proofs are less demanding, as demonstrated below.

\noindent{\it Proof of Lemma~\ref{lma5}.}  \setcounter{con}{1} 
We start with (b). The bound \Ref{lma12s4} ensures 
\begin{equation}\left|\var\(\bar{W}_{\alpha}\)-\var\(\bar{W}_{\alpha,r}\)\right|\le C_{\thecon}\qcon{}\(\alpha^{\frac{3k_0-2}{k_0}}r^{-\beta\frac{k_0-2}{k_0}}\)\vee \(\alpha^{\frac{3k_0-1}{k_0}}r^{-\beta\frac{k_0-1}{k_0}}\)\label{lma5pr01}\end{equation}
for all $r\le \alpha^{\frac{1}{d}}$ and $k'>k_0>k\ge3$. On the other hand, Lemma~\ref{lma13}~(b) says  
$$\var\(\bar{W}_{\alpha,r}\)\ge C_{\thecon}\qcon{}\alpha r^{-d}$$
for $0<R_0\le r\le \alpha^{1/d}/6$. Let $r_\alpha:=\alpha^{\frac{2k-2}{k\beta-2\beta-dk}}$, the assumption $\beta>(3k-2)d/(k-2)$ ensures that $r_\alpha<\alpha^{1/d}/6$ for large $\alpha$ and $k_0>k$ guarantees $\left|\var\(\bar{W}_{\alpha}\)-\var\(\bar{W}_{\alpha,r_\alpha}\)\right|\ll \var\(\bar{W}_{\alpha,r_\alpha}\)$ for large $\alpha$, hence $\var\(\bar{W}_{\alpha}\)\ge C_{\thecon}\qcon{}\alpha r_\alpha^{-d}=O\(\alpha^{\frac{k\beta-2\beta-3dk+2d}{k\beta-2\beta-dk}}\)$, completing the proof. 

For the proof of (a), we can proceed to replace $\bar{W}_\alpha$ with $W_\alpha$ and $\bar{W}_{\alpha,r}$ with $W_{\alpha,r}$ as in the proof of (b). \qed 

The proof of Lemma~\ref{lma5} enables us to get slightly better bounds for $\var\(W_{\alpha,r}\)$ and $\var\(\bar{W}_{\alpha,r}\)$.

\begin{lma}\label{remark3} 
	\begin{description}
		\item{(a)} (unrestricted case) If the score function $\eta$ satisfies the conditions of Lemma~\ref{lma5}~(a), then  
		$\var\(W_{\alpha,r}\)\ge C\(\alpha r^{-d}\)\vee\(\alpha^{\frac{k\beta-2\beta-3dk+2d}{k\beta-2\beta-dk}}\)$ for $R_0\le r\le \alpha^{1/d}/6$, where $C,R_0>0$ are independent of $\alpha$.
		\item{(b)} (restricted case) If the score function $\eta$ satisfies the conditions of Lemma~\ref{lma5}~(b), then  
		$\var\(\bar{W}_{\alpha,r}\)\ge C\(\alpha r^{-d}\)\vee\(\alpha^{\frac{k\beta-2\beta-3dk+2d}{k\beta-2\beta-dk}}\)$ for $R_0\le r\le \alpha^{1/d}/6$, where $C,R_0>0$ are independent of $\alpha$.
	\end{description}
\end{lma}

\noindent{\it Proof.} \setcounter{con}{1} 
We prove (b) only as the proof of (a) is similar. We observe that if $r=r_\alpha=\alpha^{\frac{2k-2}{k\beta-2\beta-dk}}$, then $\alpha r_\alpha^{-d}=\alpha^{\frac{k\beta-2\beta-3dk+2d}{k\beta-2\beta-dk}}$, hence for $r<r_\alpha$, the claim follows from Lemma~\ref{lma13}~(b). For $r>{O(r_\alpha)}$, \Ref{lma5pr01} ensures 
$$\left|\var\(\bar{W}_{\alpha}\)-\var\(\bar{W}_{\alpha,r}\)\right|\le C_{\thecon}\(\alpha^{\frac{3k_0-2}{k_0}}r_\alpha^{-\beta\frac{k_0-2}{k_0}}\)\vee \(\alpha^{\frac{3k_0-1}{k_0}}r_\alpha^{-\beta\frac{k_0-1}{k_0}}\)\ll\var\(\bar{W}_{\alpha}\)$$ for large $\alpha$, hence $\var\(\bar{W}_{\alpha,r}\)=O\(\var\(\bar{W}_{\alpha}\)\)=O\(\alpha^{\frac{k\beta-2\beta-3dk+2d}{k\beta-2\beta-dk}}\),$ as claimed. \qed

\noindent{\it Proof of Theorem~\ref{thm2a}.}  \setcounter{con}{1}
Let $\sigma_{\alpha}^2:=\var\(\bar{W}_{\alpha}\)$, $\sigma_{\alpha,r}^2:=\var\(\bar{W}_{\alpha,r}\)$ and
$\bar{Z}_{\alpha,r}\sim N\(\mean\bar{W}_{\alpha,r}, \sigma_{\alpha,r}^2\)$, then it follows from the triangle inequality that
\begin{equation}d_{TV}(\bar{W}_\alpha,\bar{Z}_\alpha)\le d_{TV}\(\bar{W}_\alpha,\bar{W}_{\alpha,r}\)+d_{TV}\(\bar{Z}_\alpha,\bar{Z}_{\alpha,r}\)+d_{TV}\(\bar{W}_{\alpha,r},\bar{Z}_{\alpha,r}\).\label{thm2a01}\end{equation}  
We take $R_0$ as the maximum of the $R_0$'s of Lemma~\ref{lma3}~(b), Corollary~\ref{cor1}~(b) and Lemma~\ref{lma13}~(b).
We start with exponentially stabilizing case (ii). 

(ii) The first term of \Ref{thm2a01} can be bounded using Lemma~\ref{lma105}~(b), giving
\begin{align}
	&d_{TV}\(\bar{W}_\alpha,\bar{W}_{\alpha,r}\)\le C_{\thecon}\qcon{} \alpha e^{-C_{\thecon}\qcon{}r}\le \frac1\alpha,\label{thm2a02}
\end{align}
for $r>C_{\thecon}\ln(\alpha).$ \setcounter{cproofa}{\thecon} \qcon{}

We can establish an upper bound for the second term $d_{TV}\(\bar{Z}_\alpha,\bar{Z}_{\alpha,r}\)$ of \Ref{thm2a01} using Lemma~\ref{lma10}. To this end, \Ref{lma12s3} gives
\begin{equation}
	\left|\sigma_\alpha^2-\sigma_{\alpha,r}^2\right|\le \frac1\alpha,
	\label{thm2a04}
\end{equation}
which, together with Lemma~\ref{lma4}~(b), implies
\begin{equation}
	\sigma_{\alpha,r}^2=\Omega(\alpha), \ \ \ \ \ \sigma_{\alpha}^2=\Omega(\alpha),
	\label{thm2a03}
\end{equation}
for $r>C_{\thecon}\ln(\alpha).$\setcounter{cproofc}{\thecon}\qcon{}
We combine \Ref{thm2a03} and \Ref{lma12.2} to obtain
\begin{equation}\frac{\left|\mathbb{E}\(\bar{Z}_\alpha\)-\mathbb{E}\(\bar{Z}_{\alpha,r}\)\right|}{\max(\sigma_\alpha,\sigma_{\alpha,r})}=\frac{\left|\mathbb{E}\(\bar{W}_{\alpha}\)-\mathbb{E}\(\bar{W}_{\alpha,r}\)\right|}{\max(\sigma_\alpha,\sigma_{\alpha,r})}\le O\(\alpha^{-2}\),
	\label{thm2a14}
\end{equation}
for $r>C_{\thecon} \qcon{}\ln(\alpha).$
Therefore, it follows from \Ref{thm2a04}, \Ref{thm2a03}, \Ref{thm2a14} and Lemma~\ref{lma10} that
\begin{align}
	d_{TV}(\bar{Z}_\alpha,\bar{Z}_{\alpha,r})
	\le&\sqrt{\frac{2}{\pi}}\(\frac{\left|\mathbb{E}\(\bar{Z}_\alpha\)-\mathbb{E}\(\bar{Z}_{\alpha,r}\)\right|}{\max(\sigma_\alpha,\sigma_{\alpha,r})}+\frac{\left|\var\(\bar{W}_{\alpha}\)-\var\(\bar{W}_{\alpha,r}\)\right|}{\min\(\var\(\bar{W}_{\alpha}\),\var\(\bar{W}_{\alpha,r}\)\)}\)\label{thm2a18}
	\\\le&O(\alpha^{-2})\label{thm2a05}
\end{align} for $r>C_{\thecon} \ln(\alpha).$\setcounter{cproofb}{\thecon}

For the last term of \Ref{thm2a01}, as a linear transformation does not change the total variation distance,  we can rewrite it as 
$$d_{TV}\(\bar{W}_{\alpha,r},\bar{Z}_{\alpha,r}\)=d_{TV}\(V_{\alpha,r},Z\),$$
where $V_{\alpha,r}:=\(\bar{W}_{\alpha,r}-\mean \bar{W}_{\alpha,r}\)/\sigma_{\alpha,r}$ and $Z\sim N(0,1)$. We now appeal to  
Stein's method to tackle the problem. Briefly speaking, Stein's method for normal approximation hinges on a Stein equation (see  \cite[p.~15]{CGS11})
\begin{equation}f'(w)-wf(w)=h(w)-Nh,\label{steineq1}\end{equation}
where $Nh:=\mean h(Z)$. The solution of \Ref{steineq1} satisfies (see \cite[p.~16]{CGS11})
$$\|f_h'\|:=\sup_w\left|f_h'(w)\right|\le 2\|h(\cdot)-Nh\|.$$ 
Hence, for $h=\bone_A$ with $A\in\scrB(\real)$, the solution $f_h=:f_A$ satisfies 
\begin{equation}
	\|f_h'\|\le 2.\label{steineq2}
\end{equation} 
The Stein equation \Ref{steineq1} enables us to bound $d_{TV}\(V_{\alpha,r},Z\)$ through a functional form of $V_{\alpha,r}$ only, giving
\begin{equation}\label{Stein}d_{TV}\(V_{\alpha,r},Z\)\le {\sup_{\{f;\ \|f'\|\le 2\}}}\mathbb{E}\[f'\(V_{\alpha,r}\)-V_{\alpha,r}f\(V_{\alpha,r}\)\].\end{equation} 
Recalling \Ref{lma4.7} and \Ref{lma4.7p}, we can represent $V_{\alpha,r}$ through
$V(dx):=\frac{1}{\sigma_{\alpha,r}}
\overline{\Xi}_\alpha^\ast(dx),$ giving $V_{\alpha,r}=\int_{\Gamma_\alpha}V(dx)$. 
Let $N_{x,\alpha, r}'=B(x,2r)\cap\Gamma_\alpha$ and $N_{x,\alpha, r}''=B(x,4r)\cap\Gamma_\alpha$, we have $N_{x,\alpha, r}'\subset B(x,2r)$ and $N_{x,\alpha, r}''\subset B(x,4r)$, so the volumes of $N_{x,\alpha, r}'$ and $N_{x,\alpha, r}''$
are bounded by $O(r^d)$. Define $S_{x,\alpha,r}'=\int_{N_{x,\alpha,r}'}V(dy)$ and $S_{x,\alpha,r}''=\int_{N_{x,\alpha,r}''}V(dy)$. Since $V(dx)$ is independent of $V(dy)$ if $|x-y|>2r$, $V(dx)$ is independent of $V_{\alpha,r}-S_{x,\alpha,r}'$
and $S_{x,\alpha,r}'V(dx)$ is independent of $V_{\alpha,r}-S_{x,\alpha,r}''$, $1=\var\(V_{\alpha,r}\)=\mean\int_{\Gamma_\alpha} S_{x,\alpha,r}'V(dx)$ and
\begin{align}
	&\mathbb{E}\[f'\(V_{\alpha,r}\)-V_{\alpha,r}f\(V_{\alpha,r}\)\]\nonumber
	\\=&\mathbb{E}f'\(V_{\alpha,r}\) -\E\int_{\Gamma_\alpha}\left(f(V_{\alpha,r})-f\left(V_{\alpha,r}-S_{x,\alpha,r}'\right)\right)V(dx)\nonumber
	\\=&\mathbb{E}f'\(V_{\alpha,r}\) - \E\int_{\Gamma_\alpha}\int_{0}^1f'\left(V_{\alpha,r}-uS_{x,\alpha,r}'\right)S_{x,\alpha,r}'du V(dx)\nonumber
	\\=&\mathbb{E}\int_{\Gamma_\alpha}\mathbb{E}\[f'\(V_{\alpha,r}\)-f'\left(V_{\alpha,r}-S_{x,\alpha,r}''\right)\]S_{x,\alpha, r}'V(dx)\nonumber
	\\&-\mathbb{E}\int_{\Gamma_\alpha}\int_{0}^1\(f'\left(V_{\alpha,r}-uS_{x,\alpha,r}'\right)-f'\left(V_{\alpha,r}-S_{x,\alpha,r}''\right)\)S_{x,\alpha, r}'du V(dx).\label{thm2.13}
\end{align} 
By the definition of the total variation distance, we have 
\begin{equation}d_{TV}\(V_{\alpha,r}, V_{\alpha,r}+\gamma\)=d_{TV}\(\bar{W}_{\alpha,r}, \bar{W}_{\alpha,r}+\sigma_{\alpha,r}\gamma\)\label{thm2a06}\end{equation} for any $\gamma\in \mathbb{R}$. Using Corollary~\ref{cor1}~(b) with $N_{\alpha,r}^{(1)}=N_{\alpha,r}^{(3)}:=\emptyset$ and $N_{\alpha,r}^{(2)}:=B\(N_{x,\alpha,r}'',r\)$, {for $r\le  {\qcon{}C_{\thecon}} \alpha^{\frac{1}{d}}$, }we have 
\begin{equation*}
	d_{TV}\left(\bar{W}_{\alpha,r},\bar{W}_{\alpha,r}+h_{\alpha,r}\(\Xi_{N_{\alpha, r}^{(2)}}\)\right)\le \mathbb{E}\left(\left|h_{\alpha,r}\(\Xi_{N_{\alpha, r}^{(2)}}\)\right|\vee 1\right)O\left(\alpha^{-\frac{1}{2}}r^{\frac{d}{2}}\right),
\end{equation*}
which, together with \Ref{thm2a06}, implies
\begin{align}
	&\left|\mathbb{E}\[f'\(V_{\alpha,r}\)-f'\left(V_{\alpha,r}-S_{x,\alpha,r}''\right)\]\right| \le2\|f'\|O\(\alpha^{-\frac{1}{2}}r^{\frac{d}{2}}\)\mathbb{E}\[\left|\sigma_{\alpha,r}S_{x,\alpha,r}''\right|\vee 1\]\label{thm2a07}.
\end{align}
Recalling \Ref{lma4.7p}, we have
$$\sigma_{\alpha,r}S_{x,\alpha,r}''=\int_{N_{x,\alpha, r}''} \overline{\Xi}_\alpha^\ast(dy).$$
Using the first order Palm distribution \Ref{palm4}, the third order Palm distribution \Ref{palm6} and the moment condition~\Ref{thm2.1r}, we obtain
\begin{align*}
	&\mean\[\left|{g_\alpha\(z,\Xi\)}\right|\mathbf{1}_{\bar{R}(z,\alpha)<r}\overline{\Xi}(dz)\left|{g_\alpha\(y,\Xi\)}\right|\mathbf{1}_{\bar{R}(y,\alpha)<r}\overline{\Xi}(dy)\left|{g_\alpha\(x,\Xi\)}\right|\mathbf{1}_{\bar{R}(x,\alpha)<r}\overline{\Xi}(dx)\]\\
	&\ \ \ \le \qcon{} C_{\thecon}\qcon{} \(\lambda^3 dzdydx+\lambda^2dzdx+\lambda^2dydx+\lambda dx\),\\
	&\mean\[\left|{g_\alpha\(y,\Xi\)}\right|\mathbf{1}_{\bar{R}(y,\alpha)<r}\overline{\Xi}(dy)\]\le C_{\thecon}\qcon{} \lambda dy,
	\\
	&\left|P_{y,\alpha,r}\right|\le C_{\thecon}\qcon{}  \lambda,
\end{align*} 
which, together with \Ref{lma4.5ad1} and \Ref{lma4.5ad2}, yield
\begin{align}
	&\mean\left|\overline{\Xi}_\alpha^\ast(dy)\right|\le C_{\thecon}\qcon{}  \lambda dy,\label{thm2a10}\\
	&\mean\left|\overline{\Xi}_\alpha^\ast(dy)\overline{\Xi}_\alpha^\ast(dx)\right|\le C_{\thecon}\qcon{}  \(\lambda^2 dydx+\lambda dx\),\label{thm2a11}\\
	&\mean\left|\overline{\Xi}_\alpha^\ast(dz)\overline{\Xi}_\alpha^\ast(dy)\overline{\Xi}_\alpha^\ast(dx)\right|\le C_{\thecon} \(\lambda^3 dzdydx+\lambda^2dzdx+\lambda^2dydx+\lambda dx\),\label{thm2a11-0}\\
	&\mean\int_{N_{x,\alpha, r}'} \left|\overline{\Xi}_\alpha^\ast(dy)\right|\le \mean\int_{N_{x,\alpha, r}''} \left|\overline{\Xi}_\alpha^\ast(dy)\right|\le O\(r^d\),\nonumber\\
	&\mathbb{E}\[\left|\sigma_{\alpha,r}S_{x,\alpha,r}''\right|\vee 1\]
	\le  1+\mean\int_{N_{x,\alpha, r}''} \left|\overline{\Xi}_\alpha^\ast(dy)\right|\le O\(r^d\).\label{thm2a13}
\end{align}
Combining \Ref{thm2a07}, \Ref{thm2a13} and \Ref{steineq2}, we have
\begin{align*}
	&\left|\mathbb{E}\[f'\(V_{\alpha,r}\)-f'\left(V_{\alpha,r}-S_{x,\alpha,r}''\right)\]\right| \le O\(\alpha^{-\frac{1}{2}}r^{\frac{3d}{2}}\),
\end{align*}
hence the first term of \Ref{thm2.13} can be bounded as
\begin{align}
	&\left|\mathbb{E}\int_{\Gamma_\alpha}\mathbb{E}\[f'\(V_{\alpha,r}\)-f'\left(V_{\alpha,r}-S_{x,\alpha,r}'' \right)\]S_{x,\alpha, r}' V(dx)\right|\nonumber
	\\ \le &O\(\alpha^{-\frac{1}{2}}r^{\frac{3d}{2}}\)\sigma_{\alpha,r}^{-2}\mathbb{E}\int_{\Gamma_\alpha}\int_{N_{x,\alpha, r}'}\left|\overline{\Xi}_\alpha^\ast(dy)\right|\left|\overline{\Xi}_\alpha^\ast(dx)\right|\nonumber
	\\ \le &O\(\alpha^{-\frac{1}{2}}r^{\frac{3d}{2}}\)\sigma_{\alpha,r}^{-2}\int_{\Gamma_\alpha}\(\int_{N_{x,\alpha, r}'}  \lambda dy+1\)\lambda dx = {O\(\sigma_{\alpha,r}^{-2}\alpha^{\frac{1}{2}}r^{\frac{5d}{2}}\)},\label{thm2a16}
\end{align} 
where the last inequality is from \Ref{thm2a11}. 

For the second term of \Ref{thm2.13}, we have from Corollary~\ref{cor1} with $N_{\alpha,r}^{(1)}:=B\(N_{x,\alpha,r}',r\)$, $N_{\alpha,r}^{(2)}:=B\(N_{x,\alpha,r}'',r\)$, $N_{\alpha,r}^{(3)}:=N_{x,\alpha,r}''$, for $r\le  {\qcon{}C_{\thecon}} \alpha^{\frac{1}{d}}$, we have
\begin{align}
	&\left|\mean\[\left.\int_{0}^1\(f'\left(V_{\alpha,r}-uS_{x,\alpha,r}'\right)-f'\left(V_{\alpha,r}-S_{x,\alpha,r}''\right)\)du\right|\Xi_{N_{x,\alpha,r}'}\]\right|\nonumber\\
	\le&2\int_{0}^{1}\mathbb{E}d_{TV}\left(\left.V_{\alpha,r}-uS_{x,\alpha,r}',V_{\alpha,r}-S_{x,\alpha,r}''\right|\Xi_{N_{x,\alpha,r}'} \right) du\nonumber
	\\\le&O\(\alpha^{-\frac{1}{2}}r^{\frac{d}{2}}\)\mean\(\left.\int_{N_{x,\alpha, r}''}\left|\overline{\Xi}_\alpha^\ast(dz)\right|+ 1\right|\Xi_{N_{x,\alpha,r}'}\),\nonumber
\end{align} 
hence
\begin{align}
	&\left|\mathbb{E}\int_{\Gamma_\alpha}\int_{0}^1\(f'\left(V_{\alpha,r}-uS_{x,\alpha,r}'\right)-f'\left(V_{\alpha,r}-S_{x,\alpha,r}''\right)\)S_{x,\alpha, r}'du V(dx)\right|\nonumber\\
	&=\left|\mathbb{E}\int_{\Gamma_\alpha}\mean\[\left.\int_{0}^1\(f'\left(V_{\alpha,r}-uS_{x,\alpha,r}'\right)-f'\left(V_{\alpha,r}-S_{x,\alpha,r}''\right)\)du\right|\Xi_{N_{x,\alpha,r}'}\]S_{x,\alpha, r}' V(dx)\right|\nonumber\\
	&\le O\(\alpha^{-\frac{1}{2}}r^{\frac{d}{2}}\)\mathbb{E}\int_{\Gamma_\alpha}\mean\(\left.\int_{N_{x,\alpha, r}''}\left|\overline{\Xi}_\alpha^\ast(dz)\right|+ 1\right|\Xi_{N_{x,\alpha,r}'}\)\left|S_{x,\alpha, r}'\right| |V(dx)|\nonumber\\
	&\le O\(\alpha^{-\frac{1}{2}}r^{\frac{d}{2}}\)\sigma_{\alpha,r}^{-2}\mean\int_{\Gamma_\alpha}\[\int_{N_{x,\alpha, r}''}\int_{N_{x,\alpha, r}'}\left|\overline{\Xi}_\alpha^\ast(dz)\overline{\Xi}_\alpha^\ast(dy)\right|+\int_{N_{x,\alpha, r}'}\left|\overline{\Xi}_\alpha^\ast(dy)\right|\]\left|\overline{\Xi}_\alpha^\ast(dx)\right|\nonumber\\
	&\le O\(\alpha^{-\frac{1}{2}}r^{\frac{d}{2}}\)\sigma_{\alpha,r}^{-2}\int_{\Gamma_\alpha}\(\int_{N_{x,\alpha, r}''}\int_{N_{x,\alpha, r}'}\lambda^2dzdy+\int_{N_{x,\alpha, r}''}\lambda dz+\int_{N_{x,\alpha, r}'}\lambda dy+1\)\lambda dx\nonumber\\
	&\le O\(\alpha^{-\frac{1}{2}}r^{\frac{d}{2}}\)\sigma_{\alpha,r}^{-2}O\(\alpha r^{2d}\)\nonumber\\
	&= O\(\sigma_{\alpha,r}^{-2}\alpha^{\frac{1}{2}}r^{\frac{5d}{2}}\),\label{thm2a17} \end{align}
where the second last inequality follows from \Ref{thm2a10}, \Ref{thm2a11}, \Ref{thm2a11-0}, and the last inequality is due to the fact that the volumes of $N_{x,\alpha, r}'$ and $N_{x,\alpha, r}''$ are 
bounded by $O\(r^d\)$. Recalling \Ref{Stein} and \Ref{thm2.13}, we add up the bounds of \Ref{thm2a16} and \Ref{thm2a17} to obtain
\begin{equation}
	d_{TV}\(\bar{W}_{\alpha,r},\bar{Z}_{\alpha,r}\)=d_{TV}\(V_{\alpha,r},Z\)\le O\(\sigma_{\alpha,r}^{-2}\alpha^{\frac{1}{2}}r^{\frac{5d}{2}}\).\label{thm2.25}
\end{equation} 
The proof of (ii) is completed by using \Ref{thm2a01}, taking $r={\max(C_{\thecproofa},C_{\thecproofc},C_{\thecproofb})}\ln(\alpha)$ for large $\alpha$, collecting the bounds in \Ref{thm2a02}, \Ref{thm2a05}, \Ref{thm2.25} and replacing $\sigma_{\alpha,r}^{2}=\Omega(\alpha)$, as shown in \Ref{thm2a03}.

(i) There exists an $r_1>0$ such that $\bar{W}_{\alpha,r_1}= \bar{W}_{\alpha}$ $a.s.$ for all $\alpha$, which implies $\mean \bar{W}_{\alpha,r_1}=\mean \bar{W}_{\alpha}$, $\var\(\bar{W}_{\alpha,r_1}\)=\var\(\bar{W}_{\alpha}\)$, hence $d_{TV}(\bar{W}_{\alpha},\bar{Z}_\alpha)= d_{TV}(\bar{W}_{\alpha,r_1},\bar{Z}_{\alpha,r_1})$. On the other hand, range-bound implies exponential stabilization, with $r_1$ in place of $r$, \Ref{thm2a03} and \Ref{thm2.25} still hold. However, $r_1$ is a constant independent of $\alpha$, the conclusion follows. 

(iii) We take $r=r_\alpha:= R_0\vee \alpha^{\frac{5k-4}{5dk+2\beta k -4\beta}}$. Lemma~\ref{lma105}~(b) gives
\begin{equation}d_{TV}\(\bar{W}_{\alpha},\bar{W}_{\alpha,r}\)\le O\(\alpha r^{-\beta}\)<O\(\alpha^{-\frac{\beta(k-2)[\beta(k-2)-d(15k-14)]}{(k\beta-2\beta-dk)(5dk+2\beta k-4\beta)}}\) .\label{thm2.26}
\end{equation}
Next, applying Lemma~\ref{lma5}~(b) and Lemma~\ref{remark3}~(b), we have
\begin{equation}\var\(\bar{W}_{\alpha,r}\)\wedge \var\(\bar{W}_{\alpha}\)\ge O\(\alpha^{\frac{k\beta-2\beta-3dk+2d}{k\beta-2\beta-dk}}\),\label{thm2a19}\end{equation}
which, together with \Ref{thm2a18}, \Ref{lma12.4} and \Ref{lma12s4}, yields 
\begin{equation}d_{TV}(\bar{Z}_\alpha,\bar{Z}_{\alpha,r})\le O\(\frac{\alpha^{\frac{2k-1}k}r^{-\beta\frac{k-1}k}}{\alpha^{\frac12\frac{k\beta-2\beta-3dk+2d}{k\beta-2\beta-dk}}}\vee\frac{\alpha^{\frac{3k-2}k}r^{-\beta\frac{k-2}k}}{\alpha^{\frac{k\beta-2\beta-3dk+2d}{k\beta-2\beta-dk}}}\vee\frac{\alpha^{\frac{3k-1}k}r^{-\beta\frac{k-1}k}}{\alpha^{\frac{k\beta-2\beta-3dk+2d}{k\beta-2\beta-dk}}}\)\label{thm2a20}
\end{equation}
for $R_0< r< C_{\qcon{}\thecon}\alpha^{\frac{1}{d}}$. Recalling that $\beta>\frac{(15k-14)d}{k-2}$, the dominating term of \Ref{thm2a20} is $\frac{\alpha^{\frac{3k-2}k}r^{-\beta\frac{k-2}k}}{\alpha^{\frac{k\beta-2\beta-3dk+2d}{k\beta-2\beta-dk}}}$, giving
\begin{equation}d_{TV}(\bar{Z}_\alpha,\bar{Z}_{\alpha,r})\le O\(\frac{\alpha^{\frac{3k-2}k}r^{-\beta\frac{k-2}k}}{\alpha^{\frac{k\beta-2\beta-3dk+2d}{k\beta-2\beta-dk}}}\)=O\(\alpha^{-\frac{\beta(k-2)[\beta(k-2)-d(15k-14)]}{(k\beta-2\beta-dk)(5dk+2\beta k-4\beta)}}\).\label{thm2a21}
\end{equation}
In terms of $d_{TV}(\bar{W}_{\alpha,r},\bar{Z}_{\alpha,r})$, we make use of \Ref{thm2.25} and \Ref{thm2a19} and replace $r$ with $r_\alpha$ to obtain
\begin{equation}\label{thm2.27}d_{TV}(\bar{W}_{\alpha,r},\bar{Z}_{\alpha,r})\le O(\alpha^{\frac{1}{2}} r^{\frac{5d}{2}}) O\(\alpha^{-\frac{k\beta-2\beta-3dk+2d}{k\beta-2\beta-dk}}\)=O\(\alpha^{-\frac{\beta(k-2)[\beta(k-2)-d(15k-14)]}{(k\beta-2\beta-dk)(5dk+2\beta k-4\beta)}}\).
\end{equation} Finally, the proof is completed by combining \Ref{thm2a01}, \Ref{thm2.26}, \Ref{thm2a21} and \Ref{thm2.27}. \qed

\noindent{\it Proof of Theorem~\ref{thm2}.} One can repeat the proof of Theorem~\ref{thm2a} by replacing $\bar{W}_\alpha$, $\bar{W}_{\alpha,r}$, $\bar{Z}_\alpha$, $\bar{Z}_{\alpha,r}$, $g_\alpha(x,\Xi)$ and $\bar{R}(x,\alpha)$ with $W_\alpha$, $W_{\alpha,r}$, $Z_\alpha$, $Z_{\alpha,r}$, $g(\Xi^x)$ and $R(x)$. \qed

\begin{re} {\rm~If we aim to find the order of the total variation distance between $\bar{W}_{\alpha}$ and a normal distribution instead of a normal distribution with the same mean and variance in the polynomially stabilizing case, we can get a better upper bound approximation error with a weaker condition. When $\beta>\frac{5dk-7d+\sqrt{20d^2k^2-60d^2k+49d^2}}{k-2}$, combining \Ref{thm2.25} and the fact that $d_{TV}(\bar{W}_{\alpha}, \bar{W}_{\alpha,r})\le C\alpha\lambda r^{-\beta}$, taking $r_\alpha:=\alpha^{\frac{3\beta k-7dk+4d-6\beta}{(\beta k-dk-2\beta)(5d+2\beta)}}$, we have 
			\begin{align*}
				d_{TV}(\bar{W}_{\alpha},\bar{Z}_{\alpha,r_\alpha})&\le d_{TV}(\bar{W}_{\alpha},\bar{W}_{\alpha,r_\alpha})+d_{TV}(\bar{W}_{\alpha,r_\alpha},\bar{Z}_{\alpha,r_\alpha})
				\\&\le O\(\alpha^{\frac{-\beta^2(k-2)+10\beta dk-14\beta d-5d^2k}{(\beta k-dk-2\beta)(5d+2\beta)}}\).
		\end{align*}}
\end{re}


\def\ac{{Academic Press}~}
\def\aap{{Adv. Appl. Prob.}~}
\def\ap{{Ann. Probab.}~}
\def\anap{{Ann. Appl. Probab.}~}
\def\eljp{{\it Electron.\ J.~Probab.\/}~} 
\def\jap{{J. Appl. Probab.}~}
\def\jws{{John Wiley $\&$ Sons}~}
\def\ny{{New York}~}
\def\ptrf{{Probab. Theory Related Fields}~}
\def\sp{{Springer}~}
\def\spa{{Stochastic Processes and their Applications}~}
\def\sv{{Springer-Verlag}~}
\def\tpa{{Theory Probab. Appl.}~}
\def\zw{{Z. Wahrsch. Verw. Gebiete}~}


\begin{thebibliography}{9}
	\bibitem[Avram and Bertsimas~(1993)]{AB93} Avram, F. and Bertsimas, D.~(1993). On central limit theorems in geometrical probability. \emph{\anap}\textbf{3}, 1033--1046.
	
	\bibitem[Bally and Caramellino~(2016)]{BC16}
	Bally, V. and Caramellino, L.~(2016). Asymptotic development for the CLT in total
variation distance. \emph{Bernoulli}~\textbf{22}, 2442--2485.
	
	\bibitem[Barbour, Holst and Janson(1992)]{BHJ}
	Barbour,  A. D.,  Holst, L. and  Janson, S. (1992). 
	{\em Poisson approximation}. 
	Oxford University Press.
	
	
	\bibitem[Barbour, Luczak and Xia~(2018)]{BLX18} Barbour, A. D.,  Luczak, M. J. and Xia, A. (2018). Multivariate approximation in total variation, I: equilibrium distributions of Markov jump processes. \emph{\ap}\textbf{46}, 1351--1404.
	
	
	\bibitem[Berry~(1941)]{Berry41} Berry, A. C. (1941). The accuracy of the Gaussian approximation to the sum of independent variates. \emph{Trans. Amer. Math. Soc.}~\textbf{49}, 122--136.
	
	
	\bibitem[Cai~(1980)]{C80} Cailliez, F. (1980). Forest volume estimation and yield prediction. \emph{Food 
	and Agriculture Organization of the United Nations}.
	
	\bibitem[\Ceka~(2000)]{Ceka00} \Ceka, V. (2000). Remarks on estimates in the total-variation metric. \emph{Lithuanian Mathematical Journal}~\textbf{40},
	1--13.
	
	
	\bibitem[Chen, Goldstein and Shao~(2011)]{CGS11} Chen, L. H. Y., Goldstein, L. and Shao, Q. M. (2011). {\em Normal approximation by Stein's method}. Springer-Verlag.
	
	\bibitem[Chen, Hwang and Tsai~(2003)]{CHT03}Chen, W. M., Hwang, H. K. and Tsai, T. H. (2003). Efficient maxima-finding algorithms for
random planar samples. \emph{Discrete Mathematics and Theoretical Computer Science}~\textbf{6}, 107--122.

	\bibitem[Chen and Leong~(2010)]{CL10} Chen, L. H. Y. and Leong, Y. K. (2010). From zero-bias to discretized normal approximation. Preprint.
	
	\bibitem[Chen, R\"ollin and Xia~(2020)]{CRX20}  Chen, L. H. Y., R\"ollin, A. and Xia, A.~(2020). Palm theory, random measures and Stein~couplings. \textit{\anap}(to appear). 
	
	\bibitem[Chen and Xia~(2004)]{CX04}
	Chen, L. H. Y. and Xia, A. (2004). 
	Stein's method, Palm theory and Poisson process approximation. \textit{\ap}\textbf{32}, 2545--2569. 

	\bibitem[Daley \& Vere-Jones~(2008)]{Daley08} Daley, D. J. and Vere-Jones, D. (2008). {\it An introduction to the
		theory of point processes.\/} Vol. 2, Springer, New York.
	
	\bibitem[Devroye~(1988)]{Devroye88} Devroye, L. (1988). The expected size of some graphs in computational geometry. \emph{Computers $\&$ Mathematics with Applications}~\textbf{15}, 53--64. 
		
	\bibitem[Diaconis and Freedman~(1987)]{DF87} Diaconis, P. and Freedman, D. (1987). A dozen de Finetti-style results in search of a theory. \emph{Ann.
		Inst. H. Poincar\'e Probab. Statist.}~\textbf{23}, no. 2, suppl., 397--423.
		
	\bibitem[Esseen~(1942)]{Esseen42} Esseen, C. G. (1942). On the Liapounoff limit of error in the theory of probability. \emph{Ark. Mat. Astr. Fys.}~\textbf{28A}, 1--19. 
	
	\bibitem[Fang~(2014)]{Fang14} Fang, X. (2014). Discretized normal approximation by Stein's method. \emph{Bernoulli}~\textbf{20}, 1404--1431.
	
	
	\bibitem[Feller~(1971)]{Feller71} Feller, W. (1971). {\em An introduction to probability theory and its applications}. Vol.~2, John Wiley and Sons.
	
	
	\bibitem[Goldstein and Xia~(2006)]{GX06} Goldstein, L. and Xia, A. (2006). Zero biasing and a discrete central limit theorem. \emph{\ap}\textbf{34}, 1782--1806. 
	
	\bibitem[Halmos~(1974)]{Halmos74} Halmos, P. R. (1974). {\em Measure theory}. Graduate Texts in Mathematics 18, Springer-Verlag.
		   


	\bibitem[Kallenberg~(1983)]{Kallenberg83} Kallenberg, O. (1983). {\em Random measures.\/} Academic Press, London.
	
	\bibitem[Kallenberg~(2017)]{Kallenberg17} Kallenberg, O. (2017). {\em Random measures, theory and applications.} Springer-Verlag.
	
		\bibitem[Khanteimouri et al.~(2017)]{K17}Khanteimouri, P.,  Mohades, A., Abam, M. A., Kazemi, M. R. and Sedighin, S. (2017). Efficiently computing the smallest axis-parallel squares
spanning all colors. \emph{Scientia Iranica D}~\textbf{24}, 1325--1334.
		
	
	\bibitem[Kung~(1975)]{K75} Kung, H. T., Luccio, F. and Preparata, F. P. (1975). On finding the maxima of a set of vectors. \emph{Journal of the ACM}~\textbf{22}, 469--476.
	
	\bibitem[Lachi\`eze-Rey, Schulte and Yukich~(2019)]{LSY19} Lachi\`eze-Rey, R., Schulte,  M. and Yukich, J. E. (2019). 
	Normal approximation for stabilizing functionals. 
	\emph{\anap}\textit{29}, 931--993.
	

	\bibitem[Li et al.~(2015)]{Li15} Li, C., Barclay, H., Hans, H. and Sidders, D. (2015). Estimation of log volumes: A Comparative Study. 
	\emph{Canadian Wood Fibre Centre.} 
	
	\bibitem[Lindvall~(1992)]{Lindvall92} Lindvall, T. (1992). {\em Lectures on the coupling method}. Wiley, New York.
	
	\bibitem[McGivney and Yukich(1999)]{MY99} McGivney, K. and Yukich, J. E. (1999). 
	Asymptotics for Voronoi tessellations on random samples. \textit{Stochastic Process. Appl.} \textbf{83}, 273--288.
	
        \bibitem[Mecke~(1967)]{Mecke63} Mecke, J. (1967). Zum Problem der Zerlegbarkeit station\"{a}rer rekurrenter zuf\"{a}lliger Punktfolgen.
	\emph{Mathematische Nachrichten}~\textbf{35}, 311--321.

	\bibitem[Meckes and  Meckes~(2007)]{MM07} Meckes, E. S. and  Meckes, M. W. (2007). The central limit problem for random vectors with symmetries.
	\emph{Journal of Theoretical Probability}~\textbf{20}, 697--720.
		
	
	
	\bibitem[Penrose and Yukich~(2001)]{PY01} Penrose, M. D. and Yukich, J. E. (2001). 
	Central limit theorems for some graphs in computational geometry.
	\emph{\anap}\textbf{11}, 1005--1041.
	
	\bibitem[Penrose and Yukich~(2005)]{PY05} Penrose, M. D. and Yukich, J. E. (2005). 
	Normal approximation in geometric probability. \emph{Stein's Method and Applications}, Eds. A.~D.~Barbour \& L.~H.~Y.~Chen, World Scientific
	Press, Singapore, pp 37--58. 
	
	\bibitem[R\'{e}nyi~(1962)]{R62}R\'{e}nyi, A. (1962). Th\'{e}orie des \'{e}l\'{e}ments saillants d'une suite d'observations. \emph{Annales scientifiques de l'Universit\'{e} de Clermont. Math\'{e}matiques}~\textbf{8}, 7--13. 
		
	
	\bibitem[R\"ollin~(2005)]{Rollin05} R\"ollin, A. (2005). Approximation of sums of conditionally independent variables by the translated Poisson distribution. \emph{Bernoulli}~\textbf{11}, 1115--1128. 
	
	\bibitem[R\"ollin~(2007)]{Rollin07} R\"ollin, A. (2007). Translated Poisson approximation using exchangeable pair couplings. \emph{\anap}\textbf{17}, 1596--1614. 
	
	\bibitem[R\"ollin~(2008)]{Rollin08} R\"ollin, A. (2008). Symmetric and centered binomial approximation of sums of locally dependent random variables. \emph{Electron. J. Probab.}~\textbf{13}, 756--776. 
	
	\bibitem[Schulte~(2012)]{S12} Schulte, M. (2012). Normal approximation of Poisson functionals in Kolmogorov
	distance. \emph{J. Theoret. Probab.}~\textbf{29}, 96--117.
	
	\bibitem[Schulte~(2016)]{S16} Schulte, M. (2016). A central limit theorem for the Poisson-Voronoi approximation.
	\emph{Adv. Appl. Math.}~\textbf{49}, 285--306.
	
	\bibitem[Toussaint~(1982)]{T82} Toussaint, G.T. (1982). Computational geometric problems in pattern recognition. In {\em Pattern Recognition Theory and Applications.} (Kittler, J., Fu, K. S., Pau, L. F. eds.) Springer, Dordrecht, 73--91.
	
	\bibitem[Wintner~(1938)]{Wintner38} Wintner, A. (1938). {\em Asymptotic distributions and infinite convolutions}. Edwards Brothers, Ann Arbor, MI.
	
	
	\bibitem[Xia and Yukich~(2015)]{XY15} Xia, A. and Yukich, J. (2015). Normal approximation for statistics of Gibbsian input in geometric probability. \emph{\aap}\textbf{47(4)}, 934--972. 
\end{thebibliography}
\end{document}